\documentclass{amsart}
\usepackage{amsmath,amssymb,amscd,latexsym,enumitem,amsthm, hyperref}
\usepackage[all]{xy}

\DeclareMathOperator{\dist}{dist}
\DeclareMathOperator{\her}{Her}
\DeclareMathOperator{\Q}{\mathcal{Q}}

\DeclareMathOperator{\aff}{Aff_b}
\DeclareMathOperator{\F}{\mathcal{F}}
\DeclareMathOperator{\G}{\mathcal{G}}
\DeclareMathOperator{\id}{\textup{id}}
\DeclareMathOperator{\into}{\hookrightarrow}

\DeclareMathOperator{\uhf}{\mathcal{U}}
\DeclareMathOperator{\Ell}{Ell}
\DeclareMathOperator{\homeo}{Homeo}
\DeclareMathOperator{\Lim}{\underrightarrow{\lim}}

\newcounter{number}[section]

\newenvironment{nummer}{\refstepcounter{number}{\noindent\arabic{section}.\arabic{number}}}{}

\newcommand{\bn}{\noindent \begin{nummer} \rm}
\newcommand{\en}{\end{nummer}}

\newenvironment{ntheorem}{\noindent {\sc Theorem:} \it}{}
\newenvironment{nlemma}{\noindent {\sc Lemma:} \it}{}
\newenvironment{nprop}{\noindent {\sc Proposition:} \it}{}
\newenvironment{ndefn}{\noindent {\sc Definition:} \it}{}
\newenvironment{ncor}{\noindent {\sc Corollary:} \it}{}

\newenvironment{nremark}{\noindent {\sc Remark:}}{}
\newenvironment{nremarks}{\noindent {\sc Remarks:}}{}

\newenvironment{nnotation}{\noindent {\sc Notation:} }{}

\newenvironment{nproof}{\noindent {\sc Proof:}}{\mbox{}\hfill 
\rule[-.2ex]{.25em}{1.8ex}}

\begin{document}

\title[UHF-slicing]{{\sc UHF-slicing and classification of nuclear $\mathrm{C}^{*}$-algebras}}

\author{Karen R. Strung}

\author{Wilhelm Winter}
\address{Mathematisches Institut der 
Westf\"alische Wilhelms-Universit\"at\\
Einsteinstra\ss e 62 \\
48149 M\"unster\\
Germany}

\email{wwinter@uni-muenster.de}

\email{karen.strung@uni-muenster.de}

\date{\today}
\subjclass[2000]{}
\thanks{Supported by EPSRC: EP/G014019/1, EP/I019227/1, DFG: SFB 878, GIF: 1137-30.6/2011.}

\setcounter{section}{-1}

\begin{abstract} In this paper we show that certain simple locally recursive subhomogeneous (RSH) $\mathrm{C}^{*}$-algebras are tracially approximately interval algebras after tensoring with the universal UHF algebra. This involves a linear algebraic encoding of the structure of the local RSH algebra allowing us to find a path through the algebra which looks like a discrete version of $[0,1]$ and exhausts most of the algebra. We produce an actual copy of the interval and use properties of $\mathrm{C}^{*}$-algebras tensored with UHF algebras to move the honest interval underneath the discrete version. It follows from our main result that such $\mathrm{C}^{*}$-algebras are classifiable by Elliott invariants. Our theorem requires finitely many tracial states that all induce the same state on the $K_0$-group; in particular we do not require that projections separate tracial states. We apply our results to classify some examples of $\mathrm{C}^{*}$-algebras constructed by Elliott to exhaust the invariant. We also give an alternate way to classify examples of Lin and Matui of $\mathrm{C}^{*}$-algebras of minimal dynamical systems. In this way our result can be viewed as a first step towards removing the requirement that projections separate tracial states in the classification theorem for $\mathrm{C}^*$-algebras of minimal dynamical systems given by Toms and the second named author.
\end{abstract}

\maketitle

\section{Introduction}

\noindent
The aim of Elliott's programme is to classify separable nuclear $\mathrm{C}^{*}$-algebras by their $K$-theory, tracial state spaces, and the natural pairings between these objects. Restricting to the simple case, the programme has met with some highly satisfactory successes: the classification of the approximately finite dimensional (AF) algebras, the approximately circle (A$\mathbb{T}$) algebras and the approximately homogeneous (AH) algebras with slow dimension growth. 

In the case of the AF, A$\mathbb{T}$ and AH algebras, classification was via inductive limit structures  consisting of manageable classes of $\mathrm{C}^{*}$-algebras from which maps at the level of the invariant are known to be liftable. However, in more general classes it is often difficult to obtain a suitable inductive limit structure for a given $\mathrm{C}^{*}$-algebra. A particular example is the class of $\mathrm{C}^{*}$-algebras of minimal dynamical systems of compact metrizable spaces. Though these are all known to be simple separable nuclear unital stably finite $\mathrm{C}^{*}$-algebras which satisfy the universal coefficient theorem (UCT), in the general case no inductive limit structure by tractable building blocks is known unless the system consists of a smooth manifold with a minimal diffeomorphism \cite{LinPhi:mindifflimits}. Even in this setting, the building blocks turn out to be quite complicated and these $\mathrm{C}^{*}$-algebras remain unclassified.

In 2001, Huaxin Lin introduced the notion of tracial approximation of \mbox{$\mathrm{C}^{*}$-algebras} in his paper on tracially approximately finite dimensional (TAF) algebras \cite{Lin:TAF1}. Using his classification result for TAF algebras, one is given a means of obtaining classification results by showing that a given class of $\mathrm{C}^{*}$-algebras are TAF. This requires no construction of an inductive limit structure; it might be thought of as a route towards classification that is, in a sense, more axiomatic. It might also be regarded as the finite counterpart of Kirchberg's and Phillips' celebrated classification of simple purely infinite $\mathrm{C}^{*}$-algebras.

Tracial approximation has since been generalized beyond the class of finite dimensional algebras. If $\mathcal{S}$ is a given class of separable unital $\mathrm{C}^{*}$-algebras, then a $\mathrm{C}^{*}$-algebra that is TA$\mathcal{S}$ may be thought of as a $\mathrm{C}^{*}$-algebra that is locally approximated by $\mathrm{C}^{*}$-algebras in the class $\mathcal{S}$ in trace. 

If a $\mathrm{C}^{*}$-algebra is simple and TAF,  hence classifiable by the results of Lin \cite{Lin:TAF1}, it is known that it must also have real rank zero \cite[Theorem 3.4]{Lin:TAF1}. This implies the presence of many projections. It is no surprise, then, that data from $K$-theory is more easily extracted. In the case of TAI algebras, that is, those which are tracially approximated by interval algebras, this is no longer necessarily the case; we can no longer assume suitably many projections in any such sense. Despite this,  these $\mathrm{C}^{*}$-algebras are also known to be classifiable in the presence of the UCT \cite{Lin:simpleNuclearTR1}.  Nevertheless, showing a particular class of $\mathrm{C}^{*}$-algebras are TAI remains tricky when one has relatively few projections available.

Further complications to the classification problem can be described (and partly resolved) with the aid of the Jiang--Su algebra $\mathcal{Z}$,  a simple separable unital $\mathrm{C}^{*}$-algebra which is, in a sense, invisible to the usual invariant: despite the fact that the Jiang--Su algebra is infinite dimensional, its $K$-theory and tracial state space are identical to those of $\mathbb{C}$. If $A$ is a $\mathrm{C}^{*}$-algebra that falls within the present scope of the classification programme, it turns out that the invariant of $A$ is isomorphic to the invariant of $A \otimes \mathcal{Z}$.
 
This isomorphism holds at the level of the invariant even in pathological examples where it turns out that the $\mathrm{C}^{*}$-algebra $A$ is not  $\mathcal{Z}$-stable (that is, $A \cong A \otimes \mathcal{Z}$) for example, those algebras constructed in \cite{Vil:perforation, Ror:simple, Toms:classproblem, GioKerr:Subshifts}). Thus one cannot assume a priori that a particular $\mathrm{C}^{*}$-algebra is $\mathcal{Z}$-stable and one often looks for classification ``up to $\mathcal{Z}$-stability'', or ``localized at $\mathcal{Z}$'' in the sense of \cite{Win:localizingEC}.

In \cite{Win:localizingEC} (see also \cite{Lin:localizingECappendix,LinNiu:KKlifting, Lin:asu-class, LinNiu:RationallyTAI}), the second named author showed that classification results up to $\mathcal{Z}$-stability can often be deduced from classification up to $\mathcal{U}$-stability where $\mathcal{U}$ is a UHF algebra of infinite type. This turns out to be particularly useful for those situtations when $\mathrm{C}^{*}$-algebras $A$, $B$ need not have many projections. In this case, one may tensor with some UHF algebra and work instead with the algebras $A \otimes \mathcal{U}$ and $B \otimes \mathcal{U}$. Not only does the UHF algebra bring along its many projections, we also gain many useful properties for the resulting tensor products that may not be present in the original $\mathrm{C}^{*}$-algebras such as strict comparison and $\mathcal{Z}$-stability.

The general strategy is then to tensor with a UHF algebra and show that a \mbox{$\mathrm{C}^{*}$-algebra} is of some classifiable type. In particular it is known that classification up to $\mathcal{Z}$-stability is true for simple $\mathrm{C}^{*}$-algebras satisfying the UCT which are TAF or TAI after tensoring with UHF algebras \cite{LinNiu:KKlifting, Lin:asu-class, LinNiu:RationallyTAI}. This strategy has already proven successful with regards to $\mathrm{C}^{*}$-algebras of minimal dynical systems of infinite compact metric spaces where projections separate tracial states: in the case of finite dimensional metric spaces, these are completely classified  \cite{TomsWinter:minhom, TomsWinter:PNAS}; more generally the result holds up to $\mathcal{Z}$-stability \cite{StrWin:Z-stab_min_dyn}.

Our original motivating examples come from minimal dynamical systems.  In \cite[Section 5]{Con:Thom}, Connes points out that there are minimal homeomorphisms $\alpha : S^n \to S^n$ of the $n$-dimensional sphere $S^n$ for $n=2k+1$, $k \in \mathbb{N} \setminus \{0\}$ such that the resulting $\mathrm{C}^{*}$-algebra crossed product $\mathcal{C}(S^n) \rtimes_{\alpha} \mathbb{Z}$ has no nontrivial projections. In \cite{Wind:not_uniquely_ergo} it was shown that such homeomorphisms exist with any prescribed number of invariant probabibility measures. Via the correspondence between $\alpha$-invariant Borel probabilty measures on the space and tracial states on the $C^*$-algebra, as soon as we move beyond the uniquely ergodic case, we no longer have projections separating tracial states and cannot apply the classification theorem in \cite{TomsWinter:minhom, TomsWinter:PNAS}. An inspection of the invariant (see \cite[Proposition 5.3]{LinQPhil:KthoeryMinHoms}) suggests that these $\mathrm{C}^{*}$-algebras should all be TAI after tensoring with a UHF algebra (see \cite[Sections 5 and 6]{LinNiu:RationallyTAI} for the range of the invariant of such $\mathrm{C}^{*}$-algebras).   By \cite[Theorem 4.6]{StrWin:Z-stab_min_dyn} it is enough to show that a large $\mathrm{C}^{*}$-subalgebra, denoted $(\mathcal{C}(S^n) \rtimes_{\alpha} \mathbb{Z})_y$, obtained by breaking the orbit at a point $y \in S^n$ is TAI after tensoring with $\Q$. Since such a $\mathrm{C}^{*}$-subalgebra can be written as an inductive limit of recursive subhomogeneous (RSH) $\mathrm{C}^{*}$-subalgebras (see \cite[Section 3]{LinQPhil:KthoeryMinHoms}), the problem becomes more tractable.
 
It has thus become prudent to develop techniques for dealing with cases with few projections and more complicated tracial state spaces in the setting of $\mathrm{C}^{*}$-algebras which can be locally approximated by recursive subhomogeneous $\mathrm{C}^{*}$-algebras. Our main result does this in the case that the $\mathrm{C}^{*}$-algebra can be locally approximated by RSH algebras with a decomposition into base spaces $X_0, X_1, \dots, X_R$ that can be arranged so that we can extend a projection from one space to the next in a suitable way (see Section~\ref{lifting} for more details) and under the assumption that there are only finitely many extremal tracial states all inducing the same state on the $K_0$-group. We prove that such $\mathrm{C}^{*}$-algebras are TAI after tensoring with the universal UHF algebra $\Q$ (that is, the UHF algebra satisfying $K_0(\Q) = \mathbb{Q}$). 

To show that a limit of recursive subhomogeneous algebras is TAI, we interpret the building blocks as matrix valued bundles and implement path-like structures in their base spaces. The technical difficulty comes from the fact that these base spaces are non-Hausdorff, and that such paths cannot always exist due to $K$-theory obstructions. We then enrich the $K_{0}$-groups to obtain rational vector spaces in which certain linear equations (determined by matrix sizes at the non-Hausdorff phase transitions) can be solved (provided the obstructions vanish); these solutions now show us how to arrange the paths through the non-Hausdorff base spaces.

As an application,  this gives classification, via results of Lin \cite{Lin:asu-class} (also Lin and Niu  \cite{LinNiu:RationallyTAI}), of the $\mathrm{C}^{*}$-algebras in question, provided that they can always be approximated by RSH algebras that in addition have finite topological dimension. This includes the examples of Elliott in \cite{Ell:invariant} (at least for finitely many extremal tracial states all inducing the same state on the $K_0$-group). Note in particular that these classification results cover the case where projections do not separate tracial states. (Similar results were so far only known in much more specialized situations.)

For minimal dynamical systems we show that our main theorem gives another classification of some examples given by Lin and Matui in \cite{LinMat:CantorT1} of minimal homeomorphisms on the product of the Cantor set and the circle. Such results give us confidence that similar techniques can be extended to broader classes, including those $\mathrm{C}^{*}$-algebras arising from the minimal dynamical systems $(S^n, \alpha)$ of odd dimensional spheres.

The paper is organized as follows. In Sections \ref{tracial} and \ref{rsh section} we introduce notation and definitions for tracial approximation and recursive subhomogeneous $\mathrm{C}^{*}$-algebras. Section \ref{tracial} includes results which allow us to simplify the verification of the TAI properties. Section~\ref{lifting} provides a result for lifting a projection along the stages of a recursive subhomogeneous algebra.  In Section~\ref{paths} we introduce $(\F, \eta)$-excisors and $(\F, \eta)$-bridges, which are the key tool used to cut out a large interval algebra. After developing a calculus for the $(\F, \eta)$-excisors and bridges, Sections~\ref{connected}  and \ref{excising traces} show how these are related to a given recursive subhomogeneous decomposition and how we can control the tracial weights of an $(\F, \eta)$-bridge. The main lemma in Section~\ref{linear algebra} uses linear algebra to manipulate $(\F, \eta)$-bridges to find a tracially large path through the algebra. Section~\ref{intervals} provides the technical results to find an interval that is large on all traces as well as a method for moving this interval from a general position and placing it underneath the discrete model given by the $(\F, \eta)$-path. The proof of the main result is then given in Section~\ref{outlook} where  applications and an outlook are also discussed.

\tableofcontents

\newpage
%%%%%%%%%%%%%%%%%%%%%%%%%%%%%%%%%%%%%%%%%%%%%%%%

\section{UHF-stable tracial approximation} \label{tracial}

\noindent
We begin with the general definition for a tracially approximately $\mathcal{S}$ (TA$\mathcal{S}$) algebra, where $\mathcal{S}$ can be any class of separable unital $\mathrm{C}^{*}$-algebras. For the time being, we will only be interested in the case where the class $\mathcal{S}$ consists of interval algebras. Despite the fact that both are notions of approximation for $\mathrm{C}^{*}$-algebras, the TA$\mathcal{S}$ class may be much broader than the A$\mathcal{S}$ class, that is, inductive limits of $C^*$-algebras in the class $\mathcal{S}$. For example, the class of simple unital tracially approximately interval (TAI) class contains all simple unital AH algebras with slow dimension growth, clearly not all of which are AI \cite[Section 10]{Lin:simpleNuclearTR1}.

\bn
\begin{ndefn}\textup{(cf. \cite{Lin:amenable, EllNiu:tracial_approx})}
\label{TAS}
Let $\mathcal{S}$ denote a class of separable unital  $\mathrm{C}^{*}$-algebras. Let $A$ be a simple unital $\mathrm{C}^{*}$-algebra. Then $A$ is tracially approximately $\mathcal{S}$ (or $\mathrm{TA}\mathcal{S}$) if the following holds.

For every finite subset $\mathcal{F} \subset A$, every $\epsilon > 0$, and every nonzero positive element $c \in A$, there exist a projection $p \in A$ and a unital $\mathrm{C}^{*}$-subalgebra $B \subset pAp$ with $1_{B}=p$ and $B \in \mathcal{S}$ such that:
\begin{enumerate}
\item $\| pa - ap \| < \epsilon$ for all $a \in \mathcal{F}$,
\item $\dist(pap, B) < \epsilon$ for all $a \in \mathcal{F}$,
\item $1_{A}-p$ is Murray--von Neumann equivalent to a projection in $\overline{cAc}$.
\end{enumerate}
\end{ndefn}
\en

In this paper we consider tracial approximation by the class I of interval algebras. An interval algebra is a $\mathrm{C}^{*}$-algebra $A$ of the form 
$$A = \bigoplus_{n=1}^N \mathcal{C}(X_n) \otimes M_{r_n}$$
for some $N \in \mathbb{N} \setminus \{0\}$, where $X_n = [0,1]$ or $X_n$ is a single point, and $r_n \in \mathbb{N} \setminus \{0\}$, $0 \leq n \leq N$.

We denote the universal UHF algebra by $\Q$. This is the unique UHF algebra whose $K_0$-group is isomorphic to $\mathbb{Q}$. Recall that if $\uhf$ is a UHF algebra of infinite type then it is strongly self-absorbing \cite[Example 1.14 (i)]{TomsWin:ssa}:  there are a $^*$-isomorphism $\phi : \uhf \to \uhf \otimes \uhf$ and a sequence of unitaries $(v_n)_{n \in \mathbb{N}}$ in $\uhf \otimes \uhf$ such that $\| v_n^* \phi(d) v_n - \id_{\uhf}(d)\otimes 1_{\uhf} \| \to 0$ as $n \to \infty$ for every $d \in \mathcal \uhf$ \cite[Definition 1.3 (iv)]{TomsWin:ssa}.  

Any $\mathrm{C}^{*}$-algebra in the class I can be written as a finitely presented universal $\mathrm{C}^{*}$-algebra (i.e.\ with finitely many generators and relations) and is semiprojective. In particular, any $A \in \mathrm{I}$ has stable, hence weakly stable, relations \cite{Lor:lifting}. Therefore we may make use of the following lemma, which says that to prove TAI after tensoring with $\Q$ it is enough to show that the approximating $\mathrm{C}^{*}$-algebras can always be chosen to have units that are bounded above zero in trace. The proof uses the same geometric series argument as the one given in \cite[Lemma 3.2]{Win:Z-class}.

\bn
\begin{nlemma} \label{away from 0} Let $A$ be a separable simple unital stably finite exact $\mathrm{C}^{*}$-algebra and let $\uhf$ be a UHF algebra of infinite type. Suppose $\mathcal{S}$ is a class of $\mathrm{C}^{*}$-algebras that can be finitely presented with weakly stable relations (as universal $\mathrm{C}^{*}$-algebras), contains all finite dimensional $\mathrm{C}^{*}$-algebras, and is closed under direct sums.  Then $A \otimes \uhf$ is $\mathrm{TA}\mathcal{S}$ if and only if there is an $n \in \mathbb{N}$ such that, for any $\epsilon >0$ and any finite subset $\mathcal{F} \subset A\otimes \uhf$, there exist a projection $p \in A \otimes \uhf$ and a unital $\mathrm{C}^{*}$-subalgebra $B \subset p(A \otimes \uhf)p$ and $B \in \mathcal{S}$ such that:
\begin{enumerate}
\item $\|pb - bp \| < \epsilon$ for all $b \in \mathcal{F}$,
\item $\dist(p b p, B) < \epsilon$ for all $b \in \mathcal{F}$,
\item $\tau(p) > 1/n$ for all $\tau \in T(A \otimes \uhf)$.
\end{enumerate}
\end{nlemma}

\begin{nproof} If $A \otimes \uhf$ is TA$\mathcal{S}$ then (i) and (ii) are easily satisfied from the definition of TA$\mathcal{S}$. To show (iii) with, for example, $n=2$, take a positive element $c \in A \otimes \uhf$ with $\tau(c) \leq 1/2$ for all $\tau \in T(A \otimes \uhf)$ and use the fact that $A \otimes \uhf$ has strict comparison \cite[Theorem 5.2(a)]{Ror:uhfII}.
Now let a finite subset $\F \subset A \otimes \uhf$, $\epsilon >0$ and a nonzero positive element $c \in A\otimes \uhf$ be given, and suppose that $A \otimes \uhf$ satisifies (i), (ii), (iii) with respect to some $n \in \mathbb{N}$. We show that $A \otimes \uhf$ is TA$\mathcal{S}$; the proof is almost identical to that of Lemma 3.2 of \cite{Win:Z-class}. First we note that $A \otimes \uhf$ has property (SP) (every nonzero hereditary $\mathrm{C}^{*}$-subalgebra has a nonzero projection) since $A \otimes \uhf$ has strict comparison and projections that are arbitrarily small in trace. Thus we find a projection $q \in \overline{c(A \otimes \uhf) c}$, just as in \cite[Lemma 3.2]{Win:Z-class}. 

We inductively construct $\mathrm{C}^{*}$-algebras $B_i \subset A \otimes \uhf$ with each $B_i \in \mathcal{S}$. 

As in \cite[Lemma 3.2]{Win:Z-class} the initial $B_0$ exists by assumption.  

The construction of $B_{i+1}$ from $B_i$ is similar to the construction in Lemma 3.2 of \cite{Win:Z-class}. We cannot apply Lemma 3.4 of \cite{Win:Z-class} directly, even though $A \otimes \uhf$ is simple and unital and has the comparability property, since we do not want to make the assumption that $K_0(A \otimes \uhf)_+$ has dense image in the positive affine functions $T(A \otimes \uhf)$. However, the result will still hold by choosing the projection $e$ in that proof to be of the form $(1_{A \otimes \uhf} - p) \otimes q$ for some projection $q \in \uhf$ (using the fact that $\uhf$ is strongly self-absorbing) satisfying
$$1/(t+1) < \tau_{\uhf}(q) < 1/t$$
where $\tau_{\uhf}$ is the unique tracial state on $\uhf$. The projection $e$ then satisfies the requirements of the projection in the proof, and the results of \cite[Lemma 3.4]{Win:Z-class} hold. Thus we get the finite dimensional $\mathrm{C}^{*}$-algebras $C_0, C_1$ and $D$ as in \cite[Lemma 3.2]{Win:Z-class}. 

Let $\mathcal{G} := \{x_1, \dots, x_n, 1_{B_i}\} \subset B_i$ where $x_1, \dots, x_n$ are generators for $B_i$. Let $\gamma > 0$ be as in  \cite[Lemma 3.2]{Win:Z-class}. Since $B_i$ has weakly stable relations, there is a $\tilde{\vartheta} > 0$ with the following property: If $E$ is another $\mathrm{C}^{*}$-algebra, $p \in A$ a projection and $\phi: B_i \to E$ a $^*$-homomorphism satisfying $\|p \phi(b) - \phi(b) p \| < \tilde{\vartheta}$ for all $b \in \mathcal{G}$, then there is a $^*$-homomorphism $\tilde{\phi} : B_i \to pEp$ satisfying $\| \tilde{\phi}(b) - p \phi(b) p \| < \gamma$ for all $b \in \mathcal{G}$. 

Now choose $0 < \vartheta < \min\{\gamma, \tilde{\vartheta} \}$ such that the assertion of  \cite[Proposition 3.3]{Win:Z-class}
holds for the finite dimensional algebra $D$. Set $\tilde{\F} := \F \cup \mathcal{G} \cup \kappa(D)^1$ where $\kappa :  D \to A \otimes \uhf$ is a $^*$-homomorphism given by \cite[Lemma 3.4]{Win:Z-class} and  $\kappa(D)^1$ denotes the unit ball of $\kappa(D)$. 

By hypothesis there is a $\mathrm{C}^*$-algebra $F \subset A \otimes \uhf$, $F \in \mathcal{S}$ satisfying d), e), f) of \cite[Proposition 3.3]{Win:Z-class} with respect to $\tilde{\F}$. As in \cite[Lemma 3.2]{Win:Z-class}, d) and the choice of $\vartheta$ provides the $^*$-homomorphisms 
$$\varrho : B_i \to (1_{A \otimes \uhf} - 1_F)A \otimes \uhf (1_{A \otimes \uhf} - 1_F)$$
satisfying 
$$\| \varrho(b) -  (1_{A \otimes \uhf} - 1_F)b(1_{A \otimes \uhf} - 1_F)\| < \gamma \text{ for all } b \in \mathcal{G}$$
and 
$$\bar{\kappa} : D \to  (1_{A \otimes \uhf} - 1_F)A \otimes \uhf (1_{A \otimes \uhf} - 1_F)$$
such that
$$\| \bar{\kappa} -  (1_{A \otimes \uhf} - 1_F) \kappa(d) (1_{A \otimes \uhf} - 1_F)\| < \gamma \cdot \| d \| \text{ for all } 0 \neq d \in D.$$

Set $B_{i+1} := \varrho(B_i) \oplus F.$ Then one easily checks that the same calculations given in \cite[Proposition 3.3]{Win:Z-class} can be used to complete the proof. 
\end{nproof} 
\en

%%%%%%%%%%%%%%%%reduction to taking finite subset in A
The next lemma shows we need only consider finite subsets of $A \otimes \Q$ of a simplified form, that is, essentially the only difficulty lies in approximating elements from $A$.

\bn
\begin{nlemma} \label{special F} Let $\mathcal{S}$ denote a class of separable unital  $\mathrm{C}^{*}$-algebras that is closed under tensoring with finite dimensional $\mathrm{C}^{*}$-algebras. Let $A$ be a separable unital $\mathrm{C}^{*}$-algebra with $T(A) \neq 0$ and let $0 < \eta \leq 1$ such that, for any $\epsilon > 0$ and any finite subset $\mathcal{G} \subset A$, there are a projection $p \in A \otimes \Q$ and a unital $\mathrm{C}^{*}$-subalgebra $B \subset p (A \otimes \Q) p$ with $1_B = p$ and  $B \in \mathcal{S}$ such that
\begin{enumerate}
\item $\| p(a \otimes 1_{\Q}) - (a \otimes 1_{\Q})p \| < \epsilon$ for all $a \in \mathcal{G}$,
\item $\dist(p(a \otimes 1_{\Q}) p, B) < \epsilon$ for all $a \in \mathcal{G}$, 
\item $\tau(p) \geq \eta$ for all $\tau \in T(A \otimes \Q)$.
\end{enumerate}
Then, for any $\epsilon > 0$ and any finite subset $\F \in A \otimes \Q$, there are a projection $q \in A \otimes \Q$ and a unital $\mathrm{C}^{*}$-subalgebra $C \subset q (A \otimes \Q) q$ with $1_{C} = q$ and  $C \in \mathcal{S}$ such that
\begin{enumerate}[resume]
\item $\| q a - a q \| < \epsilon$ for all $a \in \F$,
\item $\dist(qa q, C) < \epsilon$ for all $a \in \F$, 
\item $\tau(q) \geq \eta$ for all $\tau \in T(A \otimes \Q)$.
\end{enumerate}

\end{nlemma}

\begin{nproof} The proof essentially appears in the proof of \cite[Lemma 4.4]{StrWin:Z-stab_min_dyn}.
Let $\epsilon > 0$ and let $\mathcal{F} \subset A \otimes \Q$ be a finite subset. Using the identification 
\[ A \otimes \Q \cong A \otimes M_{S} \otimes \Q \cong A \otimes \Q \otimes M_{S}, \]
for $S \in \mathbb{N}$, we may assume that the finite set is of the form
\[ (\{1_A \}\otimes \{1_{\Q} \}  \otimes \mathcal{B} )\cup ( \mathcal{G} \otimes \{1_{\Q} \} \otimes  \{1_{M_{S}} \} ) \]
where $S\in \mathbb{N}$,  $\mathcal{B}$ is a finite subset of $M_{S}$ and $\mathcal{G}$ is a finite subset of $A$. We may further assume that $1_A \in \mathcal{G}$ and also that $1_{M_{S}} \in \mathcal{B}$. Then we have 
\[ \mathcal{F} = \mathcal{G}  \otimes \{ 1_{\Q} \} \otimes \mathcal{B}.\] 

By assumption, there exists a $B \in \mathcal{S}$ and a projection $p = 1_B$ satisfying properties (i) -- (iii) of the lemma for the finite set $\mathcal{G} $, with $\epsilon / \max (\{ \|b\| \mid b \in \mathcal{B} \}, 1)$ in place of $\epsilon$. 
 
Define $C = B \otimes M_S$ and $q := 1_D = p \otimes 1_{M_{S}} \in A \otimes \Q \otimes M_{S}$. The fact that $q$ and $C$ satisfy properties (iv) and (v) of the lemma for $\mathcal{\tilde{F}}$ and $\epsilon$ is shown in the proof of \cite[Lemma 4.4]{StrWin:Z-stab_min_dyn}.
 
To show (vi), simply observe that $\tau \in  T(A \otimes \Q \otimes M_{S})$ is of the form $\tau_1 \otimes \tau_2$ where $\tau_1 \in T(A \otimes \Q)$ and $\tau_2 \in T(M_{S})$. Then 
\[ \tau(q) = \tau(p \otimes 1_{M_S}) = \tau_1(p)\tau_2(1_{M_S}) = \tau(q) \geq \eta.  \] 
\end{nproof}
\en

%%%%%%%%%%%%%%%%%%%%%%%%%%%%%%%%%%%%%%%%%%%%%%%%
\section{Recursive subhomogeneous $\mathrm{C}^{*}$-algebras} \label{rsh section}

\noindent
In \cite{Phillips:recsub}, Phillips introduced the notion of recursive subhomogeneous algebras. These are subhomogeneous $\mathrm{C}^{*}$-algebras (that is, all irreducible representations are bounded in dimension) which arise as iterated pullbacks of homogeneous $\mathrm{C}^*$-algebras.

\bn
\begin{ndefn}\cite[Definition 1.1]{Phillips:recsub}
\label{rsh}
A recursive subhomogeneous (RSH) algebra is a $\mathrm{C}^{*}$-algebra with the following recursive definition:
\begin{enumerate}
\item  Let $X$ be a compact Hausdorff space and \mbox{$n \in \mathbb{N}$}. Then $\mathcal{C}(X, M_n)$ is a recursive subhomogeneous algebra.
\item Let $A$ be a recursive subhomogeneous $\mathrm{C}^{*}$-algebra, $X$ a compact Hausdorff space and \mbox{$n \in \mathbb{N}$}. Suppose $\Omega \subset X$ is a closed (possibly empty) subset, $\phi:  A \to \mathcal{C}(\Omega, M_n)$ is a unital $^{*}$-homomorphism, and let 
\[
\rho : \mathcal{C}(X, M_n) \to \mathcal{C}(\Omega, M_n)
\]
be the restriction map. Then the pullback
$$A \oplus_{\mathcal{C}(\Omega, M_n)} \mathcal{C}(X, M_n) = \{ (a, f) \in A \oplus \mathcal{C}(X, M_n) \mid \phi(a) = \rho(f) \}$$
is a recursive subhomogeneous $\mathrm{C}^{*}$-algebra.
\end{enumerate}
\end{ndefn}
\en

We will restrict to the case where the $X$ in the above definition is metrizable so that the resulting $\mathrm{C}^{*}$-algebra is separable. Note also that Definition~\ref{rsh} implies that a recursive subhomogeneous $\mathrm{C}^{*}$-algebra is unital.

For a given $\mathrm{C}^{*}$-algebra, its recursive subhomogeneous decomposition (if it exists) is not unique; for us it will be important to keep track of the actual decompositions.

If $B$ is a recursive subhomogeneous algebra, then we may write (cf. \cite[Definition 1.5]{Phillips:recsub})
$$B = \left( \dots \left( \left( C_0 \oplus_{C_1^{(0)}} C_1 \right) \oplus_{C_2^{(0)}} C_2 \right) \dots \right) \oplus_{C_R^{(0)}} C_R,$$
where $C_l = \mathcal{C}(X_l) \oplus M_{n_l}$ for some compact metrizable $X_l$ and some integer $n_l \geq 1$, $l \in \{0, \dots, R\}$ and $C_l^{(0)} = \mathcal{C}(\Omega_l) \oplus M_{n_l}$ for a closed subset $\Omega_l \subset X_l$.

For $0 \leq l \leq R$, define the $l^{\text{th}}$-stage of $B$ to be the $\mathrm{C}^{*}$-algebra obtained by truncating the recursion after the $l^{\text{th}}$ step,
$$B_l = \left( \dots \left( \left( C_0 \oplus_{C_1^{(0)}} C_1 \right) \oplus_{C_2^{(0)}} C_2 \right) \dots \right) \oplus_{C_l^{(0)}} C_l.$$
Note that $B_R =  B$.

\bn
\begin{ndefn}
\label{dK}
Let $B$ be a unital recursive subhomogeneous $\mathrm{C}^{*}$-algebra. A recursive subhomogeneous decomposition 
\[
[B_{l},X_{l},\Omega_{l},r_{l},\phi_{l}]_{l=1}^{R}
\]
for $B$ consists of compact Hausdorff spaces $\Omega_{l} \subset X_{l}$ and $r_{l}\in \mathbb{N}$, unital $\mathrm{C}^{*}$-algebras $B_{l}$ for $l \in \{1,\ldots,R\}$, and of unital $^{*}$-homomorphisms 
\[
\phi_{l}:B_{l} \to \mathcal{C}(\Omega_{l+1}) \otimes M_{r_{l+1}}
\]
for $l \in \{1,\ldots,R-1\}$, such that 
\[
\Omega_{1} = \emptyset, \; B_{1}= \mathcal{C}(X_{1}) \otimes M_{r_{1}}, \; B_{R}=B
\]
and such that we have pullback diagrams
\begin{equation}
\label{dK1}
\xymatrix{
B_{l+1} \ar@{>>}[r] \ar[d] & B_{l} \ar[d]^{\phi_{l}} \\
\mathcal{C}(X_{l+1}) \otimes M_{r_{l+1}} \ar@{>>}[r]& \mathcal{C}(\Omega_{l+1}) \otimes M_{r_{l+1}} ,
}
\end{equation}
where the lower horizontal map is restriction. We then have a canonical unital embedding
\[
\iota_{B}: B \hookrightarrow \mathcal{C}(X_{1}) \otimes M_{r_{1}} \oplus \ldots \oplus \mathcal{C}(X_{R}) \otimes M_{r_{R}},
\]
canonical quotient maps
\[
\psi_{l}: B \to B_{l}
\]
and canonical embeddings
\[
\iota_{l}: \mathcal{C}_0(X_{l} \setminus \Omega_{l}) \otimes M_{r_{l}} \hookrightarrow B_{l}
\]
for $l \in \{1,\ldots,R\}$.

We will usually assume that $X_{l+1} \setminus \Omega_{l+1} \neq \emptyset$, for otherwise the horizontal maps in \eqref{dK1} are just equalities.
\end{ndefn}
\en

\bn
\begin{nremark}
\label{r-dK}
In the situation above, if $x \in X_{l} \setminus \Omega_{l}$ for some $l \in \{1,\ldots,R\}$, then the map $(\mathrm{ev}_{x} \otimes \mathrm{id}_{M_{r_{l}}})  \circ \iota_{B} :B \to M_{r_{l}}$ is surjective.
\end{nremark}
\en

%%%%%%%%%%%%%%%%%%%%%%%%%%%%%%%%%%%%%%%%%%%%%%%%
\section{Lifting projections} \label{lifting}

\bn
\begin{ndefn}
Let $B$ be a unital recursive subhomogeneous $\mathrm{C}^{*}$-algebra with recursive subhomogeneous decomposition 
\[
[B_{l},X_{l},\Omega_{l},r_{l},\phi_{l}]_{l=1}^{R}.
\]
We say that projections can be lifted along $[B_{l},X_{l},\Omega_{l},r_{l},\phi_{l}]_{l=1}^{R}$, if for any $N \in \mathbb{N}$, any $l \in \{1,\ldots,R-1\}$ and any projection $p \in B_{l} \otimes M_{N}$ there is a projection \mbox{$\bar{p} \in B_{l+1} \otimes M_{N}$} lifting $p$.
\end{ndefn}
\en

\bn
\begin{nprop}
\label{1-dim-lifting}
Let $X$ be compact metrizable with $\mathrm{dim} X \le 1$. Let $k,r \in \mathbb{N}$, $\Omega \subset X$ a closed subspace and $p \in \mathcal{C}(\Omega,M_{r})$ a projection with constant rank $k$. 

Then there is a projection $\bar{p} \in \mathcal{C}(X,M_{r})$ extending $p$.
\end{nprop}

\begin{nproof}
It is straightforward to find a closed neighborhood $W$ of $\Omega$ and a projection 
\[
\tilde{p} \in \mathcal{C}(W,M_{r})
\]
extending $p$. Let $U \subset X$ be an open subset such that
\[
\Omega \subset U \subset W.
\]
Let $(W_{\lambda})_{\Lambda}$ be a finite collection of open subsets of $X$ such that 
\begin{itemize}
\item[(i)]  $W \subset \bigcup_{\Lambda} W_{\lambda}$
\item[(ii)] $\|\tilde{p}(x) - \tilde{p}(x')\| \le \frac{1}{2}$ whenever $x,x' \in \overline{W_{\lambda}}$ for some $\lambda \in \Lambda$
\item[(iii)] $W_{\lambda} \subset U$ if $W_{\lambda} \cap \Omega \neq \emptyset$.
\end{itemize}

From (ii) it is not hard to see that for each $\lambda$, $\tilde{p}|_{\overline{W_{\lambda}}}$ is homotopic to a constant projection of rank $k$; this yields projections
\begin{equation}
\label{1-dim-lifting-w2}
p_{\lambda} \in \mathcal{C}(\overline{W_{\lambda}} \times [0,1],M_{r}) 
\end{equation}
such that
\[
p_{\lambda}(x,t) = \left\{ \begin{array}{ll}
\tilde{p}(x), & \mbox{for } t \in [\frac{2}{3},1] \\
1_{k}, & \mbox{for } t \in [0,\frac{1}{3}]
\end{array}
\right.
\]
(where we think of $1_{k}$ as sitting in the upper left corner of $M_{r}$).

Since $\mathrm{dim}X \le 1$, there is   a finite open cover $(V_{\gamma})_{\Gamma}$  of $X$ refining the open cover consisting of $W_{\lambda}, \lambda \in \Lambda$, and $X \setminus W$, and such that 
\begin{equation}
\label{1dl1}
V_{\gamma_{0}} \cap V_{\gamma_{1}} \cap V_{\gamma_{2}} = \emptyset
\end{equation}
whenever $\gamma_{0}, \gamma_{1},\gamma_{2} \in \Gamma$ are pairwise distinct. Let
\[
(h_{\gamma})_{\Gamma}
\]
be a partition of unity subordinate to $(V_{\gamma})_{\Gamma}$. Set
\[
\Gamma':= \{\gamma \in \Gamma \mid V_{\gamma} \cap \Omega \neq \emptyset\}
\]
and
\[
\Gamma'':= \{\gamma \in \Gamma \setminus \Gamma' \mid V_{\gamma} \cap V_{\gamma'} \neq \emptyset \mbox{ for some } \gamma' \in \Gamma'\}
\]
%and
%\[
%\Gamma''':= \Gamma \setminus (\Gamma' \cup \Gamma'').
%\]
Note that by (ii) above, for any $\gamma' \in \Gamma'$ there is $\lambda(\gamma') \in \Lambda$ such that
\begin{equation}
\label{1-dim-lifting-w1}
V_{\gamma'} \subset W_{\lambda(\gamma')} \subset U \subset W.
\end{equation}

We now define $\bar{p}$, observing  that for each $x \in X$, by \eqref{1dl1} there are at most two indices $\gamma,\gamma' \in \Gamma$ such that $h_{\gamma}(x),h_{\gamma'}(x) \neq 0$. 

\underline{Case 1:}
There is only one index $\gamma \in \Gamma$ such that $h_{\gamma}(x) \neq 0$; in this case, $h_{\gamma}(x) = 1$.

\underline{Case 1a:} If $\gamma \in \Gamma'$, set 
\[
\bar{p}(x):= \tilde{p}(x);
\]
this is well defined by \eqref{1-dim-lifting-w1}.

\underline{Case 1b:} If $\gamma \in \Gamma \setminus \Gamma'$, set 
\[
\bar{p}(x):= 1_{k}.
\]

\underline{Case 2:}
There are two distinct indices $\gamma,\gamma' \in \Gamma$ such that $h_{\gamma}(x), h_{\gamma'}(x) \neq 0$. 

\underline{Case 2a:}
If $\gamma,\gamma' \in \Gamma'$, set 
\[
\bar{p}(x):= \tilde{p}(x);
\]
again, this is well defined by \eqref{1-dim-lifting-w1}.

\underline{Case 2b:}
If $\gamma,\gamma' \in \Gamma \setminus \Gamma'$, set 
\[
\bar{p}(x):= 1_{k}.
\]

\underline{Case 2c:}
If $\gamma \in \Gamma \setminus \Gamma'$, $\gamma' \in \Gamma'$, then $h_{\gamma}(x) + h_{\gamma'}(x) = 1$, so $h_{\gamma'}(x) \in [0,1]$ and by \eqref{1-dim-lifting-w1} and \eqref{1-dim-lifting-w2} we may set
\[
\bar{p}(x):= p_{\lambda(\gamma')}(x,h_{\gamma'}(x)).
\]

We have now defined a projection valued map
\[
\bar{p}:X \to M_{r}
\]
which by construction clearly extends $p$ (note that if $x \in \Omega$, then only Cases 1a and 2a occur). It remains to check that $\bar{p}$ is continuous. 

So let $x \in X$. In Case 2, there are $\gamma \neq \gamma' \in \Gamma$ with $h_{\gamma}(x),h_{\gamma'}(x) \neq 0$. But then $h_{\gamma}(y), h_{\gamma'}(y) \neq 0$ for all $y$ in some small neighborhood $V_{x}$ of $x$. In Case  2a, note that the map $y \mapsto \tilde{p}(y)$ is continuous; in Case  2b, $\bar{p}(y) = \bar{p}(x)$ for $y \in V_{x}$; in Case 2c, the map
\[
y \mapsto p_{\lambda(\gamma')} (y,h_{\gamma'}(y))
\] 
is continuous on $V_{x}$ since $h_{\gamma'}$ and $p_{\lambda(\gamma')}$ are.

In Case 1, we have $h_{\gamma}(x) = 1$. But then there is some neighborhood $V_{x}$ of $x$ such that $h_{\gamma}(y) \ge \frac{2}{3}$ for all $y \in V_{x}$, and we obtain
\[
\bar{p}(y) = \left\{
\begin{array}{ll}
\tilde{p}(y), & \mbox{if } \gamma \in \Gamma'\; \mbox{ (in Case 1a, 2a or 2c for }y \mbox{ in place of }x)\\
1_{k}, & \mbox{if } \gamma \in \Gamma \setminus \Gamma'  \; \mbox{ (in Case 1b  for }y \mbox{ in place of }x)
\end{array}
\right.
\]
for $y \in V_{x}$, whence $\bar{p}$ is continuous at $x$.
\end{nproof}
\en

\bn
\begin{ncor} \label{dim1}
Let $B$ be a unital recursive subhomogeneous $\mathrm{C}^{*}$-algebra with decomposition 
\[
[B_{l},X_{l},\Omega_{l},r_{l},\phi_{l}]_{l=1}^{R}.
\]
Assume that $\mathrm{dim} X_{l} \le 1$ for $l \ge 2$. 

Then projections can be lifted along  $[B_{l},X_{l},\Omega_{l},r_{l},\phi_{l}]_{l=1}^{R}$.
\end{ncor}

\begin{nproof}
Obvious from Proposition~\ref{1-dim-lifting} and Definition~\ref{dK}.
\end{nproof}
\en

%%%%%%%%%%%%%%%%%%%%%%%%%%%%%%%%%%%%%%%%%%%%%%%%
\section{Approximately excising approximate paths} \label{paths}

\noindent
Recall that a completely positive map has order zero when it preserves orthogonality, that is, a c.p.\ map $\phi: A \to B$ between the $\mathrm{C}^{*}$-algebras $A$ and $B$ such that, for any orthogonal positive elements $a, b \in A$ with $ab = 0$ we have $\phi(a) \phi(b) = 0$ in $B$. 

\bn
\begin{ndefn} \label{F-eta}
Let $B$ be a unital recursive subhomogeneous $\mathrm{C}^{*}$-algebra with recursive subhomogeneous decomposition 
\[
[B_{l},X_{l},\Omega_{l},r_{l},\phi_{l}]_{l=1}^{R};
\]
let $\mathcal{F} \subset B_{+}^{1}$ be finite a subset, where $B_{+}^{1}$ denotes the positive elements in the unit ball of $B$, and $\eta>0$ be given.

An $(\mathcal{F},\eta)$-excisor $(E,\rho,\sigma)$ for $B$ consists of a finite dimensional $\mathrm{C}^{*}$-algebra 
\[
\textstyle
E= \bigoplus_{l=1}^{R} E_{l},
\] 
a unital $^{*}$-homomorphism
\[
\textstyle
\rho= \oplus_{l=1}^{R} \rho_{l} : B \to \bigoplus_{l=1}^{R} E_{l} = E
\]
and an isometric c.p.\ order zero map
\[
\textstyle
\sigma = \oplus_{l=1}^{R} \sigma_{l} : \bigoplus_{l=1}^{R} E_{l} = E \to B \otimes \mathcal{Q}
\]
such that
\[
\| \sigma(1_{E})(b \otimes 1_{\mathcal{Q}}) - \sigma \rho(b)\| <\eta \mbox{ for } b \in \mathcal{F}.
\]

We say $(E,\rho,\sigma)$ is compatible with $[B_{l},X_{l},\Omega_{l},r_{l},\phi_{l}]_{l=1}^{R}$, if each $\rho_{l}$ factorizes through
\[
\xymatrix{
B \ar[d]_{\psi_{l}} \ar[r]^{\rho_{l}} & E_{l} \\
B_{l} \ar[r]^-{\check{\psi}_{l}} & \mathcal{C}(\check{X}_{l}) \otimes M_{r_{l}} \ar[u]_{\check{\rho}_{l}}
}
\]
for some compact $\check{X}_{l} \subset X_{l} \setminus \Omega_{l}$.

If $(E,\rho,\sigma)$ is as above and 
\[
\kappa:E \to \mathcal{Q}
\]
is a unital $^{*}$-homomorphism, we say $(E,\rho,\sigma,\kappa)$ is a weighted $(\mathcal{F},\eta)$-excisor compatible with $[B_{l},X_{l},\Omega_{l},r_{l},\phi_{l}]_{l=1}^{R}$.
\end{ndefn}
\en

\bn
\begin{ndefn}
\label{L}
Let $B$ be a unital recursive subhomogeneous $\mathrm{C}^{*}$-algebra with recursive subhomogeneous decomposition 
\[
[B_{l},X_{l},\Omega_{l},r_{l},\phi_{l}]_{l=1}^{R};
\]
let $\mathcal{F} \subset B_{+}^{1}$ finite and $\eta>0$ be given. Let $(E_{i},\rho_{i},\sigma_{i},\kappa_{i})$, $i \in \{0,1\}$, be weighted $(\mathcal{F},\eta)$-excisors (compatible with $[B_{l},X_{l},\Omega_{l},r_{l},\phi_{l}]_{l=1}^{R}$). 

An $(\mathcal{F},\eta)$-bridge from $(E_{0},\rho_{0},\sigma_{0},\kappa_{0})$ to $(E_{1},\rho_{1},\sigma_{1},\kappa_{1})$ (compatible with the decomposition $[B_{l},X_{l},\Omega_{l},r_{l},\phi_{l}]_{l=1}^{R}$) consists of $K \in \mathbb{N}$ and weighted $(\mathcal{F},\eta)$-excisors (each compatible with $[B_{l},X_{l},\Omega_{l},r_{l},\phi_{l}]_{l=1}^{R}$)
\[
(E_{\frac{j}{K}},\rho_{\frac{j}{K}},\sigma_{\frac{j}{K}},\kappa_{\frac{j}{K}}), \; j \in \{1,\ldots,K-1\},
\]
satisfying
\begin{equation}
\label{L1}
\| \kappa_{\frac{j}{K}} \rho_{\frac{j}{K}}(b) - \kappa_{\frac{j+1}{K}} \rho_{\frac{j+1}{K}}(b) \| < \eta \mbox{ for } b \in \mathcal{F} \mbox{ and } j \in \{0,\ldots,K-1\}.
\end{equation}
We write
\[
(E_{0},\rho_{0},\sigma_{0},\kappa_{0}) \sim_{(\mathcal{F},\eta)}  (E_{1},\rho_{1},\sigma_{1},\kappa_{1})
\]
if such an $(\mathcal{F},\eta)$-bridge exists.
\end{ndefn}
\en

\bn
\begin{nremarks}
\label{M}
(i) Clearly, the relation $\sim_{(\mathcal{F},\eta)}$ defines an equivalence relation on the set of compatible weighted $(\mathcal{F},\eta)$-excisors (with fixed $\mathcal{F}$, $\eta$ and recursive subhomogeneous decomposition).

(ii)  If $(E,\rho,\sigma,\kappa_{i})$, $i \in \{0,1\}$, are $(\mathcal{F},\eta)$-excisors with $\kappa_{0}$ and $\kappa_{1}$ unitarily equivalent, $\kappa_{0} \approx_{\mathrm{u}} \kappa_{1}$, then
\[
(E,\rho,\sigma,\kappa_{0}) \sim_{(\mathcal{F},\eta)}  (E,\rho,\sigma,\kappa_{1})
\]
since $\kappa_{0}$ and $\kappa_{1}$ are in fact homotopic.

(iii) Let $(E,\rho,\sigma,\kappa)$ and $(E',\rho',\sigma',\kappa')$ be $(\mathcal{F},\eta)$-excisors with an embedding 
\[
\iota: E' \to E
\]
and such that 
\[
\rho'= \iota \circ \rho, \;  \sigma' = \sigma \circ \iota , \; \kappa' = \kappa \circ \iota.
\]
Then
\[
(E,\rho,\sigma,\kappa)  \sim_{(\mathcal{F},\eta)} (E',\rho',\sigma',\kappa'),
\]
with an $(\mathcal{F},\eta)$-bridge of length $K=1$.
\end{nremarks}
\en

\bn
\begin{ndefn}
\label{Ia}
Let $B$ be a unital recursive subhomogeneous $\mathrm{C}^{*}$-algebra with recursive subhomogeneous decomposition 
\[
[B_{l},X_{l},\Omega_{l},r_{l},\phi_{l}]_{l=1}^{R};
\]
let $\mathcal{F} \subset B_{+}^{1}$ finite and $\eta>0$ be given. 

(i) Let $(E_{j},\rho_{j},\sigma_{j},\kappa_{j})$, $j \in \{1,\ldots,L\}$, be weighted $(\mathcal{F},\eta)$-excisors. We say they are pairwise orthogonal if there are pairwise orthogonal projections 
\[
q_{j} \in \mathcal{Q}, \; j \in \{1,\ldots,L\},
\]
such that 
\[
\sigma_{j}(E_{j}) \subset B \otimes q_{j} \mathcal{Q} q_{j} \subset_{\mathrm{her}} B \otimes \mathcal{Q}, \; j \in \{1,\ldots,L\}.
\]

(ii) Let $(E_{j},\rho_{j},\sigma_{j},\kappa_{j})$, $j \in \{1,\ldots,L\}$, be pairwise orthogonal weighted $(\mathcal{F},\eta)$-excisors, and let 
\[
\textstyle
\gamma: \bigoplus_{j=1}^{L} \mathcal{Q} \to \mathcal{Q}
\]
be a unital $^{*}$-homomorphism. 

We define the $\gamma$-direct sum 
\[
\textstyle
\bigoplus_{\gamma} (E_{j}, \rho_{j},\sigma_{j},\kappa_{j}) := (\bigoplus_{j=1}^{L} E_{j}, \bigoplus_{j=1}^{L}\rho_{j}, \bigoplus_{j=1}^{L} \sigma_{j}, \gamma \circ (\bigoplus_{j=1}^{L} \kappa_{j})),
\]
which is easily seen to be a weighted $(\mathcal{F},\eta)$-excisor. 

If the $(E_{j},\rho_{j},\sigma_{j},\kappa_{j})$ are compatible with $[B_{l},X_{l},\Omega_{l},r_{l},\phi_{l}]_{l=1}^{R}$, then so is the \mbox{$\gamma$-direct} sum.

Since, up to unitary equivalence in $\mathcal{Q}$, the maps $ \gamma \circ (\bigoplus_{j=1}^{L} \kappa_{j})$ only depend on the positive rational weights $\nu_{j}:= \tau_{\mathcal{Q}} (\gamma(1_{j}))$, we will sometimes neglect to explicitly specify $\gamma$ and write 
\[
\textstyle
\bigoplus_{j=1}^{L}  \nu_{j} \cdot \kappa_{j} 
\]
instead of $ \gamma \circ (\bigoplus_{j=1}^{L} \kappa_{j})$ and, similarly,
\[
\textstyle
\bigoplus_{j=1}^{L}  \nu_{j} \cdot (E_{j}, \rho_{j},\sigma_{j},\kappa_{j}) 
\]
instead of $\bigoplus_{\gamma} (E_{j}, \rho_{j},\sigma_{j},\kappa_{j}) $.

(iii) If $(E_{j},\rho_{j},\sigma_{j},\kappa_{j})$, $j \in \{1,\ldots,L\}$, are as above but not necessarily pairwise orthogonal, and if $\nu_{j}$, $j \in \{1,\ldots,L\}$, are positive rationals with $\sum_{j} \nu_{j} = 1$, we may choose a unital $^{*}$-homomorphism 
\[
\textstyle
\gamma : \bigoplus_{j=1}^{L} \mathcal{Q} \to \mathcal{Q}
\]
with $\tau_{\mathcal{Q}}(\gamma(1_{j}))= \nu_{j}$; then, the $(E_{j},\rho_{j},\gamma_{j} \circ \sigma_{j},\kappa_{j})$ are pairwise orthogonal and we write 
\[
\textstyle
\bigoplus_{j} \nu_{j} \cdot (E_{j},\rho_{j},\sigma_{j},\kappa_{j})
\]
for
\[
\textstyle
\bigoplus_{j} \nu_{j} \cdot (E_{j},\rho_{j},\gamma_{j} \circ \sigma_{j},\kappa_{j}).
\]
This is well-defined since $\bigoplus_{j} \gamma_{j} \circ \sigma_{j}$ only depends on the choice of $\gamma$ up to unitary equivalence.
\end{ndefn}
\en

\bn
\begin{nprop}
\label{I}
Let $B$ be a unital recursive subhomogeneous $\mathrm{C}^{*}$-algebra with recursive subhomogeneous decomposition 
\[
[B_{l},X_{l},\Omega_{l},r_{l},\phi_{l}]_{l=1}^{R};
\]
let $\mathcal{F} \subset B_{+}^{1}$ finite and $\eta>0$ be given.  Let $(E_{j},\rho_{j},\sigma_{j},\kappa_{j})$, $j \in \{1,\ldots,L\}$, be weighted $(\mathcal{F},\eta)$-excisors.

Then there are pairwise orthogonal weighted $(\mathcal{F},\eta)$-excisors $(E_{j},\rho_{j},\dot{\sigma}_{j},\kappa_{j})$, $j \in \{1,\ldots,L\}$, such that, for each $j$, 
\begin{equation}
\label{I1}
(E_{j},\rho_{j},\sigma_{j},\kappa_{j}) \sim_{(\mathcal{F},\eta)} (E_{j},\rho_{j},\dot{\sigma}_{j},\kappa_{j}).
\end{equation}

If  the $(E_{j},\rho_{j},\sigma_{j},\kappa_{j})$ are compatible with $[B_{l},X_{l},\Omega_{l},r_{l},\phi_{l}]_{l=1}^{R}$, we may choose the  $(E_{j},\rho_{j},\dot{\sigma}_{j},\kappa_{j})$ and the $(\mathcal{F},\eta)$-bridges to be compatible as well.
\end{nprop}

\begin{nproof}
Choose pairwise orthogonal nonzero projections $q_{j} \in \mathcal{Q}$, $j \in \{1,\ldots,L\}$, and isomorphisms  
\[
\theta_{j}:\mathcal{Q} \to q_{j} \mathcal{Q} q_{j};
\]
set
\[
\dot{\sigma}_{j}:= (\mathrm{id} \otimes \theta_{j}) \circ \sigma_{j}.
\]
It is clear that the $(E_{j},\rho_{j},\dot{\sigma}_{j},\kappa_{j})$ are pairwise orthogonal weighted $(\mathcal{F},\eta)$-excisors. Now \eqref{I1} holds, in fact with an $(\mathcal{F},\eta)$-bridge of length $K=1$, since passing from $\sigma_{j}$ to $\dot{\sigma}_{j}$ does not affect \eqref{L1}. Also, changing the $\sigma_{j}$ does not affect compatibility with the recursive subhomogeneous decomposition. 
\end{nproof}
\en

\bn
\begin{nprop}
\label{H}
Let $B$ be a unital recursive subhomogeneous $\mathrm{C}^{*}$-algebra with recursive subhomogeneous decomposition 
\[
[B_{l},X_{l},\Omega_{l},r_{l},\phi_{l}]_{l=1}^{R};
\]
let $\mathcal{F} \subset B_{+}^{1}$ finite and $\eta>0$ be given.  Let $(E_{j},\rho_{j},\sigma_{j},\kappa_{j})$, $j \in \{1,\ldots,L\}$, be pairwise orthogonal weighted $(\mathcal{F},\eta)$-excisors. Let $(E'_{j},\rho'_{j},\sigma'_{j},\kappa'_{j})$, $j \in \{1,\ldots,L\}$, be another set of pairwise orthogonal weighted $(\mathcal{F},\eta)$-excisors, and let 
\[
\textstyle
\gamma: \bigoplus_{j=1}^{L} \mathcal{Q} \to \mathcal{Q}
\]
be a unital $^{*}$-homomorphism. Suppose that
\[
(E_{j},\rho_{j},\sigma_{j},\kappa_{j}) \sim_{(\mathcal{F},\eta)} (E'_{j},\rho'_{j},\sigma'_{j},\kappa'_{j})
\]
for each $j \in \{1,\ldots,L\}$.

Then
\[
\textstyle
\bigoplus_{\gamma} (E_{j}, \rho_{j},\sigma_{j},\kappa_{j})  \sim_{(\mathcal{F},\eta)} \bigoplus_{\gamma} (E'_{j}, \rho'_{j},\sigma'_{j},\kappa'_{j}) .
\]

If  the $(E_{j},\rho_{j},\sigma_{j},\kappa_{j})$ and the  $(E'_{j},\rho'_{j},\sigma'_{j},\kappa'_{j})$ are compatible with the decomposition $[B_{l},X_{l},\Omega_{l},r_{l},\phi_{l}]_{l=1}^{R}$, we may choose the  $(\mathcal{F},\eta)$-bridge between the $\gamma$-direct sums to be compatible as well.
\end{nprop}

\begin{nproof}
For each $j \in \{1,\ldots,L\}$ choose an $(\mathcal{F},\eta)$-bridge between $(E_{j},\rho_{j},\sigma_{j},\kappa_{j})$ and $(E'_{j},\rho'_{j},\sigma'_{j},\kappa'_{j})$; by repeating some of the steps, if necessary, we may assume that all of these have the same length, say $K$, and are given by weighted $(\mathcal{F},\eta)$-excisors  $(E_{j,\frac{i}{K}},\rho_{j,\frac{i}{K}},\sigma_{j,\frac{i}{K}},\kappa_{j,\frac{i}{K}})$ with 
\[
(E_{j,0},\rho_{j,0},\sigma_{j,0},\kappa_{j,0}) = (E_{j},\rho_{j},\sigma_{j},\kappa_{j})
\]
and 
\[
(E_{j,1},\rho_{j,1},\sigma_{j,1},\kappa_{j,1}) = (E'_{j},\rho'_{j},\sigma'_{j},\kappa'_{j}).
\]

As in the proof of Proposition~\ref{I}, choose pairwise orthogonal nonzero projections $q_{j} \in \mathcal{Q}$, $j \in \{1,\ldots,L\}$, as well as isomorphisms 
\[
\theta_{j}:\mathcal{Q} \to q_{j} \mathcal{Q} q_{j}.
\]
Set
\[
\dot{\sigma}_{j,\frac{i}{K}}:= (\mathrm{id} \otimes \theta_{j}) \circ \sigma_{j,\frac{i}{K}},
\]
then the sums
\[
\textstyle
\bigoplus_{\gamma} (E_{j,\frac{i}{K}},\rho_{j,\frac{i}{K}},\dot{\sigma}_{j,\frac{i}{K}},\kappa_{j,\frac{i}{K}})
\]
are $(\mathcal{F},\eta)$-excisors implementing an $(\mathcal{F},\eta)$-bridge 
\[
\textstyle
\bigoplus_{\gamma} (E_{j,0},\rho_{j,0},\dot{\sigma}_{j,0},\kappa_{j,0}) \sim_{(\mathcal{F},\eta)} \bigoplus_{\gamma} (E_{j,1},\rho_{j,1},\dot{\sigma}_{j,1},\kappa_{j,1}). 
\]
As in the proof of \ref{I}, it remains to observe that
\[
\textstyle
\bigoplus_{\gamma} (E_{j},\rho_{j},\sigma_{j},\kappa_{j}) \sim_{(\mathcal{F},\eta)} \bigoplus_{\gamma} (E_{j,0},\rho_{j,0},\dot{\sigma}_{j,0},\kappa_{j,0})
\]
and
\[
\textstyle
\bigoplus_{\gamma} (E'_{j},\rho'_{j},\sigma'_{j},\kappa'_{j}) \sim_{(\mathcal{F},\eta)} \bigoplus_{\gamma} (E_{j,1},\rho_{j,1},\dot{\sigma}_{j,1},\kappa_{j,1}). 
\]
\end{nproof}
\en

\bn
\begin{ndefn}
\label{Fd}
Let $B$ be a unital recursive subhomogeneous $\mathrm{C}^{*}$-algebra with recursive subhomogeneous decomposition 
\[
[B_{l},X_{l},\Omega_{l},r_{l},\phi_{l}]_{l=1}^{R};
\]
let $\mathcal{F} \subset B_{+}^{1}$ finite and $\eta>0$ be given.  

If $(E_{j},\rho_{j},\sigma_{j},\kappa_{j})$ and $(E_{j},\rho_{j},\dot{\sigma}_{j},\kappa_{j})$, $j \in \{1,\ldots,L\}$, are as in Proposition~\ref{I}, and if 
\[
\textstyle
\gamma: \bigoplus_{j=1}^{L} \mathcal{Q} \to \mathcal{Q}
\]
is a unital $^{*}$-homomorphism, we say
\begin{equation}
\label{F1}
\textstyle
\bigoplus_{\gamma} (E_{j},\rho_{j},\dot{\sigma}_{j}, \kappa_{j})
\end{equation}
is a compatible $\gamma$-direct sum of the $(E_{j},\rho_{j},\sigma_{j},\kappa_{j})$.
\end{ndefn}
\en

\bn
\begin{nremark}
Of course, the $\gamma$-direct sum in \eqref{F1} depends on the choice of the $\dot{\sigma}_{j}$ in Proposition~\ref{I}, but for a different choice, say $\ddot{\sigma}_{j}$, it follows from Proposition~\ref{H} that 
\[
\textstyle
\bigoplus_{\gamma}(E_{j},\rho_{j},\dot{\sigma}_{j},\kappa_{j})\sim_{(\mathcal{F},\eta)} \bigoplus_{\gamma}(E_{j},\rho_{j},\ddot{\sigma}_{j},\kappa_{j}).
\]
\end{nremark}
\en

\bn
\begin{nprop}
\label{F}
Let $B$ be a unital recursive subhomogeneous $\mathrm{C}^{*}$-algebra with recursive subhomogeneous decomposition 
\[
[B_{l},X_{l},\Omega_{l},r_{l},\phi_{l}]_{l=1}^{R};
\]
let $\mathcal{F} \subset B_{+}^{1}$ finite and $\eta>0$ be given.  

Let 
\[
\textstyle
\left(E=\bigoplus_{j=1}^{L} E_{j},\rho,\sigma,\kappa \right)
\]
be an $(\mathcal{F},\eta)$-excisor and let
\[
\textstyle
\gamma: \bigoplus_{j=1}^{L} \mathcal{Q} \to \mathcal{Q}
\]
be a unital $^{*}$-homomorphism such that
\[
\tau_{\mathcal{Q}}(\gamma_{j}(1_{\mathcal{Q}})) = \tau_{\mathcal{Q}}(\kappa(1_{E_{j}})) \mbox{ for } j \in \{1,\ldots,L\}.
\]

Then there are pairwise orthogonal $(\mathcal{F},\eta)$-excisors
\[
(E_{j}, \rho_{j} = \rho|_{E_{j}}, \dot{\sigma}_{j}, \dot{\kappa}_{j}:E_{j}  \to \kappa(1_{E_{j}}) \mathcal{Q}  \kappa(1_{E_{j}}) \cong \mathcal{Q} )
\]
such that 
\begin{equation}
\label{F1a}
\textstyle
\kappa \approx_{\mathrm{u}} \gamma \circ (\bigoplus_{j=1}^{L} \dot{\kappa}_{j}) 
\end{equation}
and such that 
\begin{equation}
\label{F2}
\textstyle
(E,\rho,\sigma,\kappa) \sim_{(\mathcal{F},\eta)} \bigoplus_{\gamma}(E_{j},\rho_{j},\dot{\sigma}_{j},\dot{\kappa}_{j}).
\end{equation}

If   $(E,\rho,\sigma,\kappa)$ is compatible with $[B_{l},X_{l},\Omega_{l},r_{l},\phi_{l}]_{l=1}^{R}$, we may choose the  $\gamma$-direct sum and the $(\mathcal{F},\eta)$-bridge to be compatible as well.
\end{nprop}

\begin{nproof}
Let
\[
\zeta_{j}: \kappa(1_{E_{j}}) \mathcal{Q} \kappa(1_{E_{j}}) \to \mathcal{Q}
\]
be an isomorphism for each $j \in \{1,\ldots,L\}$, then the maps 
\[
\dot{\kappa}_{j}:= \zeta_{j} \circ \kappa|_{E_{j}} 
\] 
clearly satisfy \eqref{F1a}. Choose pairwise nonzero orthogonal projections $q_{j} \in \mathcal{Q}$, $j \in \{1,\ldots,K\}$, as well as isomorphisms 
\[
\theta_{j}: \mathcal{Q} \to q_{j} \mathcal{Q} q_{j};
\]
define 
\[
\dot{\sigma}_{j}:= (\mathrm{id}_{B} \otimes \theta_{j}) \circ \sigma|_{E_{j}}.
\]
It is then clear that the $(E_{j},\rho_{j},\dot{\sigma}_{j},\dot{\kappa}_{j})$ are pairwise orthogonal $(\mathcal{F},\eta)$-excisors and that 
\[
\textstyle
 \bigoplus_{\gamma}(E_{j},\rho_{j},\dot{\sigma}_{j},\dot{\kappa}_{j})  \sim_{(\mathcal{F},\eta)}   (E,\rho, \sigma, \gamma \circ (\bigoplus_{j=1}^{L} \dot{\kappa}_{j})).
\]
Finally, \eqref{F2} follows from \eqref{F1a} and Remark~\ref{M}(ii).
\end{nproof}
\en

\bn
We note the following lifting result,  which will imply the existence of sufficiently many $(\mathcal{F},\eta)$-excisors, cf.\ Remark~\ref{Fb}(ii) below. 

\begin{nprop}
\label{Fa}
Let $B$, $F$ be $\mathrm{C}^{*}$-algebras, $F$ finite dimensional, and $\pi:B \to F$ a surjective $^{*}$-homomorphism; let $\mathcal{F} \subset B_{+}^{1}$ finite and $\eta>0$ be given.  

Then there is an isometric c.p.\ order zero map 
\[
\sigma:F \to B
\]
such that
\begin{equation}
\label{Fa1}
\|\sigma(1_{F})b - \sigma \pi(b)\| < \eta \mbox{ for } b \in \mathcal{F}.
\end{equation}
\end{nprop}

\begin{nproof}
Since $F$ and $\mathcal{F}$ are separable, we may clearly also assume $B$ to be separable, hence $\sigma$-unital. Recall that c.p.c.\ order zero maps are projective, whence there is a c.p.\ isometric order zero lift 
\[
\dot{\sigma}:F \to B.
\]
Choose an approximate unit $(h_{n})_{n \in \mathbb{N}}$ for $\mathrm{ker}\, \pi$ which is quasicentral for $B$. Define c.p.c.\ maps
\[
\ddot{\sigma}_{n}:F \to B
\]
by
\[
\ddot{\sigma}_{n}(\, .\,):= (1_{B^{\sim}} - h_{n})^{\frac{1}{2}} \dot{\sigma}(\, .\,) (1_{B^{\sim}} -h_{n})^{\frac{1}{2}}
\]
(here, $B^{\sim}$ denotes the smallest unitization of $B$). The $\ddot{\sigma}_{n}$ clearly induce a c.p.\ isometric  order zero map
\[
\textstyle
\ddot{\sigma}:F \to B_{\infty} \prod_{\mathbb{N}} B/ \bigoplus_{\mathbb{N}} B
\]
which in turn lifts to a c.p.\ isometric order zero map 
\[
\textstyle
\bar{\sigma}: F \to \prod_{\mathbb{N}} B
\]
with components $\bar{\sigma}_{n}$. Upon dropping finitely many components and rescaling, if necessary, we may assume each $\bar{\sigma}_{n}$ to be isometric. It is now straightforward to check that, for large enough $N$, $\sigma:= \bar{\sigma}_{N}$ will satisfy \eqref{Fa1}. 
\end{nproof}
\en

\bn
\begin{nnotation}
\label{Fc}
Let $B$ be a unital recursive subhomogeneous $\mathrm{C}^{*}$-algebra with recursive subhomogeneous decomposition 
\[
[B_{l},X_{l},\Omega_{l},r_{l},\phi_{l}]_{l=1}^{R}.
\]
If $l \in \{1,\ldots,R\}$ and $x \in X_{l}$, then 
\[
(\mathrm{ev}_{x} \otimes \mathrm{id}_{M_{r_{l}}}) \circ \iota_{B}: B \to M_{r_{l}}
\]
factorizes through a sum of irreducible representations, say
\[
B \stackrel{\rho_{x}}{\longrightarrow} E_{x} \stackrel{\iota_{E_{x}}}{\longrightarrow} M_{r_{l}}.
\]
Upon fixing a unital embedding
\[
M_{r_{l}} \to \mathcal{Q}
\]
we obtain unital $^{*}$-homomorphisms
\[
B \stackrel{\rho_{x}}{\longrightarrow} E_{x} \stackrel{\kappa_{x}}{\longrightarrow} \mathcal{Q}
\]
such that $\rho_{x}$ is a sum of surjective irreducible representations and 
\[
\tau_{\mathcal{Q}} \kappa_{x} = \tau_{M_{r_{l}}} \iota_{E_{x}}.
\]
\end{nnotation}
\en

\bn
\begin{nremarks}
\label{Fb}
(i) The maps $\rho_{x}$ and $\kappa_{x}$ are uniquely determined by $x$ up to unitary equivalence.

\noindent
(ii) By Proposition~\ref{Fa}, for any finite subset $\mathcal{F} \subset B^{1}_{+}$ and $\eta>0$, and for any $x \in X_{l}$, there is an isometric c.p.\ order zero map 
\[
\sigma_{x}:E_{x} \to B
\]
such that $(E_{x},\rho_{x},\sigma_{x},\kappa_{x})$ is a weighted $(\mathcal{F},\eta)$-excisor (which is compatible with the decomposition $[B_{l},X_{l},\Omega_{l},r_{l},\phi_{l}]_{l=1}^{R}$, provided $x \in X_{l} \setminus \Omega_{l}$).

\noindent
(iii) It is clear that, if $x,x' \in X_{l}$ are such that 
\[
\|(\mathrm{ev}_{x} \otimes \mathrm{id}_{M_{r_{l}}}) \circ \iota_{B}(b) - (\mathrm{ev}_{x'} \otimes \mathrm{id}_{M_{r_{l}}}) \circ \iota_{B}(b) \| \le \eta
\]
for all $b \in \mathcal{F}$, then 
\[
(E_{x},\rho_{x},\sigma_{x},\kappa_{x}) \sim_{(\mathcal{F},\eta)} (E_{x'},\rho_{x'},\sigma_{x'},\kappa_{x'}),
\]
in fact via an $(\mathcal{F},\eta)$-bridge of length $K=1$.
\end{nremarks}
\en

\bn
\begin{nprop}
\label{G}
Let $B$ be a unital recursive subhomogeneous $\mathrm{C}^{*}$-algebra with recursive subhomogeneous decomposition 
\[
[B_{l},X_{l},\Omega_{l},r_{l},\phi_{l}]_{l=1}^{R};
\]
let $\mathcal{F} \subset B_{+}^{1}$ finite and $\eta>0$ be given.  

Fix $l \in \{1,\ldots, R\}$ and $x \in X_{l} \setminus \Omega_{l}$, and let $(E_{x},\rho_{x},\sigma_{x}, \kappa_{x})$ be an $(\mathcal{F},\eta)$-excisor, with $(E_{x}$, $\rho_{x}$, $\kappa_{x})$ as in \ref{Fc} (note that $E_{x} \cong M_{r_{l}}$ since $x \in X_{l} \setminus \Omega_{l}$). Let
\[
\textstyle
\gamma: \bigoplus_{j=1}^{L} \mathcal{Q} \to \mathcal{Q}
\] 
be a unital embedding for some $L \in \mathbb{N}$.   
 
Then there are pairwise orthogonal $(\mathcal{F},\eta)$-excisors
\[
(E_{x},\rho_{x},\sigma_{x,j}, \kappa_{x}), \; j \in \{1,\ldots,L\},
\]
such that
\begin{equation}
\label{G1}
\textstyle
(E_{x},\rho_{x},\sigma_{x}, \kappa_{x})
\sim_{(\mathcal{F},\eta)}
\bigoplus_{\gamma} (E_{x},\rho_{x},\sigma_{x,j}, \kappa_{x}) .
\end{equation}
\end{nprop}

\begin{nproof}
Choose $q_{j}$ and $\theta_{j}$ as in the proof of Proposition~\ref{F}. We take
\[
\sigma_{x,j}:= (\mathrm{id}_{B} \otimes \theta_{j}) \circ \sigma_{x}
\] 
(these correspond to the maps $\dot{\sigma}_{x}$ from Proposition~\ref{F}); as in the proof of \ref{F} one checks that 
\[
\textstyle
\bigoplus_{\gamma} (E_{x},\rho_{x},\sigma_{x,j}, \kappa_{x})  \sim_{(\mathcal{F},\eta)} (E_{x},\rho_{x},\sigma_{x}, \gamma \circ (\kappa_{x}^{\oplus L})).
\]
Now observe that $\kappa_{x} \approx_{\mathrm{u}} \gamma \circ(\kappa_{x}^{\oplus L})$, and apply Remark~\ref{M}(ii) to obtain \eqref{G1}.
\end{nproof}
\en

%%%%%%%%%%%%%%%%%%%%%%%%%%%%%%%%%%%%%%%%%%%%%%%%
\section{$(\mathcal{F},\eta)$-connected decompositions} \label{connected}

\bn
\begin{ndefn}
Let $B$ be a unital recursive subhomogeneous $\mathrm{C}^{*}$-algebra, and  let $\mathcal{F} \subset B_{+}^{1}$ finite and $\eta>0$ be given.  

A recursive subhomogeneous decomposition 
\[
[B_{l},X_{l},\Omega_{l},r_{l},\phi_{l}]_{l=1}^{R}
\]
for $B$ is $(\mathcal{F},\eta)$-connected if the following holds: 

If $l \in \{1,\ldots,R\}$ and $x,y \in X_{l}$, and if $(E_{x},\rho_{x},\sigma_{x},\kappa_{x})$ and  $(E_{y},\rho_{y},\sigma_{y},\kappa_{y})$ are $(\mathcal{F},\eta)$-excisors with $(E_{x}$, $\rho_{x}$, $\kappa_{x})$ and  $(E_{y}$, $\rho_{y}$, $\kappa_{y})$ as in \ref{Fc}, then 
\[
(E_{x},\rho_{x},\sigma_{x},\kappa_{x})   \sim_{(\mathcal{F},\eta)}(E_{y},\rho_{y},\sigma_{y},\kappa_{y}).
\]
\end{ndefn}
\en

\bn
\begin{nprop}
\label{B}
Let $B$ be a unital recursive subhomogeneous $\mathrm{C}^{*}$-algebra and  let $\mathcal{F} \subset B_{+}^{1}$ finite and $\eta>0$ be given.  

Then $B$ has an $(\mathcal{F},\eta)$-connected recursive subhomogeneous decomposition 
\[
[B_{l},X_{l},\Omega_{l},r_{l},\phi_{l}]_{l=1}^{R}.
\]
If $\mathrm{dr} B \le n$, then $X_{1}, \ldots, X_{R}$ may be chosen so that $\dim X_{l} \le n$ for $l \in \{1,\ldots,R\}$. 
\end{nprop}

\begin{nproof}
This follows immediately from  Remark~\ref{Fb}(iii) after decomposing each $X_{l}$ in to pairwise disjoint closed subsets 
\[
\textstyle
X_{l} = \coprod_{k=1}^{N_{l}} X_{l,k}
\]
such that, for each $l\in \{1,\ldots,R\}$, $k \in \{1,\ldots,N_{l}\}$ and $x_{0},x_{1} \in X_{l,k}$, there are $K \in \mathbb{N}$ and $x_{\frac{j}{K}} \in X_{l,k}$ such that 
\[
\|(\mathrm{ev}_{x_{\frac{j}{K}}} \otimes \mathrm{id}_{M_{r_{l}}}) \circ \iota_{B}(b) - (\mathrm{ev}_{x_{\frac{j+1}{K}}} \otimes \mathrm{id}_{M_{r_{l}}}) \circ \iota_{B}(b) \| \le \eta
\]
for all $j \in \{0,\ldots,K-1\}$ and $b \in \mathcal{F}$.

The last statement follows from \cite{Winter:subhomdr}, since in this case the $X_{l}$ (and hence the $X_{l,k}$) may be chosen to have dimension at most $n$.
\end{nproof}
\en

\bn
\begin{nprop}
\label{D}
Let $B$ be a unital recursive subhomogeneous $\mathrm{C}^{*}$-algebra; let $\mathcal{F} \subset B_{+}^{1}$ finite and $\eta>0$ be given and suppose  
\[
[B_{l},X_{l},\Omega_{l},r_{l},\phi_{l}]_{l=1}^{R}
\]
is an $(\mathcal{F},\eta)$-connected  recursive subhomogeneous decomposition for $B$. 

Let
\[
\textstyle
\left( E_{i}= \bigoplus_{l=1}^{R} E_{i,l}, \, \rho_{i} = \bigoplus_{l=1}^{R} \rho_{i,l}, \, \sigma_{i} = \bigoplus_{l=1}^{R} \sigma_{i,l}, \, \kappa_{i}  \right), \; i \in \{0,1\}
\]
be weighted $(\mathcal{F},\eta)$-excisors   (compatible with the decomposition) satisfying
\[
y_{l}:= \tau_{\mathcal{Q}}(\kappa_{0}(1_{E_{0,l}})) = \tau_{\mathcal{Q}}(\kappa_{1}(1_{E_{1,l}})) , \; l \in \{1,\ldots,R\}
\]
and such that each $\rho_{i,l}$ factorizes as 
\[
\rho_{i,l}: B \to \mathcal{C}(X_{l}) \otimes M_{r_{l}} \to E_{i,l}.
\]

Then (via a compatible $(\mathcal{F},\eta)$-bridge),
\[
(E_{0},\rho_{0},\sigma_{0},\kappa_{0}) \sim_{(\mathcal{F},\eta)}  (E_{1},\rho_{1},\sigma_{1},\kappa_{1}).
\]
\end{nprop}

\begin{nproof}
For each $l \in \{1,\ldots,R\}$, choose $x_{l} \in X_{l} \setminus \Omega_{l}$. Set 
\[
\textstyle
E:= \bigoplus_{l=1}^{R} E_{x_{l}}, \; \rho:= \bigoplus_{l=1}^{R} \rho_{x_{l}}, \; \kappa:= \gamma \circ \left(  \bigoplus_{l=1}^{R} \kappa_{x_{l}} \right), 
\]
where $E_{x_{l}}$, $\rho_{x_{l}}$, $\kappa_{x_{l}}$ are as in \ref{Fc} and 
\[
\textstyle
\gamma: \bigoplus_{l=1}^{R} \mathcal{Q} \to \mathcal{Q} 
\]
is a unital embedding such that 
\[
\tau_{\mathcal{Q}}(\gamma_{l}(1_{\mathcal{Q}})) = y_{l}, \; l \in \{1,\ldots,R\}.
\]
Note that $E_{x_{l}} \cong M_{r_{l}}$, since $x_{l} \in X_{l} \setminus \Omega_{l}$ by Remark~\ref{r-dK}.

By Proposition~\ref{Fa}, there is an isometric c.p.\ order zero map 
\[
\sigma:E \to B \otimes \mathcal{Q}
\]
such that $(E,\rho,\sigma,\kappa)$ is a weighted $(\mathcal{F},\eta)$-excisor which is compatible with the decomposition $[B_{l},X_{l},\Omega_{l},r_{l},\phi_{l}]_{l=1}^{R}$. By Proposition~\ref{F}, there are (compatible) pairwise orthogonal $(\mathcal{F},\eta)$-excisors of the form 
\[
(E_{x_{l}},\rho_{x_{l}}, \dot{\sigma}_{x_{l}}, \kappa_{x_{l}}), \; l \in \{1,\ldots,R\},
\]
such that (in a compatible way)
\begin{equation}
\label{Dw3}
\textstyle
(E,\rho,\sigma,\kappa) \sim_{(\mathcal{F},\eta)} \bigoplus_{\gamma} (E_{x_{l}},\rho_{x_{l}},\dot{\sigma}_{x_{l}},\kappa_{x_{l}}).
\end{equation}
Similarly, for $i \in \{0,1\}$ and $l \in \{1,\ldots,R\}$ there are (compatible) pairwise orthogonal $(\mathcal{F},\eta)$-excisors  
\[
(E_{i,l},\rho_{i,l}, \dot{\sigma}_{i,l}, \dot{\kappa}_{i,l})
\]
such that 
\[
\textstyle
\kappa_{i} \approx_{\mathrm{u}} \gamma \circ \left( \bigoplus_{l=1}^{R} \dot{\kappa}_{i,l}  \right)
\]
and such that (in a compatible way)
\begin{equation}
\label{Dw4}
\textstyle
(E_{i},\rho_{i}, \sigma_{i}, \kappa_{i}) \sim_{(\mathcal{F},\eta)}
 \bigoplus_{\gamma} (E_{i,l},\rho_{i,l}, \dot{\sigma}_{i,l}, \dot{\kappa}_{i,l}).
\end{equation}

Since $(E_{i,l},\rho_{i,l}, \dot{\sigma}_{i,l}, \dot{\kappa}_{i,l})$ is compatible with $[B_{l},X_{l},\Omega_{l},r_{l},\phi_{l}]_{l=1}^{R}$, there are $N_{i,l} \in \mathbb{N}$ and $x_{i,l,m} \in X_{l} \setminus \Omega_{l}$ for $m \in \{1,\ldots,N_{i,l}\}$ such that 
\[
\textstyle
(E_{i,l},\rho_{i,l}, \dot{\sigma}_{i,l},\dot{\kappa}_{i,l}) = \left( \bigoplus_{m=1}^{N_{i,l}} E_{x_{i,l,m}}, \;  \bigoplus_{m=1}^{N_{i,l}} \rho_{x_{i,l,m}}, \; \bigoplus_{m=1}^{N_{i,l}} \dot{\sigma}_{x_{i,l,m}}, \; \bigoplus_{m=1}^{N_{i,l}} \dot{\kappa}_{x_{i,l,m}}   \right);
\] 
note that 
\[
E_{x_{i,l,m}} \cong M_{r_{l}}
\]
for all $i,l,m$. 

Let
\[
\textstyle
\gamma_{i,l}: \bigoplus_{m=1}^{N_{i,l}}  \mathcal{Q} \to \mathcal{Q} 
\]
be a unital embedding such that 
\[
\tau_{\mathcal{Q}}(\gamma_{i,l,m}(1_{\mathcal{Q}})) = \tau_{\mathcal{Q}}(\dot{\kappa}_{i,l}(1_{E_{i,l,m}})), \; m \in \{1,\ldots,N_{i,l}\}.
\]
By Proposition~\ref{F}, there are pairwise orthogonal $(\mathcal{F},\eta)$-excisors
\[
(E_{x_{i,l,m}}, \rho_{x_{i,l,m}}, \ddot{\sigma}_{x_{i,l,m}}, \ddot{\kappa}_{x_{i,l,m}}), \; m \in \{1,\ldots,N_{i,l}\},
\]
such that 
\[
\textstyle
\dot{\kappa}_{i,l} \approx_{\mathrm{u}} \gamma_{i,l} \circ \left( \bigoplus_{m=1}^{N_{i,l}} \ddot{\kappa}_{i,l,m}  \right)
\]
and such that 
\begin{equation}
\label{Dw1}
\textstyle
(E_{i,l},\rho_{i,l},\dot{\sigma}_{i,l}, \dot{\kappa}_{i,l}) \sim_{(\mathcal{F},\eta)}  \bigoplus_{\gamma_{i,l}} (E_{x_{i,l,m}}, \rho_{x_{i,l,m}}, \ddot{\sigma}_{x_{i,l,m}}, \ddot{\kappa}_{x_{i,l,m}}    )  . 
\end{equation}
By Proposition~\ref{G}, there are pairwise orthogonal $(\mathcal{F},\eta)$-excisors 
\[
(E_{x_{l}}, \rho_{x_{l}}, \dot{\sigma}_{x_{l},m}, \kappa_{x_{l}}), \; m \in \{1,\ldots,N_{i,l}\},
\]
such that
\begin{equation}
\label{Dw2}
\textstyle
(E_{x_{l}},\rho_{x_{l}}, \dot{\sigma}_{x_{l}},\kappa_{x_{l}}) \sim_{(\mathcal{F},\eta)} \bigoplus_{\gamma_{i,l}} (E_{x_{l}}, \rho_{x_{l}}, \dot{\sigma}_{x_{l},m}, \kappa_{x_{l}}).
\end{equation}
Since $[B_{l},X_{l},\Omega_{l},r_{l},\phi_{l}]_{l=1}^{R}$ is $(\mathcal{F},\eta)$-connected, for each $i,l,m$ we have
\[
(E_{x_{l}},\rho_{x_{l}},\dot{\sigma}_{x_{l},m}, \kappa_{x_{l}}) \sim_{(\mathcal{F},\eta)} (E_{x_{i,l,m}}, \rho_{x_{i,l,m}}, \ddot{\sigma}_{x_{i,l,m}}, \ddot{\kappa}_{x_{i,l,m}}). 
\]
By Proposition~\ref{H}, we have
\[
\textstyle
\bigoplus_{\gamma_{i,l}} (E_{x_{l}},\rho_{x_{l}}, \dot{\sigma}_{x_{l},m}, \kappa_{x_{l}}) \sim_{(\mathcal{F},\eta)}  \bigoplus_{\gamma_{i,l}} (E_{x_{i,l,m}}, \rho_{x_{i,l,m}}, \ddot{\sigma}_{x_{i,l,m}}, \ddot{\kappa}_{x_{i,l,m}}    )  
\]
which in turn yields
\[
(E_{x_{l}},\rho_{x_{l}}, \dot{\sigma}_{x_{l}},\kappa_{x_{l}})     \sim_{(\mathcal{F},\eta)}     (E_{i,l},\rho_{i,l},\dot{\sigma}_{i,l}, \dot{\kappa}_{i,l})
\]
by \eqref{Dw1} and \eqref{Dw2}.

Again by Proposition~\ref{H}, together with \eqref{Dw3} and \eqref{Dw4} this gives 
\[
(E,\rho,\sigma,\kappa) \sim_{(\mathcal{F},\eta)} (E_{i}, \rho_{i},\sigma_{i},\kappa_{i}), \; i \in \{0,1\},
\]
from which we obtain
\[
(E_{0},\rho_{0},\sigma_{0},\kappa_{0}) \sim_{(\mathcal{F},\eta)} (E_{1},\rho_{1},\sigma_{1},\kappa_{1}),
\]
as desired. Of course all the $(\mathcal{F},\eta)$-bridges above may be chosen to be compatible with the given recursive subhomogeneous decomposition.
\end{nproof}
\en

%%%%%%%%%%%%%%%%%%%%%%%%%%%%%%%%%%%%%%%%%%%%%%%%
\section{Excising traces} \label{excising traces}

\bn
\begin{nnotation}
\label{local-trace-notation}
Let $B$ be a separable unital recursive subhomogeneous $\mathrm{C}^{*}$-algebra with (separable) recursive subhomogeneous decomposition 
\[
[B_{l},X_{l},\Omega_{l},r_{l},\phi_{l}]_{l=1}^{R};
\]
let $\tau \in T(B)$ be a tracial state. We inductively define positive tracial functionals 
\[
\tau_{l},\bar{\tau}_{l}:B_{l} \to \mathbb{C}, \; l \in \{1,\ldots,R\}
\]
as follows:

For each $l$, let $0\le h_{l} \le 1$ be a strictly positive element of $\mathcal{C}_{0}(X_{l} \setminus \Omega_{l})$. Set
\[
\tau_{R} := \tau: B \cong B_{R} \to \mathbb{C}.
\]

If $\tau_{l}:B_{l} \to \mathbb{C}$ has been constructed, set 
\[
\bar{\tau}_{l}(b):= \lim_{n \to \infty} \tau_{l} ((h_{l}^{\frac{1}{n}} \otimes 1_{M_{r_{l}}}) b), \; b \in B_{l}.
\]
(On positive elements $b$, the limit is over a bounded increasing sequence, hence exists; but then the limit also exists for general $b$). 

If $\tau_{l}, \bar{\tau}_{l}$ have been constructed, set
\[
\tau_{l-1}(b)  = \tau_{l}(\hat{b}) - \bar{\tau}_{l}(\hat{b}), \; b \in B_{l-1},
\]
where $\hat{b} \in B_{l}$ is a lift for $b$. It is easy to see that $\tau_{l}, \bar{\tau}_{l}$, $l \in \{1,\ldots,R\}$, are well-defined positive functionals which do not depend on the choice of the $h_{l}$, that $\bar{\tau}_{l} \le \tau_{l}$, that
\begin{equation}
\label{local-trace-notationw1}
y^{\tau}_{l} := \tau_{l}(1_{B_{l}}) - \tau_{l-1}(1_{B_{l-1}}) = \|\bar{\tau}_{l}\| (\le 1)
\end{equation}
and that
\begin{equation}
\label{local-trace-notationw2}
\textstyle
\sum_{l=1}^{R} y^{\tau}_{l} = 1.
\end{equation}
Call the $y^{\tau}_{l}$ the weights of $\tau$ with respect to the decomposition $[B_{l},X_{l},\Omega_{l},r_{l},\phi_{l}]_{l=1}^{R}$.

\bigskip
Now suppose
\[
W \subset X_{l} \setminus \Omega_{l}
\]
is a  subset closed in $X_{l}$. Let 
\[
0 \le g_{l} \le 1
\]
be a strictly positive element for $\mathcal{C}_{0}(X_{l} \setminus W)$, with  $g|_{\Omega_{l}} \equiv 1$. It is not hard to check that, for all $b \in B_{l}$,
\[
\lim_{n \to \infty} \tau_{l}(((1-g_{l}^{\frac{1}{n}})\otimes 1_{M_{r_{l}}}) b)  = \lim_{n \to \infty} \bar{\tau}_{l}(((1-g_{l}^{\frac{1}{n}}) \otimes 1_{M_{r_{l}}}) b)
\]
(and, in particular, that the limits exist). As above, one may define a positive tracial functional 
\[
\tilde{\tau}_{W}: \mathcal{C}(W) \otimes M_{r_{l}} \to \mathbb{C} 
\]
by
\[
\tilde{\tau}_{W}(b) = \lim_{n \to \infty} \tau_{l} ((1-g_{l}^{\frac{1}{n}})\hat{b}),
\]
where $\hat{b} \in B_{l}$ is a lift of $b \in    \mathcal{C}(W) \otimes M_{r_{l}}$. 
\end{nnotation}
\en

\bn
\begin{nprop}
\label{A}
Let $B$ be a separable unital recursive subhomogeneous \mbox{$\mathrm{C}^{*}$-algebra} with (separable) recursive subhomogeneous decomposition 
\[
[B_{l},X_{l},\Omega_{l},r_{l},\phi_{l}]_{l=1}^{R};
\]
let $\tau \in T(B)$ be a tracial state and let $\mathcal{F} \subset B^{1}_{+}$ finite and $\eta>0$ be given. 

Then, for any $\gamma >0$, there is an $(\mathcal{F},\eta)$-excisor
\[
\textstyle
\left( E=\bigoplus_{l=1}^{R} E_{l}, \rho = \bigoplus_{l=1}^{R} \rho_{l}, \sigma = \bigoplus_{l=1}^{R} \sigma_{l}   \right)
\] 
which is compatible with $[B_{l},X_{l},\Omega_{l},r_{l},\phi_{l}]_{l=1}^{R}$ and such that 
\[
(\bar{\tau}_{l} \otimes  \tau_{\mathcal{Q}}) \circ (\psi_{l} \otimes \mathrm{id}_{\mathcal{Q}}) \circ \sigma_{l} (1_{E_{l}}) \ge y^{\tau}_{l} - \gamma, \; l \in \{1,\ldots,R\}.
\]
\end{nprop}

\begin{nproof}
It is straightforward to find, for each $l \in \{1,\ldots, R\}$, an $N_{l} \in \mathbb{N}$ and subsets $W_{l} \subset X_{l} \setminus \Omega_{l} $
satisfying the following:
\begin{itemize}
\item[(i)] Each $W_{l}$ is a disjoint union $W_{l} = \coprod_{n=1}^{N_{l}} W_{l,n}$ of closed subsets $W_{l,n} \subset X_{l}$, each containing a point $w_{l,n} \in W_{l,n}$,
\item[(ii)]  $X_{l} \setminus W_{l}$ is an open neighborhood of $\Omega_{l}$,
\item[(iii)] $y^{\tau}_{l} = \|\bar{\tau}_{l}\| \ge \|\tilde{\tau}_{W_{l}}\| \ge \|\bar{\tau}_{l}\| - \gamma$ (see \ref{local-trace-notation} for notation),
\item[(iv)] for any $b \in \mathcal{F}$, $l\in \{1,\ldots,R\}$, $n \in\{1,\ldots,N_{l}\}$ and $w,w'\in W_{l,n}$, 
\[
\|\mathrm{ev}_{w} \pi_{W_l}(b) - \mathrm{ev}_{w'} \pi_{W_l}(b)\| < \eta/2,
\]
where  
\[
\pi_{W_{l}}: B_{l} \to \mathcal{C}(W_{l}) \otimes M_{r_{l}}
\]
denote the canonical surjections.
\end{itemize}

Define
\begin{equation}
\label{Aw1}
\textstyle
E_{l}:= \bigoplus_{1}^{N_{l}} M_{r_{l}},
\end{equation}
\[
\textstyle
\rho_{l}:= \bigoplus_{n=1}^{N_{l}} \mathrm{ev}_{w_{l,n}} :B \to E_{l}
\]
and 
\[
\textstyle
\tilde{\sigma}_{l} :=  \bigoplus_{n=1}^{N_{l}}1_{W_{l,n}} \otimes \mathrm{id}_{M_{r_{l}}}: E_{l} \to \bigoplus_{n=1}^{N_{l}} \mathcal{C}(W_{l,n}) \otimes M_{r_{l}} \cong \mathcal{C}(W_{l}) \otimes M_{r_{l}}.
\]
Note that 
\[
\textstyle
\tilde{\sigma}:= \bigoplus_{l = 1}^{R} \tilde{\sigma}_{l}: \bigoplus_{l=1}^{R} E_{l} \to \bigoplus_{l=1}^{R} \mathcal{C}(W_{l}) \otimes M_{r_{l}}
\]
is a $^{*}$-homomorphism, hence in particular c.p.\ order zero.

Let 
\[
\textstyle
\pi: \bigoplus_{l=1}^{R}     \pi_{W_{l}} \circ \psi_{l} :B \to \bigoplus_{l=1}^{R}  \mathcal{C}(W_{l}) \otimes M_{r_{l}}
\]

Using projectivity of c.p.c.\ order zero maps together with an approximate unit for $\mathrm{ker}\, \pi \lhd B$ which is quasicentral for $B$, it is not hard to find a c.p.c.\ order zero lift 
\[
\textstyle
\sigma = \bigoplus_{l=1}^{R}  \sigma_{l} : \bigoplus_{l=1}^{R}  E_{l} \to B
\]
with the right properties; the argument is essentially the same as in the proof of Proposition~\ref{Fa}, so we omit the details.
\end{nproof}
\en

%%%%%%%%%%%%%%%%%%%%%%%%%%%%%%%%%%%%%%%%%%%%%%%%
\section{$(\mathcal{F},\eta)$-bridges via linear algebra} \label{linear algebra}

\bn
\begin{nprop}
\label{J}
Let $B$ be a separable unital recursive subhomogeneous $\mathrm{C}^{*}$-algebra and let $\mathcal{F} \subset B^{1}_{+}$ finite and $\eta>0$ be given. Suppose $B$ has an $(\mathcal{F},\eta)$-connected recursive subhomogeneous decomposition 
\[
[B_{l},X_{l},\Omega_{l},r_{l},\phi_{l}]_{l=1}^{R}
\]
along which projections can be lifted and such that $X_{l} \setminus \Omega_{l} \neq \emptyset$ for $l \ge 1$. 

Let $\tau_{0},\tau_{1} \in T(B)$ be tracial states with 
\[
(\tau_{0})_{*} = (\tau_{1})_{*}
\] 
(as states on the ordered $\mathrm{K}_{0}(B)$) and let $0< \bar{\beta} \le 1$ be given. 

Then there are $x_{l} \in X_{l} \setminus \Omega_{l}$ for $l \in \{1,\ldots,R\}$ and pairwise orthogonal weighted $(\mathcal{F},\eta)$-excisors 
\[
(E_{x_{l}}, \pi_{x_{l}}, \sigma_{x_{l}},\kappa_{x_{l}}), \; l \in \{1,\ldots,R\},
\] 
as well as unital embeddings
\[
\textstyle
\gamma_{0},\gamma_{1}, \tilde{\gamma}: \bigoplus_{l=1}^{R} \mathcal{Q} \to \mathcal{Q} 
\]
and
\[
\bar{\gamma}: \mathcal{Q} \oplus \mathcal{Q} \to \mathcal{Q}
\]
such that, for 
\[
\textstyle
E:= \bigoplus_{l=1}^{R} E_{x_{l}}, \; \pi:= \bigoplus_{l=1}^{R} \pi_{x_{l}}, \; \sigma:= \bigoplus_{l=1}^{R} \sigma_{x_{l}}, 
\]
\begin{equation}
\label{Jw22a}
\textstyle
\bar{\kappa}_{i}:= \gamma_{i} \circ \left(\bigoplus_{l=1}^{R} \kappa_{x_{l}} \right),\; \bar{\kappa}:= \tilde{\gamma} \circ \left( \bigoplus_{l=1}^{R} \kappa_{x_{l}}  \right),
\end{equation}
the weighted $(\mathcal{F},\eta)$-excisors $(E,\pi,\sigma, \bar{\kappa}_{i})$, $i \in \{0,1\}$, and $(E,\pi,\sigma,\bar{\kappa})$ satisfy
\begin{equation}
\label{Jw23}
(E,\pi,\sigma,\bar{\gamma} \circ (\bar{\kappa}_{0} \oplus \bar{\kappa}))  \sim_{(\mathcal{F},\eta)}  (E,\pi,\sigma,\bar{\gamma} \circ (\bar{\kappa}_{1} \oplus \bar{\kappa}))
\end{equation}
and such that 
\begin{equation}
\label{Jw12}
\bar{y}_{i,l} := \tau_{\mathcal{Q}} (\gamma_{i} ( 1_{l}   ))
\end{equation}
satisfy 
\begin{equation}
\label{Jw24}
|\bar{y}_{i,l} - y^{\tau_{i}}_{l} | < \bar{\beta}
\end{equation}
for $i \in \{0,1\}$, $l \in \{1,\ldots,R\}$, where $y_{l}^{\tau_{i}}$ is defined as in \ref{local-trace-notationw1}.
\end{nprop}

\begin{nproof}
Choose $x_{l} \in X_{l} \setminus \Omega_{l}$, $l \in \{1,\ldots,R\}$; by Remark~\ref{Fb}(ii) and Proposition~\ref{I}  there are  pairwise orthogonal weighted $(\mathcal{F},\eta)$-excisors $(E_{x_{l}},\pi_{x_{l}},\sigma_{x_{l}},\kappa_{x_{l}})$; note that $E_{x_{l}} \cong M_{r_{l}}$ for all $l$. 

\underline{Claim 1:} For $l \in \{2, \ldots,R\}$, there are $L^{(l)} \in \mathbb{N}$ and pairwise disjoint nonempty subsets 
\[
\Omega_{l,1}, \ldots,\Omega_{l,L^{(l)}} \subset \Omega_{l}
\]
and 
\[
\nu_{m,k}^{(l)} \in \mathbb{Q}_{+}, \; m \in \{1,\ldots,R\}, \; k \in \{1,\ldots,L^{(l)}\},
\]
such that the following hold:
\begin{enumerate}
\item[a)] $\sum_{m=1}^{R} \nu_{m,k}^{(l)} = 1$ for $k \in \{1,\ldots,L^{(l)}\}$ and $\nu^{(l)}_{m,k}=0$ if $m \ge l$,
\item[b)] $\bigcup_{k=1}^{L^{(l)}} \Omega_{l,k} = \Omega_{l}$,
\item[c)] for each $x \in \Omega_{l,k}$, $k \in \{1,\ldots,L^{(l)}\}$, there are finite subsets 
\[
Y_{l,x,m} \subset X_{m} \setminus \Omega_{m}, \; m \in \{1,\ldots,l-1\},
\]
and, for each $y \in Y_{l,x,m}$, there is a positive integer 
\[
\mu_{l,x,y} \in \mathbb{N}
\]
such that
\begin{equation}
\label{J2}
\textstyle
\pi_{x} \approx_{\mathrm{u}} \bigoplus_{m=1}^{l-1} \left( \bigoplus_{y \in Y_{l,x,m}} \left(  \bigoplus_{1}^{\mu_{l,x,y}} \pi_{y}   \right)  \right)
\end{equation}
and
\begin{equation}
\label{J4}
\textstyle
\sum_{y \in Y_{l,x,m}} \mu_{l,x,y} \cdot r_{m} = \nu_{m,k}^{(l)} \cdot r_{l}, \; m \in \{1,\ldots,l-1\};
\end{equation}
moreover, we have
\begin{equation}
\label{J1a}
\textstyle
(E_{x},\pi_{x},\sigma_{x},\kappa_{x}) \sim_{(\mathcal{F},\eta)}  
\bigoplus_{\substack{m \in \{1,\ldots,l-1\} \\ y \in Y_{l,x,m}}} \frac{r_{m} \mu_{l,x,y}}{r_{l}} \cdot (E_{y},\pi_{y},\sigma_{y},\kappa_{y}).
\end{equation}
\end{enumerate}

\underline{Proof of Claim 1:} 
Note that we do not rule out $\Omega_{l} = \emptyset$. In this case, we set $L^{(l)} = 0$ and there is nothing to show.

Now for each $l \in \{2,\ldots,R\}$ and $x \in \Omega_{l}$, $\pi_{x}$ is unitarily equivalent to a direct sum of irreducible representations of $B_{l-1}$. More precisely, there are finite subsets $Y_{l,x,m} \subset X_{m} \setminus \Omega_{m}$, $m \in \{1,\dots,l-1\}$, and for each $y \in Y_{l,x,m}$ there is $\mu_{l,x,y} \in \mathbb{N}$ such that 
\[
\textstyle
\pi_{x} \approx_{\mathrm{u}} \left(  \bigoplus_{m=1}^{l-1} \left( \bigoplus_{y \in Y_{l,x,m}} \left( \bigoplus_{1}^{\mu_{l,x,y}} \pi_{y}   \right)   \right)   \right).
\]
The ranks of the representations of $B_{l-1}$ (with multiplicities) add up to the rank of $\pi_{x}$, so that
\begin{equation}
\label{J1}
\textstyle
\sum_{m=1}^{l-1} \left( \sum_{y \in Y_{l,x,m}}  \mu_{l,x,y} \cdot r_{m} \right) = r_{l}.
\end{equation}
From this it follows that there are only finitely many, say $L^{(l)}$,  values for tuples of the form 
\[
(\mu_{l,x,y})_{m \in \{1,\ldots,l-1\} , y \in Y_{l,x,m}}
\]
where $x$ ranges over $\Omega_{l}$. Decompose $\Omega_{l}$ into $L^{(l)}$ pairwise disjoint nonempty subsets $\Omega_{l,k}$, $k\in \{1,\ldots, L^{(l)}\}$, such that the maps
\[
x \mapsto (\mu_{l,x,y})_{m \in \{1,\ldots,l-1\} , y \in Y_{l,x,m}}
\]
are constant on each $\Omega_{l,k}$. For $k \in \{1,\ldots,L^{(l)}\}$ and $m \in \{1,\ldots,l-1\}$ set 
\begin{equation}
\label{Jw1}
\textstyle
\nu_{m,k}^{(l)} := \sum_{y \in Y_{l,x,m}}   \frac{r_{m}\mu_{l,x,y}}{r_{l}};
\end{equation}
set
\[
\nu_{m,k}^{(l)} := 0 \mbox{ for } m \ge l.
\]
Then property a) of Claim 1 holds by \eqref{J1};  b) and c) hold by construction.

\underline{Claim 2:} 
For $l \in \{1,\ldots,R\}$ and $k \in \{1,\ldots,L^{(l)}\}$ let
\[
\kappa^{(l)}_{k}:E \to \mathcal{Q}
\]
be a unital $^{*}$-homomorphism such that 
\begin{equation}
\label{Jw2}
\tau_{\mathcal{Q}} \circ \kappa^{(l)}_{k}(1_{E_{x_{m}}}) = \nu_{m,k}^{(l)} , \; m \in \{1,\ldots,R\}
\end{equation}
(such $\kappa^{(l)}_{k}$ exist by Claim 1a)).

Then
\begin{equation}
\label{J3c}
(E,\pi,\sigma,\kappa_{x_{l}}) \sim_{(\mathcal{F},\eta)} (E,\pi,\sigma,\kappa_{k}^{(l)}),
\end{equation}

where we have slightly misused notation by writing $\kappa_{x_{l}}$ for the (canonical) extension of $\kappa_{x_{l}}: E_{x_{l}} \to \mathcal{Q}$ to all of $E$. Moreover (cf.\ \ref{Ia} for notation),
\begin{equation}
\label{J3d}
\textstyle
(E,\pi,\sigma,\kappa_{k}^{(l)}) \sim_{(\mathcal{F},\eta)} \bigoplus_{m=1}^{R} \nu_{m,k}^{(l)} \cdot (E,\pi,\sigma,\kappa_{x_{m}}).
\end{equation}

\underline{Proof of Claim 2:} 
Take $\Omega_{l,k}$ and $\nu^{(l)}_{m,k}$, $m \in \{1,\ldots,l-1\}$ as in Claim 1; fix $x \in \Omega_{l,k}$ and let $Y_{l,x,m}$, $\mu_{l,x,y}$ be as in Claim 1c).

Note that since our recursive subhomogeneous decomposition is $(\mathcal{F},\eta)$-connected, we have
\begin{equation}
\label{J3a}
(E_{x},\pi_{x},\sigma_{x},\kappa_{x}) \sim_{(\mathcal{F},\eta)} (E_{x_{l}},\pi_{x_{l}},\sigma_{x_{l}},\kappa_{x_{l}})
\end{equation}
and, for each $m \in \{1,\ldots,l-1\}$ and $y \in Y_{l,x,m}$,
\begin{equation}
\label{J3}
(E_{y},\pi_{y},\sigma_{y},\kappa_{y}) \sim_{(\mathcal{F},\eta)} (E_{x_{m}},\pi_{x_{m}},\sigma_{x_{m}},\kappa_{x_{m}}),
\end{equation}
with notation as in \ref{Fc}.

Moreover note that
\begin{equation}
\label{J3b}
 (E_{x_{l}},\pi_{x_{l}},\sigma_{x_{l}},\kappa_{x_{l}})  \sim_{(\mathcal{F},\eta)} (E,\pi,\sigma,\kappa_{x_{l}})
\end{equation}
by Remark~\ref{M}(iii). It follows from \eqref{J2} that
\[
\textstyle
\kappa_{x} \circ \pi_{x} \approx_{\mathrm{u}} \bigoplus_{\substack{ m \in \{1,\ldots,l-1\}\\ y \in Y_{l,x,m}}} \left(\frac{r_{m}}{r_{l}}\mu_{l,x,y} \right) \cdot \kappa_{y} \circ \pi_{y},
\]
cf.\ \ref{Ia} for notation.

But then by  Proposition~\ref{H}   we have
\begin{eqnarray*}
\textstyle
\lefteqn{(E_{x},\pi_{x},\sigma_{x},\kappa_{x})  }\\
\textstyle
& \stackrel{\eqref{J1a}}{\sim_{(\mathcal{F},\eta)}}   &\textstyle  \bigoplus_{\substack{ m \in \{1,\ldots,l-1\}\\ y \in Y_{l,x,m}}} \left(\frac{r_{m}}{r_{l}}\mu_{l,x,y} \right) \cdot (E_{y}, \pi_{y},\sigma_{y},\kappa_{y}) \\
&  \stackrel{\eqref{J3}}{\sim_{(\mathcal{F},\eta)}} & \textstyle  \bigoplus_{\substack{ m \in \{1,\ldots,l-1\}\\ y \in Y_{l,x,m}}} \left(\frac{r_{m}}{r_{l}}\mu_{l,x,y} \right) \cdot (E_{x_{m}}, \pi_{x_{m}},\sigma_{x_{m}},\kappa_{x_{m}}) \\
&   \stackrel{\eqref{Jw1}}{\sim_{(\mathcal{F},\eta)}} & \textstyle  \bigoplus_{ m \in \{1,\ldots,l-1\}} \nu_{m,k}^{(l)} \cdot (E_{x_{m}}, \pi_{x_{m}},\sigma_{x_{m}},\kappa_{x_{m}}) \\
&  \stackrel{\eqref{J3b}}{\sim_{(\mathcal{F},\eta)}} &  \textstyle  \bigoplus_{ m \in \{1,\ldots,l-1\}} \nu_{m,k}^{(l)} \cdot (E, \pi,\sigma,\kappa_{x_{m}}) \\
&   \stackrel{\eqref{Jw2}}{\sim_{(\mathcal{F},\eta)}} & (E,\pi,\sigma,\kappa_{k}^{(l)}).
\end{eqnarray*}
Combining this with \eqref{J3a} and \eqref{J3b} now yields  \eqref{J3c}; it also shows \eqref{J3d}.  We have now  verified Claim 2.

\underline{Claim 3:}
Let $p \in B$ be  a projection such that 
\begin{equation}
\label{J6}
\frac{1}{r_{l}} \cdot \mathrm{rank} (p|_{X_{l}}) \equiv \xi_{l} \in \mathbb{Q}
\end{equation} 
is constant for each $l \in \{1,\ldots,R\}$. 

Then, for each $l \in \{2,\ldots,R\}$, the $\xi_{l}$ satisfy the relations
\begin{equation}
\label{J5}
\textstyle
\xi_{l} = \sum_{m=1}^{R} \nu_{m,k}^{(l)}\cdot \xi_{m}  , \; k \in \{1,\ldots,L^{(l)}\},
\end{equation} 
where the $\nu_{m,k}^{(l)}$ are as in Claim 1.

\underline{Proof of Claim 3:}
For $l \in \{2,\ldots,R\}$ and $k \in \{1,\ldots,L^{(l)}\}$ choose $x \in \Omega_{l,k} \subset X_{l}$ and let $Y_{l,x,m}$ and $\mu_{l,x,y}$  be as in Claim 1c). 

We have
\begin{eqnarray*}
\xi_{l} & \stackrel{\eqref{J6}}{=} & \textstyle \frac{1}{r_{l}} \cdot \mathrm{rank}(\pi_{x}(p )) \\
&  \stackrel{\eqref{J2}}{=}  &\textstyle  \frac{1}{r_{l}} \cdot \left( \sum_{m=1}^{l-1}  \left(  \sum_{y \in Y_{l,x,m}}  \left(  \sum_{1}^{\mu_{l,x,y}}  \mathrm{rank}(\pi_{y}(p ))  \right) \right)  \right)\\
&  \stackrel{\eqref{J6}}{=}  & \textstyle  \frac{1}{r_{l}} \cdot \left( \sum_{m=1}^{l-1}  \left(  \sum_{y \in Y_{l,x,m}} \mu_{l,x,y}  \cdot r_{m} \cdot \xi_{m} \right)  \right)\\
&  \stackrel{\eqref{J4}}{=}  &   \textstyle    \sum_{m=1}^{l-1} \nu_{m,k}^{(l)} \cdot \xi_{m} ,
\end{eqnarray*}
so \eqref{J5} holds and Claim 3 is proven.

Before moving on to Claim 4, let us set
\[
\textstyle
\bar{L}:= \sum_{l=2}^{R} L^{(l)}
\]
and define $\bar{L} \times R$ matrices
\[
T_{+}:=\begin{pmatrix}
\nu_{1,1}^{(2)} & 0 & \ldots \\
\vdots & \vdots \\
\nu_{1,L^{(2 )}}^{(2)} & 0 & \ldots\\
\vdots \\
\nu_{1,1}^{(l)} & \ldots & \nu_{l-1,1}^{(l)} & 0 & \ldots\\
\vdots & & \vdots & \vdots \\
\nu_{1,L^{(l)}}^{(l)} & \ldots & \nu_{l-1,L^{(l)}}^{(l)} & 0 & \ldots\\
\vdots \\
\nu_{1,1}^{(R )} && \ldots && \nu_{R-1,1}^{(R )} & 0 \\
\vdots & &&& \vdots & \vdots \\
\nu_{1,L^{(R )}}^{(R )} && \ldots && \nu_{R-1,L^{(R )}}^{(R )} & 0 \\
\end{pmatrix}
\]
and 
\[
T_{-}:=\begin{pmatrix}
0 & 1 & \ldots \\
\vdots & \vdots \\
0& 1 & \ldots\\
\vdots \\
0 && \ldots &&0& 1 & \ldots\\
\vdots &&& & \vdots & \vdots \\
0 & &\ldots && 0 & 1 & \ldots\\
\vdots \\
0 &&& \ldots &&&0 & 1 \\
\vdots & &&& \vdots & \vdots \\
0 &&& \ldots &&& 0 & 1 \\
\end{pmatrix}
\]
and note that, with these definitions, $\xi_{1},\ldots,\xi_{R}$ satisfy the equation \eqref{J5} for $l \in \{2,\ldots,R\}$, $k \in \{1,\ldots,L^{(l)}\}$ if and only if 
\begin{equation}
\label{J6a}
\begin{pmatrix}
\xi_{1}\\
\vdots\\
\xi_{R}
\end{pmatrix}
\in \mathrm{ker}(T_{+} - T_{-}).
\end{equation}

\underline{Claim 4:}
Suppose we have
\[
\underline{\xi}= \begin{pmatrix}
\xi_{1}\\
\vdots\\
\xi_{R}
\end{pmatrix}
\in \mathrm{ker}(T_{+} - T_{-}) \cap \mathbb{N}^{R}.
\]
Then there are $\bar{N} \in \mathbb{N}$ and a projection
\[
\textstyle
p \in B \otimes M_{\bar{N}} \subset \bigoplus_{l=1}^{R} \mathcal{C}(X_{l}) \otimes M_{r_{l}} \otimes M_{\bar{N}}
\]
such that
\begin{equation}
\label{claim4-1}
\frac{1}{r_{l}} \cdot \mathrm{rank}(p|_{X_{l}}) \equiv \xi_{l}
\end{equation}
for $l \in \{1,\ldots,R\}$.

\underline{Proof of Claim 4:}
Take a trivial projection $p_{1}$ in $B_{1} = \mathcal{C}(X_{1}) \otimes M_{r_{1}} \otimes M_{\bar{N}}$ (for $\bar{N}$ large enough) with rank $r_{1} \xi_{1}$.

Now suppose we have constructed projections $p_{1},\ldots,p_{l}$ in $B_{1},\ldots,B_{l}$, respectively, such that 
\begin{equation}
\label{Jw3}
\frac{1}{r_{m}} \cdot \mathrm{rank}(p_{l'}|_{X_{m}}) \equiv \xi_{m} \mbox{ for } 1 \le m \le l' \le l
\end{equation}
and
\[
\psi'_{l'}(p_{l'+1}) = p_{l'} \mbox{ for } l' \in \{1,\ldots,l-1\},
\]
where 
\[
\psi'_{l'}:B_{l'+1} \twoheadrightarrow B_{l'}
\]
denotes the canonical surjection, cf.\ \eqref{dK1}. 

If $\Omega_{l+1} = \emptyset$, then 
\[
B_{l+1} \cong \mathcal{C}(X_{l+1}) \otimes M_{r_{l+1}} \oplus B_{l}
\]
and we may define
\[
p_{l+1}:= q_{l+1} \oplus p_{l},
\]
where
\[
q_{l+1} \in \mathcal{C}(X_{l+1}) \otimes M_{r_{l+1}} \otimes M_{\bar{N}}
\]
is a trivial projection with rank $\xi_{l+1} r_{l+1}$.  

If $\Omega_{l+1} \neq \emptyset$, then $\phi_{l} \otimes \mathrm{id}_{M_{\bar{N}}} (p_{l})$ is a projection in $\mathcal{C}(\Omega_{l+1}) \otimes M_{r_{l+1}} \otimes M_{\bar{N}}$ and, for $x \in \Omega_{l+1}$, we have
\begin{eqnarray*}
\mathrm{rank}((\phi_{l} \otimes \mathrm{id}_{M_{\bar{N}}})(p_{l})(x)) & \stackrel{\eqref{J2}}{=} & \textstyle \sum_{m=1}^{l} \sum_{y \in Y_{l+1,x,m}} \mu_{l+1,x,y} \cdot \mathrm{rank}(p_{l}|_{X_{m}}) \\
& \stackrel{\eqref{Jw3}}{=} & \textstyle \sum_{m=1}^{l} \sum_{y \in Y_{l+1,x,m}} \mu_{l+1,x,y} \cdot \xi_{m} r_{m} \\
& \stackrel{\eqref{J4}}{=} & \textstyle \sum_{m=1}^{l} \nu_{m,k}^{(l+1)} \cdot r_{l+1} \xi_{m} \\
& \stackrel{\eqref{J5}}{=} & \xi_{l+1} r_{l+1}.
\end{eqnarray*}
But then by hypothesis, $(\phi_{l} \otimes \mathrm{id}_{M_{\bar{N}}})(p_{l})$ lifts to a projection $p'_{l+1}$ in $\mathcal{C}(X_{l+1}) \otimes M_{r_{l+1}} \otimes M_{\bar{N}}$; by changing $p'_{l+1}$ on those components of $X_{l+1}$ which do not intersect $\Omega_{l+1}$, if necessary, we may assume that $p'_{l+1}$ has constant rank $\xi_{l+1} r_{l+1}$ on $X_{l+1}$. Now
\[
p_{l+1} := p'_{l+1} \oplus p_{l} \in B_{l+1} \otimes M_{\bar{N}} \subset \mathcal{C}(X_{l+1}) \otimes M_{r_{l+1}} \otimes M_{\bar{N}}
\]
satisfies 
\[
\frac{1}{r_{m}} \cdot \mathrm{rank}(p_{l+1}|_{X_{m}}) \equiv \xi_{m}
\]
for $1 \le m \le l+1$. Proceed inductively  to construct $p_{1},p_{2},\ldots,p_{R}$, then 
\[
p:= p_{R}
\] 
will be as desired.  This proves Claim 4.

Let $\underline{\xi}$, $\bar{N}$ and $p$ be as in Claim 4. Since $(\tau_{0})_{*} = (\tau_{1})_{*}$, we have 
\[
(\tau_{0} \otimes \mathrm{tr}_{M_{\bar{N}}}) (p ) =  (\tau_{1} \otimes \mathrm{tr}_{M_{\bar{N}}})(p ),
\]
whence
\[
\textstyle
\sum_{l=1}^{R} \xi_{l} \cdot y_{l}^{\tau_{0}} = \sum_{l=1}^{R} \xi_{l} \cdot y_{l}^{\tau_{1}},
\]
cf.\ \ref{local-trace-notation}. But this just means that 
\[
\langle \underline{\xi}, \underline{y}^{(0)} \rangle = \langle \underline{\xi}, \underline{y}^{(1)} \rangle
\]
or, equivalently, 
\begin{equation}
\label{J7}
\underline{\xi} \perp (\underline{y}^{(0)} - \underline{y}^{(1)}) \mbox{ in } \mathbb{R}^{R},
\end{equation}
where
\begin{equation}
\label{Jw16}
\underline{y}^{(i)} = \begin{pmatrix} y^{\tau_{i}}_{1} \\ \vdots \\ y^{\tau_{i}}_{R} \end{pmatrix}, \; i \in \{0,1\}.
\end{equation}
Since $1_{B} \in B$ is a projection with $\frac{1}{r_{l}} \cdot \mathrm{rank} (1_{B}|_{X_{l}})= 1$ for all $l$, we see from Claim 3 and \eqref{J6a} that
\begin{equation}
\label{Jw5}
\underline{r} := \begin{pmatrix} 1\\ \vdots \\ 1   \end{pmatrix} \in \mathrm{ker}(T_{+} - T_{-}) \cap \mathbb{Z}^{R}.
\end{equation}
But positive integer multiples of $\underline{r}$ are also in $ \mathrm{ker}(T_{+} - T_{-}) \cap \mathbb{N}^{R}$, from which follows that
\begin{equation}
\label{Jw4}
 \mathrm{ker}(T_{+} - T_{-}) \cap \mathbb{Z}^{R} =  \mathrm{ker}(T_{+} - T_{-}) \cap \mathbb{N}^{R} -  \mathrm{ker}(T_{+} - T_{-}) \cap \mathbb{N}^{R}.
\end{equation}
Moreover, Claim 4 and \eqref{Jw4} imply 
\[
 \mathrm{ker}(T_{+} - T_{-}) \cap \mathbb{Z}^{R} \perp (\underline{y}^{(0)} -  \underline{y}^{(1)} ) \mbox{ in } \mathbb{R}^{R},
\]
whence 
\[
 \mathrm{ker}(T_{+} - T_{-}) \cap \mathbb{Q}^{R} \perp (\underline{y}^{(0)} -  \underline{y}^{(1)} ) \mbox{ in } \mathbb{R}^{R};
\]
since $T_{+}$ and $T_{-}$ have only rational coefficients, it follows that $\mathrm{ker}(T_{+}-T_{-}) \cap \mathbb{Q}^{R}$ is dense in $\mathrm{ker}(T_{+}-T_{-})$, whence 
\[
 \mathrm{ker}(T_{+} - T_{-})  \perp (\underline{y}^{(0)} -  \underline{y}^{(1)} ) \mbox{ in } \mathbb{R}^{R}.
\]

By elementary linear algebra we have 
\[
(\mathrm{ker}(T_{+} - T_{-}))^{\perp} = \mathrm{Im}(T_{+} - T_{-})^{*},
\]
so there is $\underline{\zeta} \in \mathbb{R}^{\bar{L}}$ such that 
\[
(T_{+}-T_{-})^{*} \underline{\zeta} = \underline{y}^{(0)} - \underline{y}^{(1)}.
\]
We may then write 
\[
\underline{\zeta} = \underline{\zeta}_{+} - \underline{\zeta}_{-} \mbox{ with } \underline{\zeta}_{+},\underline{\zeta}_{-} \in \mathbb{R}_{+}^{\bar{L}}
\]
to obtain the equation 
\[
\underline{y}^{(0)} + T_{+}^{*} \underline{\zeta}_{-} + T_{-}^{*} \underline{\zeta}_{+} = \underline{y}^{(1)} + T_{+}^{*} \underline{\zeta}_{+} + T_{-}^{*} \underline{\zeta}_{-} 
\]
in  $ \mathbb{R}^{R}$, in which all vectors and matrices have only positive entries.

We wish to interpret the entries of the $\underline{y}^{(i)}$, $\underline{\zeta}_{+}$ and $\underline{\zeta}_{-}$ as coefficients of sums of $(\mathcal{F},\eta)$-bridges. To this end, we have to approximate them by rationals. Let us first set  
\begin{equation}
\label{Jw19}
\alpha:= \frac{\bar{\beta}}{8R} (\le 1).
\end{equation}

\underline{Claim 5:} 
There are 
\[
g^{(i)}= (g^{(i)}_{l})_{l \in \{1,\dots,R\}} \in \mathbb{Q}^{R}, \; i \in \{0,1\},
\]
\[
z_{+} = (z_{+,k}^{(l)})_{\substack{l \in \{2,\ldots,R\}\\k \in \{1,\ldots,L^{(l)}\}}} \in \mathbb{Q}_{+}^{\bar{L}}
\]
and 
\[
z_{-} = (z_{-,k}^{(l)})_{\substack{l \in \{2,\ldots,R\}\\k \in \{1,\ldots,L^{(l)}\}}} \in \mathbb{Q}_{+}^{\bar{L}}
\]
satisfying
\begin{equation}
\label{Jw15}
\| g^{(i)} - \underline{y}^{(i)} \|_{\mathrm{max}} \le \alpha , \; i \in \{0,1\},
\end{equation}
\[
\| z_{+} - \underline{\zeta}_{+} \|_{\mathrm{max}}, \; \| z_{-} - \underline{\zeta}_{-} \|_{\mathrm{max}} \le \alpha , 
\]
\begin{equation}
\label{J8}
\| g^{(0)} + T_{+}^{*} z_{-} + T_{-}^{*} z_{+} - (g^{(1)} + T_{+}^{*} z_{+} + T_{-}^{*} z_{-})  \|_{\mathrm{max}} \le \alpha, 
\end{equation}
\begin{equation}
\label{J9}
\langle \underline{r},g^{(0)} \rangle = \langle \underline{r},g^{(1)} \rangle
\end{equation}
(with $\underline{r}$ as in \eqref{Jw5}), and
\begin{equation}
\label{J10}
\langle \underline{r}, T_{+}^{*} z_{-} + T_{-}^{*}z_{+} \rangle = \langle \underline{r}, T_{+}^{*} z_{+} + T_{-}^{*}z_{-} \rangle.
\end{equation}

\underline{Proof of Claim 5:} Easy.

We now set 
\begin{equation}
\label{Jw8}
\underline{v}^{(0)} := T_{+}^{*} z_{-} + T^{*}_{-}z_{+}, \; \underline{v}^{(1)}:= T^{*}_{+}z_{+} + T^{*}_{-} z_{-},
\end{equation}
\begin{equation}
\label{Jw9}
\textstyle
G:= \langle \underline{r},g^{(0)} \rangle = \sum_{l=1}^{R} g_{l}^{(0)}  \stackrel{\eqref{J9}}{=} \langle \underline{r},g^{(1)} \rangle,
\end{equation}
\begin{equation}
\label{Jw6}
\textstyle
Z_{+}:= \langle \underline{r}, T^{*}_{-}z_{+} \rangle = \sum_{\substack{m \in \{2,\ldots,R\}\\ k \in \{1,\ldots,L^{(m)}\}}} z_{+,k}^{(m)}
\end{equation}
and
\begin{equation}
\label{Jw7}
\textstyle
Z_{-}:= \langle \underline{r}, T^{*}_{-}z_{-} \rangle = \sum_{\substack{m \in \{2,\ldots,R\}\\ k \in \{1,\ldots,L^{(m)}\}}} z_{-,k}^{(m)}.
\end{equation}
Note that 
\begin{equation}
\label{Jw17}
|G-1| \stackrel{\eqref{Jw16},\eqref{local-trace-notationw2}}{=} |\langle \underline{r}, g^{(0)} - \underline{y}^{(0)} \rangle | \stackrel{\eqref{Jw15}}{\le} R \alpha.
\end{equation}

For any 
\[
z = (z_{k}^{(m)})_{\substack{m \in \{2,\ldots,R\}, \\ k \in \{1,\ldots,L^{(m)}\}}} \in \mathbb{R}^{\bar{L}}
\]
we compute 
(observing that $\nu^{(m)}_{l,k}=0$ if $m \le l$)
\begin{eqnarray}
\langle \underline{r}, T^{*}_{+}z\rangle & = & \textstyle \sum_{l=1}^{R} \sum_{m=2}^{R} \sum_{k=1}^{L^{(m)}} \nu^{(m)}_{l,k} \cdot z^{(m)}_{k} \nonumber \\
& = &  \textstyle \sum_{m=2}^{R} \sum_{k=1}^{L^{(m)}} \left( \sum_{l=1}^{R} \nu^{(m)}_{l,k} \right) \cdot z^{(m)}_{k} \nonumber \\
& \stackrel{\mathrm{Claim~1a)}}{=} &  \textstyle \sum_{m=2}^{R} \sum_{k=1}^{L^{(m)}}  z^{(m)}_{k} \nonumber \\
& = & \langle \underline{r}, T^{*}_{-}z \rangle, \label{Jw20}
\end{eqnarray}
so that in particular
\[
\langle \underline{r}, T^{*}_{+}z_{+} \rangle = Z_{+} = \langle \underline{r}, T^{*}_{-}z_{+} \rangle
\]
and
\[
\langle \underline{r}, T^{*}_{+}z_{-} \rangle = Z_{-} = \langle \underline{r}, T^{*}_{-}z_{-} \rangle,
\]
see \eqref{Jw6}, \eqref{Jw7}.

We set
\begin{equation}
\label{Jw11}
V:= \langle \underline{r}, \underline{v}^{(0)} \rangle \stackrel{\eqref{Jw8}}{=}  \langle \underline{r}, \underline{v}^{(1)} \rangle  \stackrel{\eqref{Jw6},\eqref{Jw7}}{=} Z_{+} + Z_{-}.
\end{equation}
By \eqref{J8}, we may choose $w_{+}, w_{-} \in \mathbb{Q}_{+}^{R}$ such that 
\begin{equation}
\label{Jw21}
-(g^{(0)} + v^{(0)}) + (g^{(1)} +v^{(1)}) = w_{+} - w_{-}
\end{equation}
and such that
\begin{equation}
\label{Jw18}
\|w_{+}\|_{\mathrm{max}}, \; \|w_{-}\|_{\mathrm{max}} \le \alpha.
\end{equation}
Set
\begin{equation}
\label{Jw10}
W:= \langle \underline{r}, w_{+}\rangle \stackrel{\eqref{J9},\eqref{J10}}{=} \langle \underline{r}, w_{-}\rangle.
\end{equation}
Note also that
\begin{eqnarray*}
|G+W-1| & \le & \frac{1}{2} |\langle \underline{r}, g^{(0)} \rangle + \langle  \underline{r}, g^{(1)} \rangle \\
& & + \langle \underline{r}, w_{+} \rangle + \langle  \underline{r}, w_{-} \rangle \\
& & - \langle \underline{r}, \underline{y}^{(0)} \rangle - \langle  \underline{r}, \underline{y}^{(1)} \rangle | \\
& \le & \frac{1}{2}( |\langle \underline{r}, g^{(0)} -\underline{y}^{(0)} \rangle | + |\langle  \underline{r}, g^{(1)}  - \underline{y}^{(1)}\rangle | \\
& & + |\langle \underline{r}, w_{+} \rangle | + |\langle  \underline{r}, w_{-} \rangle | \\
& \le & 2R \alpha \, ( = \bar{\beta}/4 < 1/2),
\end{eqnarray*}
whence
\begin{equation}
\label{Jw14}
\textstyle \frac{1}{G+W} \le 1+4R \alpha.
\end{equation}
(The estimate is only nontrivial for $0 \le 1-G+W \le 2R\alpha (< 1/2)$; in this case we use that $1/(1-\theta) \le 1+2 \theta$ for $0 \le \theta \le 1/2$.)

Let $e_l$ denote the unit of the $l^{\mathrm{th}}$ copy of $\Q$ in $\oplus_{l=1}^R \Q$. We now choose unital $^{*}$-homomorphisms
\[
\textstyle
\gamma_{0},\, \gamma_{1},\, \tilde{\gamma}: \bigoplus_{l=1}^{R} \mathcal{Q} \to \mathcal{Q}
\]
and 
\[
\bar{\gamma}: \mathcal{Q} \oplus \mathcal{Q} \to \mathcal{Q}
\]
such that 
\begin{equation}
\label{Jw13}
\textstyle \tau_{\mathcal{Q}} \circ \gamma_{0}(e_{l}) = \frac{g^{(0)}_{l} + w_{+,l}}{G+W},
\end{equation}
\[ \textstyle
\tau_{\mathcal{Q}} \circ \gamma_{1}(e_{l}) = \frac{g^{(1)}_{l} + w_{-,l}}{G+W},
\]
\[\textstyle
\tau_{\mathcal{Q}} \circ \tilde{\gamma}(e_{l}) = \frac{g^{(0)}_{l} + v_{l}^{(0)}}{G+V},
\]
for $l \in \{1,\ldots,R\}$ and 
\begin{equation}
\label{Jw22}
\textstyle \tau_{\mathcal{Q}} \circ \bar{\gamma}((1,0)) = \frac{G+W}{2G+V+W},
\end{equation}
\[
\textstyle \tau_{\mathcal{Q}} \circ \bar{\gamma}((0,1)) = \frac{G+V}{2G+V+W};
\]
these exist by \eqref{Jw9}, \eqref{Jw10} and \eqref{Jw11}.

Next observe that
\[
\textstyle \bar{y}_{0,l} \stackrel{\eqref{Jw12}}{:=}  \tau_{\mathcal{Q}}(\gamma_{0}(e_{l})) \stackrel{\eqref{Jw13}}{=} \frac{g^{(0)}_{l} + w_{+,l}}{G+W},
\]
so 
\begin{eqnarray*}
|\bar{y}_{0,l} - y^{\tau_{0}}_{l}| & = &\textstyle  \left| \frac{1}{G+W} (g^{(0)}_{l} + w_{+,l}) - y_{l}^{\tau_{0}} \right| \\
& = & \textstyle \frac{1}{G+W} |g^{(0)}_{l} + w_{+,l} - (G+W) y_{l}^{\tau_{0}} | \\
& \le &\textstyle  \frac{1}{G+W} (|g^{(0)}_{l} - y^{\tau_{0}}_{l}| + |G-1| y^{\tau_{0}}_{l} + w_{+,l} + W y^{\tau_{0}}_{l}) \\
& \stackrel{\eqref{local-trace-notationw1}}{\le} &\textstyle  \frac{1}{G+W} (|g^{(0)}_{l} - y^{\tau_{0}}_{l}| + |G-1| + 2W ) \\
& \le & (1+4R\alpha) (\alpha + 3R\alpha) \\
& \stackrel{\eqref{Jw19}}{<} & 8 R \alpha \\
&  \stackrel{\eqref{Jw19}}{=} & \bar{\beta}
\end{eqnarray*}
and \eqref{Jw24} holds. Here, for the third inequality we have used \eqref{Jw14}, \eqref{Jw15}, \eqref{Jw16}, \eqref{Jw17} and \eqref{Jw18}.

Set
\[
z:= (z^{(l)}_{k})_{\substack{l \in \{2,\ldots,R\} \\ k \in \{1,\ldots, L^{(l)} \}}} \in \mathbb{Q}_{+}^{\bar{L}}
\]
and
\[
Z:= \langle \underline{r},T^{*}_{+}z\rangle;
\]
note that
\begin{equation*}
Z= \langle \underline{r},T^{*}_{-}z\rangle;
\end{equation*}
by \eqref{Jw20}. We then compute
\begin{eqnarray}
\lefteqn{\textstyle \bigoplus_{m=1}^{R}  \frac{1}{Z} (T_{+}^{*}z)_{m} \cdot (E,\pi,\sigma,\kappa_{x_{m}}) } \nonumber \\
& = &   \textstyle \bigoplus_{m=1}^{R} \left( \sum_{l=2}^{R} \sum_{k=1}^{L^{(l)}}  \nu_{m,k}^{(l)} \frac{z_{k}^{(l)}}{Z}     \right) \cdot (E,\pi,\sigma,\kappa_{x_{m}})    \nonumber  \\
& \stackrel{\ref{Ia}\mathrm{(iii)},\ref{H}}{\sim_{(\mathcal{F},\eta)}} &   \textstyle \bigoplus_{l=2}^{R}  \bigoplus_{k=1}^{L^{(l)}}  \frac{z_{k}^{(l)}}{Z}  \cdot \left( \bigoplus_{m=1}^{R}  \nu_{m,k}^{(l)}     \cdot (E,\pi,\sigma,\kappa_{x_{m}})   \right) \nonumber  \\
& \stackrel{\eqref{J3d}}{\sim_{(\mathcal{F},\eta)}} &   \textstyle \bigoplus_{l=2}^{R}  \bigoplus_{k=1}^{L^{(l)}}  \frac{z_{k}^{(l)}}{Z}  \cdot  (E,\pi,\sigma,\kappa_{k}^{(l)})    \nonumber \\
& \stackrel{\eqref{J3c}}{\sim_{(\mathcal{F},\eta)}} &  \textstyle \bigoplus_{l=2}^{R}  \bigoplus_{k=1}^{L^{(l)}}  \frac{z_{k}^{(l)}}{Z}  \cdot  (E,\pi,\sigma,\kappa_{x_{l}})    \nonumber \\
& = &  \textstyle \bigoplus_{m=1}^{R}    \frac{1}{Z}  (T_{-}^{*}z)_{m} \cdot  (E,\pi,\sigma,\kappa_{x_{m}})  .
\label{J11}
\end{eqnarray}

As a consequence, we obtain
\begin{eqnarray}
\lefteqn{\textstyle
\bigoplus_{m=1}^{R} \frac{1}{V} \underline{v}^{(0)}_{m} \cdot (E,\pi,\sigma,\kappa_{x_{m}})} \nonumber \\
 & = & \textstyle \bigoplus_{m=1}^{R} \frac{1}{V}  ((T^{*}_{+}z_{-})_{m} + (T^{*}_{-}z_{+})_{m})   \cdot (E,\pi,\sigma,\kappa_{x_{m}}) \nonumber \\
 & \sim_{(\mathcal{F},\eta)} &  \textstyle \frac{Z_{-}}{V} \cdot \left( \bigoplus_{m=1}^{R} \frac{1}{Z_{-}}  (T^{*}_{+}z_{-})_{m}    \cdot (E,\pi,\sigma,\kappa_{x_{m}}) \right) \nonumber \\
&&   \textstyle \oplus \, \frac{Z_{+}}{V} \cdot \left( \bigoplus_{m=1}^{R} \frac{1}{Z_{+}}   (T^{*}_{-}z_{+})_{m}   \cdot (E,\pi,\sigma,\kappa_{x_{m}})\right) \nonumber \\
 & \stackrel{\eqref{J11}}{\sim_{(\mathcal{F},\eta)}} &    \textstyle \frac{Z_{-}}{V} \cdot \left( \bigoplus_{m=1}^{R} \frac{1}{Z_{-}}  (T^{*}_{-}z_{-})_{m}    \cdot (E,\pi,\sigma,\kappa_{x_{m}}) \right) \nonumber \\
&&   \textstyle  \oplus \, \frac{Z_{+}}{V} \cdot \left( \bigoplus_{m=1}^{R} \frac{1}{Z_{+}}   (T^{*}_{+}z_{+})_{m}   \cdot (E,\pi,\sigma,\kappa_{x_{m}})\right)   \nonumber \\
 & \sim_{(\mathcal{F},\eta)} &\textstyle \bigoplus_{m=1}^{R} \frac{1}{V}  ((T^{*}_{-}z_{-})_{m} + (T^{*}_{+}z_{+})_{m})   \cdot (E,\pi,\sigma,\kappa_{x_{m}}) \nonumber \\
 & = & \textstyle
\bigoplus_{m=1}^{R} \frac{1}{V} \underline{v}^{(1)}_{m} \cdot (E,\pi,\sigma,\kappa_{x_{m}}) \label{J12}.
\end{eqnarray}

We finally compute
\begin{eqnarray*}
\lefteqn{(E,\pi,\sigma,\bar{\gamma}\circ (\bar{\kappa}_{0} \oplus \bar{\kappa}))} \\
& \stackrel{\eqref{Jw22},\eqref{Jw13},\eqref{Jw22a},\ref{M}\mathrm{(iii)}}{\sim_{(\mathcal{F},\eta)}} & \textstyle \bigoplus_{m=1}^{R} \frac{1}{2G+V+W} (g^{(0)}_{m} + w_{+,m} +g^{(0)}_{m} + \underline{v}_{m}^{(0)}) \cdot (E,\pi,\sigma,\kappa_{x_{m}})\\
& \stackrel{\eqref{Jw21}}{=} & \textstyle \bigoplus_{m=1}^{R} \frac{1}{2G+V+W} (g^{(1)}_{m} + \underline{v}_{m}^{(1)}+ w_{-,m} +g^{(0)}_{m} ) \cdot (E,\pi,\sigma,\kappa_{x_{m}})\\
& \stackrel{\ref{H}}{\sim_{(\mathcal{F},\eta)}} & \textstyle \frac{2G+W}{2G+V+W} \cdot \left(\bigoplus_{m=1}^{R} \frac{1}{2G+W} (g^{(1)}_{m} + w_{-,m} +g^{(0)}_{m} ) \cdot (E,\pi,\sigma,\kappa_{x_{m}})  \right)\\
& &\oplus \, \textstyle \frac{V}{2G+V+W} \cdot \left( \bigoplus_{m=1}^{R} \frac{1}{V}  \underline{v}_{m}^{(1)} \cdot (E,\pi,\sigma,\kappa_{x_{m}})  \right)\\
& \stackrel{\ref{H},\eqref{J12}}{\sim_{(\mathcal{F},\eta)}} & \textstyle \bigoplus_{m=1}^{R} \frac{1}{2G+V+W} (g^{(1)}_{m} + w_{-,m} +g^{(0)}_{m} + \underline{v}_{m}^{(0)}) \cdot (E,\pi,\sigma,\kappa_{x_{m}})\\
& \stackrel{\eqref{Jw22},\eqref{Jw13},\eqref{Jw22a}}{\sim_{(\mathcal{F},\eta)}} & (E,\pi,\sigma,\bar{\gamma}\circ (\bar{\kappa}_{1} \oplus \bar{\kappa})),
\end{eqnarray*}
thus establishing \eqref{Jw23}.
\end{nproof}
\en

\bn
\begin{nprop}
\label{K}
Let $B$ be a separable unital recursive subhomogeneous \mbox{$\mathrm{C}^{*}$-algebra} with recursive subhomogeneous decomposition $[B_{l},X_{l},\Omega_{l},r_{l},\phi_{l}]_{l=1}^{R}$ 
and let $\mathcal{F} \subset B^{1}_{+}$ finite and $\eta, \delta >0$ be given.

Let $(E,\pi,\sigma,\bar{\kappa}_{0})$, $(E,\pi,\sigma,\bar{\kappa}_{1})$ and $(E,\pi,\sigma,\bar{\kappa})$ be weighted $(\mathcal{F},\eta)$-excisors and let 
\[
\bar{\gamma}: \mathcal{Q} \oplus \mathcal{Q} \to \mathcal{Q}
\]
be a unital embedding such that
\begin{equation}
\label{Kw1}
(E,\pi,\sigma,\bar{\gamma} \circ (\bar{\kappa}_{0} \oplus \bar{\kappa})) \sim_{(\mathcal{F},\eta)} (E,\pi,\sigma,\bar{\gamma} \circ (\bar{\kappa}_{1} \oplus \bar{\kappa})) 
\end{equation}
(compatible with the decomposition).

Then there is a unital embedding 
\[
\gamma: \mathcal{Q} \oplus \mathcal{Q} \to \mathcal{Q}
\]
such that
\[
(E,\pi,\sigma,\gamma \circ (\bar{\kappa}_{0} \oplus \bar{\kappa})) \sim_{(\mathcal{F},\eta)} (E,\pi,\sigma,\gamma \circ (\bar{\kappa}_{1} \oplus \bar{\kappa})) 
\]
(also compatible with the decomposition) and 
\[
|\tau_{\mathcal{Q}}(\gamma((p,0))) - \tau_{\mathcal{Q}}( p)| < \delta
\]
for every projection $p \in \mathcal{Q}$, in particular
\[
\tau_{\mathcal{Q}}(\gamma((1_{\mathcal{Q}},0))) > 1-\delta.
\]
\end{nprop}

\begin{nproof}
Choose $N \in \mathbb{N}$ so large that 
\[
\frac{\tau_{\mathcal{Q}}(\bar{\gamma}((0,1)))}{N \cdot \tau_{\mathcal{Q}}(\bar{\gamma}((1,0))) + \tau_{\mathcal{Q}}(\bar{\gamma}((0,1)))} < \delta
\]
and a unital embedding
\[
\theta: \mathbb{C}^{N+1} \otimes \mathcal{Q} \to \mathcal{Q} 
\]
such that 
\[
\tau_{\mathcal{Q}}(\theta(e_{i} \otimes 1_{\Q})) = \frac{\tau_{\mathcal{Q}}(\bar{\gamma}((1,0)))}{N \cdot \tau_{\mathcal{Q}}(\bar{\gamma}((1,0))) + \tau_{\mathcal{Q}}(\bar{\gamma}((0,1)))}
\]
for $i \in \{1,\ldots,N\}$ and
\[
\tau_{\mathcal{Q}}(\theta(e_{N+1}\otimes 1_{\Q})) = \frac{\tau_{\mathcal{Q}}(\bar{\gamma}((0,1)))}{N \cdot \tau_{\mathcal{Q}}(\bar{\gamma}((1,0))) + \tau_{\mathcal{Q}}(\bar{\gamma}((0,1)))}.
\]

Define
\[
\gamma:\mathcal{Q} \oplus \mathcal{Q}  \to \mathcal{Q}
\]
by
\[
\textstyle
\gamma:= \theta \circ \left(\left(\left( \sum_{i=1}^{N} e_{i}   \right) \otimes \mathrm{id}_{\mathcal{Q}} \right)   \oplus \left( e_{N+1} \otimes \mathrm{id}_{\mathcal{Q}}   \right)\right)
\]
and
\[
\gamma_{j} :\mathcal{Q} \oplus \mathcal{Q} \oplus \mathcal{Q} \oplus \mathcal{Q} \to \mathcal{Q}
\]
by
\begin{eqnarray*}
\gamma_{j} & := & \textstyle  \theta \circ \left(\left(\left( \sum_{i=1}^{j-1} e_{i}   \right) \otimes \mathrm{id}_{\mathcal{Q}} \right)   \oplus \left( e_{j} \otimes \mathrm{id}_{\mathcal{Q}}   \right)\right. \\
&& \textstyle \oplus \left. \left(\left( \sum_{i=j+1}^{N} e_{i}   \right) \otimes \mathrm{id}_{\mathcal{Q}} \right)   \oplus \left( e_{N+1} \otimes \mathrm{id}_{\mathcal{Q}}   \right)\right)
\end{eqnarray*}
for $j \in \{1,\ldots,N\}$. 

We clearly have 
\[
\tau_{\mathcal{Q}}(\gamma((1_{\mathcal{Q}},0)))  > 1-\delta.
\]
Note also that
\begin{equation}
\label{K3}
\gamma \circ (\bar{\kappa}_{0} \oplus \bar{\kappa}) = \gamma_{N} \circ (\bar{\kappa}_{0} \oplus \bar{\kappa}_{0} \oplus \bar{\kappa}_{1} \oplus \bar{\kappa})
\end{equation}
and
\begin{equation}
\label{K4}
\gamma \circ (\bar{\kappa}_{1} \oplus \bar{\kappa}) = \gamma_{1} \circ (\bar{\kappa}_{0} \oplus \bar{\kappa}_{1} \oplus \bar{\kappa}_{1} \oplus \bar{\kappa}),
\end{equation}
since
\[
\gamma_{N}((0,0,x,0)) = \gamma_{1}((x,0,0,0)) = 0
\]
for $x \in \mathcal{Q}$. 

We furthermore have 
\begin{equation}
\label{K1}
\gamma_{j} \circ (\bar{\kappa}_{0} \oplus \bar{\kappa}_{0} \oplus \bar{\kappa}_{1} \oplus \bar{\kappa}) = \gamma_{j+1} \circ  (\bar{\kappa}_{0} \oplus \bar{\kappa}_{1} \oplus \bar{\kappa}_{1} \oplus \bar{\kappa})
\end{equation}
for $j \in \{1,\ldots,N-1\}$. From \eqref{Kw1} and \ref{H} we obtain
\begin{equation}
\label{K2}
(E,\pi,\sigma, \gamma_{j} \circ (\bar{\kappa}_{0} \oplus \bar{\kappa}_{0} \oplus \bar{\kappa}_{1} \oplus \bar{\kappa}) ) \sim_{(\mathcal{F},\eta)} (E,\pi,\sigma, \gamma_{j} \circ (\bar{\kappa}_{0} \oplus \bar{\kappa}_{1} \oplus \bar{\kappa}_{1} \oplus \bar{\kappa}) ) .
\end{equation}

We now have
\begin{eqnarray*}
(E,\pi,\sigma,\gamma \circ(\bar{\kappa}_{0} \oplus \bar{\kappa}))& \stackrel{\eqref{K3}}{=} & (E,\pi,\sigma, \gamma_{N} \circ (\bar{\kappa}_{0} \oplus \bar{\kappa}_{0} \oplus \bar{\kappa}_{1} \oplus \bar{\kappa}) ) \\
&\stackrel{\eqref{K2}}{\sim_{(\mathcal{F},\eta)}}& (E,\pi,\sigma, \gamma_{N} \circ (\bar{\kappa}_{0} \oplus \bar{\kappa}_{1} \oplus \bar{\kappa}_{1} \oplus \bar{\kappa}) ) \\
& \vdots&  \\
& \stackrel{\eqref{K1}}{=} & (E,\pi,\sigma, \gamma_{j} \circ (\bar{\kappa}_{0} \oplus \bar{\kappa}_{0} \oplus \bar{\kappa}_{1} \oplus \bar{\kappa}) )\\
& \sim_{(\mathcal{F},\eta)} & (E,\pi,\sigma, \gamma_{j} \circ (\bar{\kappa}_{0} \oplus \bar{\kappa}_{1} \oplus \bar{\kappa}_{1} \oplus \bar{\kappa}) )\\
& \vdots & \\
& = & (E,\pi,\sigma, \gamma_{1} \circ (\bar{\kappa}_{0} \oplus \bar{\kappa}_{0} \oplus \bar{\kappa}_{1} \oplus \bar{\kappa}) )\\
& \sim_{(\mathcal{F},\eta)} & (E,\pi,\sigma, \gamma_{1} \circ (\bar{\kappa}_{0} \oplus \bar{\kappa}_{1} \oplus \bar{\kappa}_{1} \oplus \bar{\kappa}) )\\
&  \stackrel{\eqref{K4}}{=} & (E,\pi,\sigma,\gamma \circ(\bar{\kappa}_{1} \oplus \bar{\kappa})).
\end{eqnarray*}
\end{nproof}
\en

\bn
\begin{nprop}
\label{C}
Let $B$ be a separable unital recursive subhomogeneous $\mathrm{C}^{*}$-algebra 
and let $\mathcal{F} \subset B^{1}_{+}$ finite and $0<\eta, \beta \le1$ be given. 

Suppose $B$ has an $(\mathcal{F},\eta)$-connected  recursive subhomogeneous decomposition 
\[
[B_{l},X_{l},\Omega_{l},r_{l},\phi_{l}]_{l=1}^{R}
\]
along which projections can be lifted and such that  $X_{l} \setminus \Omega_{l} \neq \emptyset$ for $l \ge 1$. 

Let $\tau_{0},\tau_{1} \in T(B)$ be tracial states with 
\[
(\tau_{0})_{*} = (\tau_{1})_{*}
\] 
(as states on the ordered $\mathrm{K}_{0}(B)$). 

Then there are $x_{l} \in X_{l} \setminus \Omega_{l}$ for $l \in \{1,\ldots,R\}$ and pairwise orthogonal  $(\mathcal{F},\eta)$-excisors 
\[
(E_{x_{l}}, \pi_{x_{l}}, \sigma_{x_{l}}), \; l \in \{1,\ldots,R\};
\] 
in this case,
\[
E_{x_{l}} \cong M_{r_{l}}, \; l \in \{1,\ldots,R\}
\]
and
\begin{equation}
\label{Cw3}
\textstyle
\left( E:= \bigoplus_{l=1}^{R} E_{x_{l}}, \; \pi:= \bigoplus_{l=1}^{R} \pi_{x_{l}}, \; \sigma:= \bigoplus_{l=1}^{R} \sigma_{x_{l}} \right)
\end{equation}
is an $(\mathcal{F},\eta)$-excisor. 

Furthermore, there are unital embeddings 
\[
\kappa_{i}: E \to \mathcal{Q}, \; i \in \{0,1\},
\]
such that 
\[
(E,\pi,\sigma,\kappa_{0}) \sim_{(\mathcal{F},\eta)} (E,\pi,\sigma,\kappa_{1})
\]
and such that
\begin{equation}
\label{Cw4}
y_{i,l}:= \tau_{\mathcal{Q}}(\kappa_{i}(1_{E_{x_{l}}}))
\end{equation}
satisfy
\begin{equation}
\label{C1}
|y_{i,l} - y^{\tau_{i}}_{l}| < \beta
\end{equation}
for $i \in \{0,1\}$, $l \in \{1,\ldots,R\}$, where the $y_{l}^{\tau_{i}}$ are as in \eqref{local-trace-notationw1}. 
\end{nprop}

\begin{nproof}
Apply Proposition~\ref{J} with
\[
\bar{\beta}:= \frac{\beta}{3}
\]
to obtain $x_{l} \in X_{l} \setminus \Omega_{l}$, pairwise orthogonal weighted $(\mathcal{F},\eta)$-excisors 
\[
(E_{x_{l}},\pi_{x_{l}},\sigma_{x_{l}},\kappa_{x_{l}}), \; l \in \{1,\ldots,R\},
\]
and unital embeddings
\[
\textstyle
\gamma_{0},\gamma_{1},\tilde{\gamma}: \bigoplus_{l=1}^{R} \mathcal{Q} \to \mathcal{Q}
\]
and
\[
\bar{\gamma}: \mathcal{Q} \oplus \mathcal{Q} \to \mathcal{Q}.
\]
Apply Proposition~\ref{K} with 
\[
\delta:= \frac{\beta}{3}
\]
and with
\[
\textstyle
\bar{\kappa}_{i}:= \gamma_{i} \circ \left( \bigoplus_{l=1}^{R} \kappa_{x_{l}}  \right), \; i \in \{0,1\},
\]
and
\[
\textstyle
\bar{\kappa}:= \tilde{\gamma} \circ \left( \bigoplus_{l=1}^{R} \kappa_{x_{l}}  \right)
\]
to obtain a unital embedding 
\[
\gamma: \mathcal{Q} \oplus \mathcal{Q} \to \mathcal{Q}
\]
such that 
\begin{equation}
\label{Cw1}
|\tau_{\mathcal{Q}}(\gamma((p,0))) - \tau_{\mathcal{Q}}( p)| < \frac{\beta}{3}
\end{equation}
for every projection $p \in \mathcal{Q}$, whence in particular 
\begin{equation}
\label{Cw2}
\tau_{\mathcal{Q}}(\gamma((0,1_{\mathcal{Q}}))) < \frac{\beta}{3},
\end{equation}
and such that
\[
\kappa_{i}:= \gamma \circ (\bar{\kappa}_{i} \oplus \bar{\kappa}),\; i \in \{0,1\}
\]
satisfy
\[
(E,\pi,\sigma,\kappa_{0}) \sim_{(\mathcal{F},\eta)} (E,\pi,\sigma,\kappa_{1}).
\]
With
\[
y_{i,l}= \tau_{\mathcal{Q}}(\kappa_{i}(1_{E_{x_{l}}})) = \tau_{\mathcal{Q}}(\gamma(\gamma_{i}(\kappa_{x_{l}}(1_{E_{x_{l}}}))  \oplus \tilde{\gamma} (\kappa_{x_{l}}(1_{E_{x_{l}}}))))
\]
and 
\[
\bar{y}_{i,l}:= \tau_{\mathcal{Q}} (\gamma_{i}(\kappa_{x_{l}}(1_{E_{x_{l}}})))
\]
we have 
\begin{eqnarray*}
|y_{i,l} - y^{\tau_{i}}_{l}| & \le & |y_{i,l} -  \bar{y}_{i,l}| + | \bar{y}_{i,l} - y^{\tau_{i}}_{l}|  \\
& \le & |\tau_{\mathcal{Q}} (\gamma(\gamma_{i}(\kappa_{x_{l}} (1_{E_{x_{l}}})) \oplus 0 )) - \tau_{\mathcal{Q}}(\gamma_{i}(\kappa_{x_{l}}(1_{E_{x_{l}}})))| \\
&& + \tau_{\mathcal{Q}} (\gamma(0 \oplus \tilde{\gamma }(\kappa_{x_{l}}(1_{E_{x_{l}}})))) \\
&& + | \bar{y}_{i,l} - y^{\tau_{i}}_{l} | \\
& \stackrel{\eqref{Cw1},\eqref{Cw2},\eqref{Jw24}}{\le} & \frac{\beta}{3} +  \frac{\beta}{3} +  \frac{\beta}{3}.
\end{eqnarray*}
\end{nproof}
\en

\bn
\begin{nlemma} 
\label{E}
Let $B$ be a separable unital recursive subhomogeneous $\mathrm{C}^{*}$-algebra 
and let $\mathcal{F} \subset B^{1}_{+}$ finite and $\eta >0$ be given. 

Suppose $B$ has an $(\mathcal{F},\eta)$-connected  recursive subhomogeneous decomposition 
\[
[B_{l},X_{l},\Omega_{l},r_{l},\phi_{l}]_{l=1}^{R}
\]
along which projections can be lifted and such that  $X_{l} \setminus \Omega_{l} \neq \emptyset$ for $l \ge 1$. 

Let $\tau^{(0)},\ldots,\tau^{(n-1)} \in T(B)$ be $n$ faithful tracial states with 
\[
(\tau^{(0)})_{*} =\ldots =  (\tau^{(n-1)})_{*}
\] 
(as states on the ordered $\mathrm{K}_{0}(B)$).

Then there are 
\[
0=K_{0} < K_{1 }< \ldots < K_{n-1} = K \in \mathbb{N}
\]
and pairwise orthogonal weighted $(\mathcal{F},\eta)$-excisors 
\[
(Q_{\frac{j}{K}},\rho_{\frac{j}{K}},\sigma_{\frac{j}{K}},\kappa_{\frac{j}{K}}), \; j \in \{0,\ldots,K\},
\]
implementing $(\mathcal{F},\eta)$-bridges 
\begin{eqnarray*}
(Q_{\frac{K_{0}}{K}},\rho_{\frac{K_{0}}{K}},\sigma_{\frac{K_{0}}{K}},\kappa_{\frac{K_{0}}{K}}) &  \sim_{(\mathcal{F},\eta)}  \ldots \sim_{(\mathcal{F},\eta)} & (Q_{\frac{K_{m}}{K}},\rho_{\frac{K_{m}}{K}},\sigma_{\frac{K_{m}}{K}},\kappa_{\frac{K_{m}}{K}}) \\ 
& \sim_{(\mathcal{F},\eta)} \ldots \sim_{(\mathcal{F},\eta)} &  (Q_{\frac{K_{n-1}}{K}},\rho_{\frac{K_{n-1}}{K}},\sigma_{\frac{K_{n-1}}{K}},\kappa_{\frac{K_{n-1}}{K}}) ,
\end{eqnarray*}
and such that, for each projection $q \in Q_{\frac{K_{m}}{K}}$, $m \in \{0,\ldots,n-1\}$,
\begin{equation}
\label{E9}
(\tau^{(m)} \otimes \tau_{\mathcal{Q}}) \sigma_{\frac{K_{m}}{K}}(q) \ge \frac{1}{n+1} \cdot \tau_{\mathcal{Q}} \kappa_{\frac{K_{m}}{K}}(q).
\end{equation}
\end{nlemma}

\begin{nproof}
Let us first prove the lemma for $n = 2$. Choose
\[
0< \bar{\alpha}, \beta, \delta <\frac{1}{n}
\]
such that
\begin{equation}
\label{E7}
\left( \frac{1}{n} - \delta   \right) \cdot (1-2\bar{\alpha}) \ge \frac{1}{n+1}
\end{equation}
and
\begin{equation}
\label{E1}
\beta < \bar{\alpha} \cdot \frac{y_{l}^{\tau^{(i)}}}{4} \mbox{ for all } i \in \{0,1\},\, l \in \{1,\ldots,R\}
\end{equation}
(this is possible since $X_{l} \setminus \Omega_{l} \neq \emptyset$ and the traces are faithful, whence $y_{l}^{\tau^{(i)}}>0$). 

Let $(E,\pi,\sigma,\kappa_{i})$ and $y_{i,l} = \tau_{\mathcal{Q}(\kappa_{i}(1_{E_{x_{l}}}))}$ for $i \in \{0,1\}$, $l\in \{1,\ldots,R\}$, be as in Proposition~\ref{C}. 

Choose
\[
0<\gamma<\beta,
\]
then
\begin{eqnarray}
y_{i,l} - \gamma - \beta \ge y_{i,l} - 2 \beta  & \stackrel{\eqref{E1}}{\ge} & y_{i,l} - \bar{\alpha}  \cdot \frac{y_{l}^{\tau^{(i)}}}{2} \nonumber \\
& \stackrel{\eqref{E1}}{\ge} &  y_{i,l} - \bar{\alpha}  \cdot (y_{l}^{\tau^{(i)}} - \beta) \nonumber  \\
&  \stackrel{\eqref{C1}}{\ge} & (1-\bar{\alpha}) \cdot y_{i,l}. \label{Ew1}
\end{eqnarray}
By Proposition~\ref{A}, there are $(\mathcal{F},\eta)$-excisors 
\[
(\dot{E}_{i},\dot{\rho}_{i},\dot{\sigma}_{i}),\; i \in \{0,1\},
\]
compatible with the recursive subhomogeneous decomposition, with
\[
 \textstyle
\dot{E}_{i} = \bigoplus_{l=1}^{R} \dot{E}_{i,l}
\]
and each $\dot{E}_{i,l}$ a direct sum of copies of $M_{r_{l}}$, cf.\ \eqref{Aw1},  and such that 
\begin{eqnarray}
\lefteqn{(\bar{\tau}_{i,l} \otimes \tau_{\mathcal{Q}}) \circ (\psi_{l} \otimes \mathrm{id}_{\mathcal{Q}}) \circ \dot{\sigma}_{i,l}(1_{\dot{E}_{i,l}})} \nonumber \\
 & \ge & y_{l}^{\tau^{(i)}} - \gamma \nonumber \\
&\stackrel{\eqref{C1}}{\ge} & y_{i,l} - \gamma - \beta \nonumber \\
& \stackrel{\eqref{Ew1}}{\ge} & (1-\bar{\alpha}) \cdot y_{i,l}. \label{Ew2}
\end{eqnarray}

Choose unital $^{*}$-homomorphisms
\[
\dot{\kappa}_{i}:\dot{E}_{i} \to \mathcal{Q}, \; i \in \{0,1\},
\]
such that
\[
\tau_{\mathcal{Q}} \circ \dot{\kappa}_{i}(1_{\dot{E}_{i,l}}) = y_{i,l}
\]
and
\[
(\tau^{(i)} \otimes \tau_{\mathcal{Q}}) \circ \dot{\sigma}_{i}(q) \ge (\bar{\tau}_{i,l} \otimes \tau_{\mathcal{Q}}) \circ (\psi_{l} \otimes \mathrm{id}_{\mathcal{Q}}) \circ \dot{\sigma}_{i,l} (q) \ge (1-\bar{\alpha}) \cdot \tau_{\mathcal{Q}} \circ \dot{\kappa}_{i}(q)
\]
for all projections $q \in \dot{E}_{i,l}$; it follows that
\begin{equation}
\label{Ew5}
(\tau^{(i)} \otimes \tau_{\mathcal{Q}}) \circ \dot{\sigma}_{i}(q) \ge (1-\bar{\alpha}) \cdot \tau_{\mathcal{Q}} \circ \dot{\kappa}_{i}(q)
\end{equation}
for all projections $q \in \dot{E}_{i}$.

Now by Proposition~\ref{D}, we have
\begin{equation}
\label{Ew3}
(\dot{E}_{i}, \dot{\rho}_{i}, \dot{\sigma}_{i}, \dot{\kappa}_{i}) \sim_{(\mathcal{F},\eta)} (E,\pi,\sigma,\kappa_{i})
\end{equation}
for $i \in \{0,1\}$. By Proposition~\ref{C}, 
\[
(E,\pi,\sigma,\kappa_{0}) \sim_{(\mathcal{F},\eta)} (E,\pi,\sigma,\kappa_{1}),
\]
so by transitivity,
\begin{equation}
\label{Ew4}
(\dot{E}_{0}, \dot{\rho}_{0}, \dot{\sigma}_{0}, \dot{\kappa}_{0}) \sim_{(\mathcal{F},\eta)} (\dot{E}_{1}, \dot{\rho}_{1}, \dot{\sigma}_{1}, \dot{\kappa}_{1}),
\end{equation}
with a bridge consisting of $(\mathcal{F},\eta)$-excisors
\[
(\dot{E}_{\frac{j}{K}}, \dot{\rho}_{\frac{j}{K}}, \dot{\sigma}_{\frac{j}{K}}, \dot{\kappa}_{\frac{j}{K}}),\; j \in \{0,\ldots, K\}, 
\]
for some $K \in \mathbb{N}$.

Choose pairwise orthogonal projections
\[
q_{0},q_{1/K},\ldots,q_{1} \in \mathcal{Q}
\]
such that
\[
\textstyle
\sum_{j=0}^{K} q_{\frac{j}{K}} = 1_{\mathcal{Q}}
\]
and
\[
\tau_{\mathcal{Q}}(q_{0}) = \tau_{\mathcal{Q}}(q_{1}) = \frac{1}{2} - \delta
\]
and
\[
\tau_{\mathcal{Q}}(q_{j/K}) = \frac{2\delta}{K-1}, \; j \in \{1,\ldots,K-1\}.
\]
Choose a $^{*}$-isomorphism
\[
\theta:\mathcal{Q} \otimes \mathcal{Q} \to \mathcal{Q}
\]
and define, for $j \in \{0,\ldots,K\}$,
\begin{eqnarray*}
Q_{\frac{j}{K}}& := & \dot{E}_{\frac{j}{K}}, \\
\rho_{\frac{j}{K}}(\, .\,)& := &  \dot{\rho}_{\frac{j}{K}} (\, .\,),\\
\sigma_{\frac{j}{K}}(\, .\,)& := & (\mathrm{id}_{B} \otimes \theta) \circ (\dot{\sigma}_{\frac{j}{K}}(\, .\,) \otimes q_{\frac{j}{K}}),\\
\kappa_{\frac{j}{K}}& := & \dot{\kappa}_{\frac{j}{K}}.
\end{eqnarray*}

We check that the $(Q_{\frac{j}{K}}, \rho_{\frac{j}{K}},\sigma_{\frac{j}{K}},\kappa_{\frac{j}{K}})$, $j \in \{0,\ldots,K\}$, have the right properties:

Each $\sigma_{\frac{j}{K}}$ is an isometric c.p.\ order zero map since $\dot{\sigma}_{\frac{j}{K}}$ is and since $q_{\frac{j}{K}}$ is nonzero. The $\sigma_{\frac{j}{K}}$ have pairwise orthogonal images, since the $q_{\frac{j}{K}}$ are pairwise orthogonal. 

For $i \in \{0,1\}$ and $q \in Q_{i}$ a projection, we have
\begin{eqnarray}
(\tau^{(i)} \otimes \tau_{\mathcal{Q}}) (\sigma_{i}(q)) & = & (\tau^{(i)} \otimes \tau_{\mathcal{Q}}) (\mathrm{id}_{B} \otimes \theta) (\dot{\sigma}_{i} (q) \otimes q_{i}) \nonumber  \\
& = & (\tau^{(i)} \otimes \tau_{\mathcal{Q}} \otimes \tau_{\mathcal{Q}}) ( \dot{\sigma}_{i}(q) \otimes q_{i}) \nonumber \\
& = & \tau_{\mathcal{Q}}(q_{i}) \cdot (\tau^{(i)} \otimes \tau_{\mathcal{Q}}) (\dot{\sigma}_{i}(q)) \nonumber \\
& \ge & (1/2 - \delta) \cdot (1-\bar{\alpha}) \cdot \tau_{\mathcal{Q}} \circ \kappa_{i}(q) \nonumber \\
& \ge & \frac{1}{2+1}  \cdot \tau_{\mathcal{Q}} \circ \kappa_{i}(q) . \label{E10}
\end{eqnarray}

For $j \in \{0,\ldots,K\}$ and $b \in \mathcal{F}$,
\begin{eqnarray*}
\lefteqn{ \| \sigma_{\frac{j}{K}}(1_{Q_{\frac{j}{K}}}) ( b \otimes 1_{\mathcal{Q}}) - \sigma_{\frac{j}{K}} \rho_{\frac{j}{K}}(b) \|  }\\
& = &  \| ((\mathrm{id}_{B}  \otimes \theta ) (\dot{\sigma}_{\frac{j}{K}}(1_{Q_{\frac{j}{K}}})   \otimes q_{\frac{j}{K}} )) ((\mathrm{id}_{B}  \otimes \theta)(b \otimes 1_{\mathcal{Q}}  \otimes 1_{\mathcal{Q}}))\\
&& - (\mathrm{id}_{B} \otimes \theta)(\dot{\sigma}_{\frac{j}{K}} \dot{\rho}_{\frac{j}{K}}(b) \otimes q_{\frac{j}{K}} ) \| \\
& = & \| (\dot{\sigma}_{\frac{j}{K}}(1_{Q_{\frac{j}{K}}})  (b \otimes 1_{\mathcal{Q}})) \otimes q_{\frac{j}{K}} - \dot{\sigma}_{\frac{j}{K}} \dot{\rho}_{\frac{j}{K}} (b) \otimes q_{\frac{j}{K}} \| \\
& < & \eta,
\end{eqnarray*}
so each $(Q_{\frac{j}{K}}, \rho_{\frac{j}{K}},\sigma_{\frac{j}{K}},\kappa_{\frac{j}{K}})$ is a weighted $(\mathcal{F},\eta)$-excisor (which is clearly compatible with the given recursive subhomogeneous decomposition).

\bigskip

We now turn to the case of arbitrary $n$. Fix $(E,\pi,\sigma)$ as in the first part of the proof, cf.\ \eqref{Cw3}. Choose $\bar{\alpha}, \beta, \delta$ as above; we may in addition assume that
\begin{equation}
\label{E8}
(\tau^{(m)} \otimes \tau_{\mathcal{Q}}) \sigma(q) \ge \frac{n \beta}{\delta} 
\end{equation}
for all $m \in \{0,\ldots,n-1\}$ and all nonzero projections $q \in E$. 

We now apply the first part of the proof to each pair $\tau^{(m)}, \tau^{(m+1)}$, $m \in \{0,\ldots, n-2\}$. This yields for each $m  \in \{0,\ldots,n-1\}$ and $i \in \{0,1\}$ $(\mathcal{F},\eta)$-excisors 
\[
(E,\pi,\sigma,\kappa_{i}^{(m)})
\]
and
\[
y^{(m)}_{i,l} \stackrel{\eqref{Cw4}}{=} \tau_{\mathcal{Q}}(\kappa_{i}^{(m)}(1_{M_{r_{l}}})), \, l \in \{1,\ldots,R\},
\]
such that 
\[
|y_{i,l}^{(m)} - y_{l}^{\tau^{(m+i)}}| \stackrel{\eqref{C1}}{<} \beta
\]
as well as $(\mathcal{F},\eta)$-excisors 
\[
(\dot{E}_{i}^{(m)},\dot{\rho}_{i}^{(m)},\dot{\sigma}_{i}^{(m)},\dot{\kappa}_{i}^{(m)})
\]
with
\begin{eqnarray}
(E,\pi,\sigma,\kappa_{0}^{(m)}) & \stackrel{\eqref{Ew3}}{\sim_{(\mathcal{F},\eta)}} & (\dot{E}_{0}^{(m)},\dot{\rho}_{0}^{(m)},\dot{\sigma}_{0}^{(m)},\dot{\kappa}_{0}^{(m)}) \nonumber \\
& \stackrel{\eqref{Ew4}}{\sim_{(\mathcal{F},\eta)}}  & (\dot{E}_{1}^{(m)},\dot{\rho}_{1}^{(m)},\dot{\sigma}_{1}^{(m)},\dot{\kappa}_{1}^{(m)}) \nonumber \\
& \stackrel{\eqref{Ew3}}{\sim_{(\mathcal{F},\eta)}} &  (E,\pi,\sigma,\kappa_{1}^{(m)}) \label{E4}
\end{eqnarray}
and with 
\[
(\tau^{(m+i)} \otimes \tau_{\mathcal{Q}}) \circ \dot{\sigma}_{i}^{(m)} (q) \stackrel{\eqref{Ew5}}{\ge} (1-\bar{\alpha}) \cdot \tau_{\mathcal{Q}} \circ \dot{\kappa}_{i}^{(m)}(q) 
\]
for all projections $q \in \dot{E}_{i}^{(m)}$, $m \in \{0,\ldots,n-2\}$, $i \in \{0,1\}$. 

But then it is not hard to find unital $^{*}$-homomorphisms 
\[
\kappa^{(m)},\hat{\kappa}^{(m)}_{1},\hat{\kappa}^{(m+1)}_{0}: E \to \mathcal{Q}
\]
such that 
\begin{equation}
\label{E2}
\textstyle
\frac{1-n\beta/2}{1-(n-1)\beta/2}\cdot \kappa^{(m)} \oplus \frac{\beta/2}{1-(n-1)\beta/2} \cdot \hat{\kappa}_{1}^{(m)} \approx_{\mathrm{u}}  \kappa_{1}^{(m)}
\end{equation}
and
\begin{equation}
\label{E3}
\textstyle
\frac{1-n\beta/2}{1-(n-1)\beta/2}\cdot \kappa^{(m)} \oplus \frac{\beta/2}{1-(n-1)\beta/2} \cdot \hat{\kappa}_{0}^{(m+1)} \approx_{\mathrm{u}}  \kappa_{0}^{(m+1)}
\end{equation}
(here, we use notation as in \ref{Ia}(ii) to denote weighted sums of $^{*}$-homomorphisms $E \to \mathcal{Q}$).

Combining \eqref{E4}, \eqref{E2} and \eqref{E3} with Remark~\ref{M}(ii), one checks  that 
\begin{eqnarray}
\lefteqn{\textstyle (1-(n-1) \frac{\beta}{2})\cdot (E,\pi,\sigma,\kappa_{0}^{(m)}) }\nonumber  \\
\lefteqn{\textstyle \oplus \left( \bigoplus_{m' \in \{0,\ldots,n-1\} \setminus \{m\}}   \frac{\beta}{2} \cdot (E,\pi,\sigma,\hat{\kappa}_{0}^{(m')}) \right) } \nonumber \\
& \sim_{(\mathcal{F},\eta)} & \textstyle (1-(n-1)\frac{\beta}{2})\cdot (E,\pi,\sigma,\kappa_{0}^{(m+1)}) \nonumber \\
&& \textstyle \oplus \left( \bigoplus_{m' \in \{0,\ldots,n-1\} \setminus \{m+1\}}   \frac{\beta}{2} \cdot (E,\pi,\sigma,\hat{\kappa}_{0}^{(m')}) \right) .   \label{E5}
\end{eqnarray}

Combining \eqref{E4} with \eqref{E5} we see that, for all $m \in \{0,\ldots,n-2\}$, 
\begin{eqnarray}
\lefteqn{\textstyle (1-(n-1) \frac{\beta}{2})\cdot (\dot{E}_{0}^{(m)},\dot{\rho}_{0}^{(m)},\dot{\sigma}^{(m)}_{0},\dot{\kappa}_{0}^{(m)}) }\nonumber  \\
\lefteqn{\textstyle \oplus \left( \bigoplus_{m' \in \{0,\ldots,n-1\} \setminus \{m\}}   \frac{\beta}{2} \cdot (E,\pi,\sigma,\hat{\kappa}_{0}^{(m')}) \right) } \nonumber \\
& \sim_{(\mathcal{F},\eta)} & \textstyle (1-(n-1)\frac{\beta}{2})\cdot (\dot{E}_{0}^{(m+1)},\dot{\rho}_{0}^{(m+1)},\dot{\sigma}^{(m+1)}_{0},\dot{\kappa}_{0}^{(m+1)}) \nonumber \\
&& \textstyle \oplus \left( \bigoplus_{m' \in \{0,\ldots,n-1\} \setminus \{m+1\}}   \frac{\beta}{2} \cdot (E,\pi,\sigma,\hat{\kappa}_{0}^{(m')}) \right) .   \label{E6}
\end{eqnarray}
Note that, for any projection $q \in \dot{E}^{(m)}_{i}$,
\begin{eqnarray*}
\lefteqn{\textstyle (1- (n-1) \frac{\beta}{2}) \cdot (\tau^{(m+i)} \otimes \tau_{\mathcal{Q}}) \circ \dot{\sigma}_{i}^{(m)}(q) }\\
& \ge &  \textstyle (1- (n-1) \frac{\beta}{2})  (1-\bar{\alpha}) \cdot  \tau_{\mathcal{Q}} \circ \dot{\kappa}^{(m)}_{i}(q) \\
& \ge & (1 - 2 \bar{\alpha}) \cdot \tau_{\mathcal{Q}} \circ \dot{\kappa}^{(m)}_{i} (q).
\end{eqnarray*}

We may therefore assume that there are
\[
0=K_{0}<  K_{1} < \ldots < K_{n-1} = K \in \mathbb{N}
\]
and an $(\mathcal{F},\eta)$-bridge consisting of $(\mathcal{F},\eta)$-excisors 
\[
(Q_{\frac{j}{K}}, \rho_{\frac{j}{K}},\bar{\sigma}_{\frac{j}{K}},\kappa_{\frac{j}{K}}), \, j \in \{0,\ldots,K\}
\]
with 
\begin{eqnarray*}
\lefteqn{(Q_{\frac{K_{m}}{K}}, \rho_{\frac{K_{m}}{K}},\bar{\sigma}_{\frac{K_{m}}{K}},\kappa_{\frac{K_{m}}{K}}) }\\
& = & \textstyle (1-(n-1) \frac{\beta}{2})\cdot (\dot{E}_{0}^{(m)},\dot{\rho}_{0}^{(m)},\dot{\sigma}^{(m)}_{0},\dot{\kappa}_{0}^{(m)}) \nonumber  \\
&  & \textstyle \oplus \left( \bigoplus_{m' \in \{0,\ldots,n-1\} \setminus \{m\}}   \frac{\beta}{2} \cdot (E,\pi,\sigma,\hat{\kappa}_{0}^{(m')}) \right)  \nonumber 
\end{eqnarray*}
for $m \in \{0,\ldots,n-1\}$. 

Choose pairwise orthogonal projections 
\[
q_{0}, q_{\frac{1}{K}},\ldots, q_{1} \in \mathcal{Q}
\]
such that
\[
\textstyle
\sum_{j} q_{\frac{j}{K}} = 1_{\mathcal{Q}}
\]
and such that each $q_{\frac{K_{m}}{K}}$ can be written as a sum of two projections 
\[
q_{\frac{K_{m}}{K}} = q'_{\frac{K_{m}}{K}} + q''_{\frac{K_{m}}{K}}
\]
with 
\[
\tau_{\mathcal{Q}}(q'_{\frac{K_{m}}{K}}) = 1/n - \delta,
\]
\[
\tau_{\mathcal{Q}}(q''_{\frac{K_{m}}{K}}) = \delta/2
\]
and such that all other projections have the same tracial value $n \delta/2K$.

As in the first part of the proof, choose an isomorphism 
\[
\theta: \mathcal{Q} \otimes \mathcal{Q} \to \mathcal{Q}
\]
and define 
\[
\sigma_{\frac{j}{K}} := (\mathrm{id}_{B} \otimes \theta) \circ (\bar{\sigma}_{\frac{j}{K}}(\,.\,) \otimes q_{\frac{j}{K}}).
\]
Then the
\[
(Q_{\frac{j}{K}}, \rho_{\frac{j}{K}},\sigma_{\frac{j}{K}},\kappa_{\frac{j}{K}})
\]
clearly are  $(\mathcal{F},\eta)$-excisors implementing an $(\mathcal{F},\eta)$-bridge. \eqref{E9} is now checked in a similar manner as \eqref{E10}, using \eqref{E7} and \eqref{E8}.
\end{nproof}
\en

%%%%%%%%%%%%%%%%%%%%%%%%%%%%%%%%%%%%%%%%%%%%%%%%
\section{Tracially large intervals} \label{intervals}

\noindent
Let $a, b \in A_+$. Recall that the element $a$ is Cuntz subequivalent to $b$, written $a \precsim b$, if there is a sequence $(z_n)_{n \in \mathbb{N}} \subset A$ such that $\|z_n^*bz_n - a\| \to 0$ as $n \to \infty$.

For any tracial state $\tau$ we may define a dimension function on the positive elements of $M_{\infty}(A)$ by
$$d_{\tau}(a) = \lim_{n \to \infty} \tau(a^{1/n}).$$
We say that a simple unital $\mathrm{C}^{*}$-algebra $A$ has strict comparison (of positive elements) if $d_{\tau}(a) < d_{\tau}(b)$ for all $\tau \in T(A)$ implies that $a \precsim b$.  

The technical foundation for our main result, Theorem~\ref{TAI}, is laid out in Theorem~\ref{main result}. There we must find an interval in the $\mathrm{C}^{*}$-algebra that is large on all traces and can be moved into position under the ``discrete'' version of the interval that will come from the $(\F, \eta)$-bridges of the previous sections. The interval is twisted into place using a partial isometry obtained from strict comparison. To do this, we will tracially match the endpoints of the $(\F, \eta)$-bridges to functions in a partition of unity. This requires that the finitely many traces be separated along the interval.  The next lemma shows that we can find an interval with each trace almost concentrated at distinct points.

In Proposition~\ref{get s}, we find a positive element which acts as an ``almost'' partial isometry which takes order zero maps to order zero maps. In Theorem~\ref{main result} such an element will be perturbed into an honest partial isometry (dependent on the finite subset $\F \subset A \otimes \Q$ and $\epsilon> 0$) and its support projection will be the unit for the approximating interval algebra. Proposition~\ref{get s} below will furnish this unit with the appropriate properties so that (i) and (ii) of Definition~\ref{TAS} are satisfied. 

The next two lemmas show how to find an interval that is large on all traces. In Lemma~\ref{tracially orthogonal} we follow the techniques of Kishimoto in \cite[Theorem 4.5]{Kis:RohlinUHF}, Matui and Sato in \cite[Lemma 3.2]{MatSat:SC_Z-absorb} and Toms, White and the second author in \cite[Lemma 3.4]{TomWhiWin:findim_Trace} to move positive contractions of given tracial sizes that are approximately tracially orthogonal to positive contractions which are norm orthogonal and remain approximately the same tracial size as the original elements. In Lemma~\ref{interval map} we line up pairwise orthogonal elements, which, using Lemma~\ref{tracially orthogonal}, can be of a specified tracial size, in such a way as to generate an interval. 

\bn
\begin{nlemma} \label{tracially orthogonal}
For every $\epsilon > 0$ and every $k \in \mathbb{N}$ there is $\delta > 0$ such that if $A$ is a separable unital $\mathrm{C}^*$-algebra with $T(A) \neq \emptyset$ and $a_0, \dots, a_k \in A$ are positive contractions satisfying $\tau(a_i a_{i'}) < \delta \text{ for all } \tau \in T(A), i \neq i',$ then there are pairwise orthogonal positive contractions $b_0, \dots, b_k \in A$ satisfying $0 <  \tau(a_i - b_i) < \epsilon \text{ for all } \tau \in T(A).$
\end{nlemma}

\begin{nproof} 
First of all, there is a $0 < \delta_0 < 1$ such that if $A$ is a $\mathrm{C}^*$-algebra and $e_0, \dots, e_k \in A_+$ are contractions satisfying  $\|e_i e_{i'} \| < \delta_0$ when $i \neq i' \in \{0, \dots, k\}$ then there are contractions $\tilde{e}_0, \dots, \tilde{e}_k \in A_+$ such that $\| \tilde{e}_i - e_{i} \| < \epsilon/2$ and $\tilde{e}_i \tilde{e}_{i'} = 0$ for every $i \neq i' \in \{0, \dots, k\}$ \cite[Lemma 2.5.15]{Lin:amenable}. 

Define a continuous function $f_{\delta_0} : (0 ,\infty] \to [0,1]$ by
\[ f_{\delta_0}(t) = \min (1 , \frac{t}{\delta_0}).\]
Note that $(1 - f_{\delta_0}(t))t \leq \delta_0$ for every $t \geq 0$.

Let $A$ be a separable unital $\mathrm{C}^*$-algebra with nonempty tracial state space $T(A)$. Choose 
\[ 0< \delta < \frac{\epsilon \cdot \delta_0}{2 k} \]
and suppose that $a_0, \dots, a_k \in A$ are positive contractions with $\tau(a_i a_{i'}) < \delta$ for every $\tau \in T(A)$ whenever $i \neq i' \in \{0, \dots, k\}$.

For each $i \in \{0, \dots, k\}$ define
\begin{equation}
\textstyle g_i = a_i^{1/2} \left(\sum_{\substack{i' \in \{0, \dots, k\} \\ i' \neq i}} a_{i'} \right) a_i^{1/2}.\label{g_i def}
\end{equation}
Then
\[ \textstyle \tau(g_i)  = \tau\left(\sum_{\substack{i' \in \{0, \dots, k\} \\ i' \neq i}} a_i a_{i'}\right) < k  \delta < \frac{\epsilon \cdot \delta_0}{2}\]
for every $\tau \in T(A)$.

For each $i \in \{0, \dots, k\}$ define positive contractions 
\begin{equation}
\label{x_i def} x_{i} = a_i^{1/2}(1 - f_{\delta_0}(g_i)) a_i^{1/2}.
\end{equation}

We have that  $f_{\delta_0}(t) \leq \frac{t}{\delta_0}$ for every $t \in [0, k-1]$ from which it follows that 
\begin{eqnarray*}
0 \leq \tau(a_i - x_{i}) &=& \tau(a_i^{1/2} f_{\delta_0}(g_i) a_i^{1/2}) \\
&\leq& \tau(f_{\delta_0}(g_i)) \\
&\leq& \frac{ \tau(g_i)}{\delta_0} \\
&<& \frac{\epsilon}{2},
\end{eqnarray*}
for every $\tau \in T(A)$.

We compute
\begin{eqnarray*}
\|x_{i'} x_i \|^2 &=& \| x_i x_{i'}^2 x_i \| \\
&\leq& \textstyle \| x_i (\sum_{\substack{j \in \{0, \dots, k\} \\ j \neq i}} x_j) x_i \| \\
&\stackrel{\eqref{x_i def}}{=}& \textstyle \|a_i^{1/2}(1- f_{\delta_0}(g_i)) a_i^{1/2} (\sum_{\substack{j \in \{0, \dots, k\} \\ j \neq i}} x_j) a_i^{1/2}(1 - f_{\delta_0}(g_i))a_i^{1/2}\|\\
&\stackrel{\eqref{x_i def}}{\leq}& \textstyle \|a_i^{1/2}(1- f_{\delta_0}(g_i)) a_i^{1/2} (\sum_{\substack{j \in \{0, \dots, k\} \\ j \neq i}} a_j^{1/2}(1- f_{\delta_0}(g_j))a_j^{1/2}) a_i^{1/2}\\
&& (1 - f_{\delta_0}(g_i))a_i^{1/2}\|\\
&\leq& \textstyle \|a_i^{1/2}(1- f_{\delta_0}(g_i))a_i^{1/2} (\sum_{\substack{j \in \{0, \dots, k\} \\ j \neq i}}a_j)a_i^{1/2}(1 - f_{\delta_0}(g_i))a_i^{1/2} \|\\
&\stackrel{\eqref{g_i def}}{=}& \|a_i^{1/2}(1- f_{\delta_0}(g_i))g_i(1 - f_{\delta_0}(g_i))a_i^{1/2} \|\\
&\leq& \|(1 - f_{\delta_0}(g_i) )g_i \|\\
&\leq& \delta_0.
\end{eqnarray*}

By the choice of $\delta_0$ there are $b_0, \dots, b_1 \in A$ pairwise orthogonal positive contractions with $\| b_i - x_i \| < \epsilon/2$. Thus 
\begin{eqnarray*}
\tau(a_i - b_i) &=& \tau(a_i - x_i) + \tau(x_i - b_i) \\
&<& \epsilon/2 + \| b_i - x_i \| \cdot \tau(1_A) \\
&<& \epsilon.
\end{eqnarray*}
\end{nproof}
\en

\bn For $0 \leq \beta_1  < \beta_2 \leq 1$, define functions $f_{\beta_1, \beta_2}$ and $g_{\beta_1, \beta_2} \in C_0((0,1])$ by

 $$ f_{\beta_1, \beta_2}(t)= \left\{
\begin{array}{cl}
0, & 0 \leq t \leq \beta_1 \\
\text{linear,} & \beta_1 \leq t \leq \beta_2\\
t, &  \beta_2 \leq t \leq 1.
\end{array} \right. $$

 $$ g_{\beta_1, \beta_2}(t)= \left\{
\begin{array}{cl}
0, & 0 \leq t \leq \beta_1 \\
\text{linear,} & \beta_1 \leq t \leq \beta_2\\
1, &  \beta_2 \leq t \leq 1.
\end{array} \right. $$

Note that  if $\beta_1 < \beta_2 < \beta_3 \leq 1$ then 
\begin{equation} \label{gf=f} 
g_{\beta_1, \beta_2} f_{\beta_2, \beta_3} = f_{\beta_2, \beta_3} g_{\beta_1, \beta_2}  = f_{\beta_2, \beta_3}.
\end{equation}

\begin{nlemma} 
\label{interval map}
Let $A$ be a separable simple unital nuclear $\mathrm{C}^{*}$-algebra with exactly $n > 0$ extreme tracial states $\tau_0, \dots, \tau_{n-1} \in T(A)$.  For $i \in \{0 , \dots, n-1\}$, define continuous functions on $[0,1]$ by
 
 $$\gamma_i(t) = \left\{ \begin{array}{ccl} 0, & \, & t \in [0,1]\cap ((-\infty, \frac{i-1}{n-1}] \cup [\frac{i+1}{n-1}, \infty))\\  1, & \, & t = \frac{i}{n-1} \\ \text{ linear, } && \text{elsewhere.} \end{array} \right. $$
 
Then for any $\delta > 0$ there is a $^{*}$-homomorphism 
$$\phi : \mathcal{C}([0,1]) \to A \otimes \Q$$
such that for $i \in \{0, \dots, n-1\}$
$$ \tau_{i} \otimes \tau_{\Q} (\phi(\gamma_i))  \geq 1 - \delta,$$
and
$$ 0 < \tau_{j} \otimes \tau_{\Q} (\phi(\gamma_i)) < \delta$$
when $j \neq i$.
\end{nlemma}

%%%%%%%%%PROOF%%%%%%%%%%%%%%%%%%% PROOF
\begin{nproof}
Choose $0 < \beta < \min (\frac{1}{2}, \frac{\delta}{3})$  and from Lemma~\ref{tracially orthogonal} obtain $\delta_0$ for $n-1$ in place of $k$ and $\epsilon < \frac{\delta}{6}$.

Let $\aff(T(A \otimes \Q))$ denote the set of $\mathbb{R}$-valued bounded affine functions on the tracial state space $T(A \otimes \Q)$.
For each $i \in \{0, \dots, n-1\}$ define continuous functions $\tilde{h}_i$ on the extreme boundary of $T(A \otimes \Q)$ by
$$\tilde{h}_i (\tau_i \otimes \tau_{\Q} ) = 1$$
and
\[ \textstyle 0 < \tilde{h}_i(\tau_j) \leq \min( \frac{\delta}{6}, \delta_0) \text{ when } i \neq j.\]
Since the extreme boundary of $T(A \otimes \Q)$ has only finitely many points and hence is compact, each $\tilde{h}_i$ extends to a continuous affine function $h_i \in \aff(T(A \otimes \Q))$ satisfying $0 < h_i(\tau) \leq 1$ for all $\tau \in T(A \otimes \Q)$ \cite[Theorem II.3.12]{Alf:convex}.

Note that the $h_i$ are not only continuous but are also strictly positive. Since $A$ is simple and unital, by \cite[Corollary 3.10]{BroPerToms:cuntz-semigroup}, there are positive contractions $a_i \in A_+$ satisfying 
\[\tau(a_i) = h_i(\tau) \text{ for all } \tau \in T(A \otimes \Q).\] 
This gives 
\[ \tau(a_i a_{i'} ) < \delta_0 \text{ for all } \tau \in T(A \otimes \Q) \text{ and } i \neq i',\]
whence the previous lemma allows us to obtain pairwise orthogonal positive contractions $y_0, \dots, y_{n-1} \in A\otimes \Q$ such that $\tau_i \otimes \tau_{\Q}(y_i) \geq 1 - \frac{\delta}{3}$ and $\tau_{i'} \otimes \tau_{\Q}(y_i) \leq  \frac{\delta}{3}$ for $i \neq i' \in \{0, \dots, n-1\}$.

Define the following positive elements:
\begin{eqnarray*}
\tilde{b}_{n-1} &=& y_{n-1}\\
b_{n-1} &=& f_{\beta, 2 \beta}(\tilde{b}_{n-1}) \\
\tilde{b}_{n-2} &=& g_{0, \beta}(\tilde{b}_{n-1}) + y_{n-2} \\
b_{n-2} &=&  f_{\beta, 2 \beta}(\tilde{b}_{n-2})\\
&\vdots&\\
\tilde{b}_{1} &=& g_{0, \beta}(\tilde{b}_{2}) + y_1 \\
b_{1} &=& f_{\beta, 2\beta}(\tilde{b}_{1})\\
b_{0} &=& 1. 
\end{eqnarray*} 

Then we have 
\begin{equation} \tag{$\mathcal{R}$}
\|b_i \| \leq 1 \text{ and } b_i b_{i-1} = b_{i-1} b_i = b_i \text{ for every }i \in \{1, \dots, n-1\}. 
\end{equation} Thus we obtain the map
$$ \phi : \mathcal{C}([0,1]) \to \mathrm{C}^*(b_0, \dots, b_{n-1})$$
satisfying
$$\phi(1_{\mathcal{C}([0,1])}) = b_{0} \text{ and }  \phi(g_{\frac{i-1}{n-1}, \frac{i}{n-1}}) = b_i \text{ for } i \in \{1, \dots, n-1\},$$
since $\mathcal{C}([0,1])$  can be written as the universal \mbox{$\mathrm{C}^{*}$-algebra} generated by positive contractions satisfying the relations ($\mathcal{R}$).
 
For each $i = 1, \dots, n-2$ we have that
$$\gamma_i = g_{\frac{i-1}{n-1}, \frac{i}{n-1}} - g_{\frac{i}{n-1} , \frac{i+1}{n-1}}$$
and also
$$\gamma_0 = 1 - g_{0, \frac{1}{n-1}}, \qquad \gamma_{n-1} = g_{\frac{n-2}{n-1}, 1}.$$
Thus $\phi(\gamma_i) = b_i - b_{i+1}$.

We note that 
\begin{eqnarray*}
\tau_i \otimes \tau_{\Q}(b_{i+1}) &\leq&  \tau_i\otimes \tau_{\Q}(\tilde{b}_{i+1}) \\
&=&  \tau_i\otimes \tau_{\Q}(g_{0, \beta}(\tilde{b}_{i+2})) + \tau_i\otimes \tau_{\Q}( y_{i+1})\\
&\leq& \textstyle \frac{1}{\beta} \tau_i\otimes \tau_{\Q}(\tilde{b}_{i+2}) + \tau_i\otimes \tau_{\Q}( y_{i+1})\\
%&\leq& \textstyle \frac{1}{\beta^2}  \tau_i( \tilde{b}_{i+3}) + \frac{1}{\beta} \tau( y_{i+2}) + \tau(y_{i+1})\\
&\vdots&\\
%&\leq& \textstyle \frac{1}{\beta^{n-i-1}} \tau_i(y_{n-1}) + \frac{1}{\beta^2} \tau_i(y_{n-2})  + \cdots + \frac{1}{\beta}  \tau_i(y_{i+2}) + \tau_i(y_{i+1} )\\
&\leq& \textstyle \frac{1}{\beta^{n-i-1}} \cdot \frac{\delta}{3} + \frac{1}{\beta^{n-i-2}} \cdot \frac{\delta}{3}  + \cdots + \frac{1}{\beta} \cdot\frac{\delta}{3} + \frac{\delta}{3}\\
&=&\textstyle \frac{\delta}{3} \cdot(1 - \frac{1}{\beta^{n-i}})/(1 - \frac{1}{\beta})\\
&<& \textstyle \frac{\delta}{3}.
\end{eqnarray*}

It follows that
\begin{eqnarray*}
\tau_i \otimes \tau_{\Q} (\phi(\gamma_i)) &=& \tau_i \otimes \tau_{\Q} (b_i) - \tau_i \otimes \tau_{\Q}  (b_{i+1}) \\
&\geq& \textstyle \tau_i \otimes \tau_{\Q} (y_i) - \beta  - \frac{\delta}{3}\\
&\geq& 1 - \delta.
\end{eqnarray*}

Since $\sum_{j=0}^{n-1} \gamma_j  = 1$, whenever $j \neq i$ we get 
$$\gamma_j \leq 1 - \gamma_i,$$
whence
\[
\tau_i \otimes \tau_{\Q}(\phi(\gamma_j)) \leq 1 - \tau_i \otimes \tau_{\Q}(\phi(\gamma_i)) \leq \delta.
\]
\end{nproof}
\en

%%%%%%%%%%%%%%%%%%%% getting s

\bn
Recall that if $F$ and $A$ are separable $C^*$-algebras with $F$ unital and $\sigma : F \to A$ is a c.p.\ order zero map, we may define a functional calculus for $\sigma$ as follows. Let $\pi_{\sigma}$ denote the supporting $^*$-homomorphism of $\sigma$. Then for $f\in C([0,1])$, we define $f(\sigma)(x) = f(\sigma(1_F))\pi_{\sigma}(x)$, and $f(\sigma)$ is a well-defined c.p.\ order zero map \cite[Corollary 3.2]{WinZac:order-zero}, \cite[1.3]{Winter:dr-Z-stable}.

\begin{nprop} \label{get s}
Let $A$ be a separable simple unital nuclear $\mathrm{C}^{*}$-algebra with stable rank one and strict comparison. Let $F$ be a finite dimensional $\mathrm{C}^{*}$-algebra. Let $0 < \alpha, \epsilon < 1$  and suppose
$$\theta, \sigma: F \to A$$
are c.p.\ order zero maps satisfying
$$ \tau(\sigma(p))  -  d_{\tau}(\theta(p)) \geq \alpha$$
for every nonzero projection $p \in F$ and for every $\tau \in T(A)$. Then, for $0 < \beta_1 < \alpha/2$, there exists $s \in A$ satisfying, with $\beta_1 < \beta_2 < 1$, the following:
\begin{enumerate}
\item  $s^*s \in \her(f_{\beta_1, \beta_2}(\sigma(1_F)))$, \label{s*s location}
\item$(\theta(x) - \epsilon)_+ s s^* = s s^*  (\theta(x) - \epsilon)_+ =  (\theta(x) - \epsilon)_+$ for every $x \in F$ \label{theta commutes}
\item $s^* (\theta(x)- \epsilon)_+ s = g_{0, \beta_1}(\sigma)(x)s^* (\theta(1_F) - \epsilon)_+ s = s^* (\theta(1_F) - \epsilon)_+ s g_{0, \beta_1}(\sigma)(x)$ for every $x \in F$. \label{g(sigma)}
\end{enumerate}
\end{nprop}

\begin{nproof} 
Let $\epsilon$ and  $\alpha$ be given and let $\theta, \sigma : F \to A$ be c.p.\ order zero maps satisfying the statement of the proposition. Denote the supporting $^{*}$-homomorphisms for the c.p.\ order zero maps $\theta$ and $\sigma$ as $\pi_{\theta}$ and $\pi_{\sigma}$, respectively.  By the functional calculus, $f_{\beta_1, \beta_2}(\sigma)$ is a well-defined order zero map for any choice of $0 < \beta_1  < \beta_2 < 1$.

We claim that $d_{\tau}(\theta(p)) <  {\tau}(\sigma(f_{\beta_1, \beta_2})(p))$ for every projection $p \in F$ and every $\tau \in T(A)$.

Note that if $\sigma$ is a $^{*}$-homomorphism then $f_{\beta_1, \beta_2}(\sigma) = \sigma$ for any choice of $0 < \beta_1 < \beta_2 < 1$. In this case, the claim follows immediately. 

Otherwise, we have 
\begin{equation*}
\textstyle
f_{\beta_1, \beta_2}(t) \geq t - \frac{\beta_1}{1- \beta_1} \cdot (1 - t) \text{ for all } t \in [0,1],
\end{equation*}
thus
\[
 f_{\beta_1, \beta_2} (\sigma)(p) \geq \sigma(p) - \textstyle{\frac{\beta_1}{1- \beta_1}} \cdot (1 - \sigma(p)).
  \]
Note that if $p \in F$ is nonzero then $\tau(\sigma(p)) \neq 0$ for any $\tau \in T(A \otimes \Q)$ since $A$ is simple. So $\tau(\sigma(p)) > 0$ for every $\tau \in T(A \otimes \Q)$. Also, since $\alpha < 1$ we also have $\beta_1 < \frac{1}{2}$. Thus 
\begin{eqnarray}
\nonumber \tau(f_{\beta_1, \beta_2}(\sigma)(p)) &\geq&
\nonumber \textstyle \tau(\sigma(p)) - \frac{\beta_1}{1 - \beta_1} \cdot d_{\tau} (1 - \sigma(p))\\
\nonumber &>& \tau(\sigma(p)) - 2 \beta_1 \cdot d_{\tau} (1- \sigma(p))\\
\nonumber &{>}& \textstyle \tau(\sigma(p)) - \alpha\\
\label{dt} &\geq& d_\tau(\theta(p))
\end{eqnarray}
for every $\tau \in T(A \otimes \Q)$ and every nonzero projection $p \in F$, proving the claim.

Write
$$F = M_{r_1} \oplus \cdots \oplus M_{r_L},$$
and for $l = 1, \dots, L$, let $e_{i,j}^{(l)}$ denote the partial isometry in $F$ corresponding to the $(i,j)^{\text{th}}$ matrix unit in $M_{r_l}$. For $1 \leq l \leq L$, by (\ref{dt}) we have that
\[ d_{\tau}(\theta(e_{1,1}^{(l)})) <\tau(f_{\beta_1, \beta_2}(\sigma)(e_{1,1}^{(l)}))\text{ for all } \tau \in T(A),\]
so by strict comparison it follows that $\theta(e_{1,1}^{(l)})\precsim f_{\beta_1, \beta_2}(\sigma)(e_{1,1}^{(l)})$. By \cite[Proposition 2.4]{Ror:uhfII}, there are unitaries $u_l \in A$ such that 
\begin{equation} \label{u_l} u_l (g_{\epsilon/2, \epsilon}(\theta)(e_{1,1}^{(l)}))  u_l^* \in \her( f_{\beta_1,\beta_2}(\sigma)(e_{1,1}^{(l)})).
\end{equation}
Let 
\begin{equation} \label{d_l}
d_l = (g_{\epsilon/2, \epsilon}(\theta)(e_{1,1}^{(l)}))^{1/2} u_l^*.
\end{equation}
Then $d_l$ satisfies
\begin{equation} \label{d_l 2}
d_l d_l^*(\theta(e_{1,1}^{(l)}) - \epsilon)_+ 
= g_{\epsilon/2, \epsilon}(\theta)(e_{1,1}^{(l)}))(\theta(e_{1,1}^{(l)}) - \epsilon)_+
=(\theta(e_{1,1}^{(l)}) - \epsilon)_+, 
\end{equation}
and similarly $(\theta(e_{1,1}^{(l)}) - \epsilon)_+ d_ld_l^* = (\theta(e_{1,1}^{(l)}) - \epsilon)_+$.

Furthermore, since $d_l^* (\theta(e_{1,1}^{(l)}) - \epsilon)_+ d_l \in \her( f_{\beta_1, \beta_2}(\sigma)(e_{1,1}^{(l)}))$ by (\ref{u_l}), we have
\begin{eqnarray}
\nonumber g_{0, \beta_1}(\sigma)(e_{1,1}^{(l)}) d_l^* (\theta(e_{1,1}^{(l)}) - \epsilon)_+ d_l  &\stackrel{\eqref{gf=f}}{=}& d_l^* (\theta(e_{1,1}^{(l)}) - \epsilon)_+ d_l \\
\label{g commutes} &\stackrel{\eqref{gf=f}}{=}& d_l^* (\theta(e_{1,1}^{(l)}) - \epsilon)_+ d_l g_{0, \beta_1}(\sigma)(e_{1,1}^{(l)}).
\end{eqnarray}

Set
\begin{equation} \label{s}
\textstyle s = \sum_{l=1}^L \sum_{k=1}^{r_l} \pi_{\theta}(e_{k,1}^{(l)}) d_l \pi_{\sigma}(e_{1,k}^{(l)}).
\end{equation}
Note that since $d_l^* d_l \in \her(f_{\beta_1, \beta_2}(\sigma)(e_{1,1}^{(l)}))$ we have that 
\[
d_l \pi_{\sigma}(e_{1,k}^{(l)}) \stackrel{\eqref{gf=f}}{=} d_l g_{0, \beta_1}(\sigma)(e_{1,1}^{(l)})\pi_{\sigma}(e_{1,k}^{(l)})  = d_l g_{0, \beta_1}(\sigma)(e_{1,k}^{(l)}) \in A,
\]
and similarly, since $d_l d_l^* \in \her(g_{\epsilon/2, \epsilon}(\theta)(e_{1,1}^{(l)}))$ we have
\[
\pi_{\theta}(e_{k,1}^{(l)}) d_l \stackrel{\eqref{gf=f}}{=} \pi_{\theta}(e_{k,1}^{(l))}) g_{\epsilon/4, \epsilon/2}(\theta)(e_{1,1}^{(l)}) d_l =  g_{\epsilon/4, \epsilon/2}(\theta)(e_{k,1}^{(l)}) d_l \in A
\]
thus
\[
\pi_{\theta}(e_{k,1}^{(l)}) d_l \pi_{\sigma}(e_{1,k}^{(l)}) = g_{\epsilon/4, \epsilon/2}(\theta)(e_{k,1}^{(l)}) d_l g_{0, \beta_1}(\sigma)(e_{1,k}^{(l)}) \in A,
\]
and hence $s \in A$. 

Since the hereditary $C^*$-subalgebras $\her(f_{\beta_1, \beta_2}(\sigma)(e_{1,1}^{(l)}))$ are pairwise orthogonal, we have that  $d_l d_{l'}^* = 0$ when $l \neq l'$ and 
\begin{equation} \label{ss*}
\textstyle ss^* = \sum_{l=1}^L \sum_{k=1}^{r_l} \pi_{\theta}(e_{k,1}^{(l)}) d_l d_l^* \pi_{\sigma}(e_{1,k}^{(l)}).
\end{equation}

We have that $s^*s \in \her(f_{\beta_1, \beta_2}(\sigma(\oplus_{l=1}^L e_{1,1}^{(l)})) \subset \her(f_{\beta_1, \beta_2}(\sigma(1_F))$, showing \ref{s*s location}.

For \ref{theta commutes}, it is obviously enough to show that $(\theta(e_{i,j}^{(l)}) - \epsilon)_+  ss^*= (\theta(e_{i,j}^{(l)}) - \epsilon)_+ = ss^*(\theta(e_{i,j}^{(l)}) - \epsilon)_+$ for arbitrary $i, j, l$.  Furthermore, since $\theta$ is order zero, it is clear that $(\theta(e_{i,j}^{(l)}) - \epsilon)_+   \pi_{\theta}(e_{k,1}^{(l')}) = 0$ when $l \neq l'$. Thus
\begin{eqnarray*}
\lefteqn{(\theta(e_{i,j}^{(l)}) - \epsilon)_+  ss^*} \\
& \stackrel{\eqref{ss*}}{=}& (\theta(e_{i,j}^{(l)}) - \epsilon)_+ (\textstyle{\sum}_{k=1}^{r_l} \pi_{\theta}(e_{k,1}^{(l)}) d_l d_l^* \pi_{\theta}(e_{1,k}^{(l)}))\\
&=& (\theta(e_{i,j}^{(l)}) - \epsilon)_+ \pi_{\theta}(e_{j,1}^{(l)}) d_l d_l^* \pi_{\theta}(e_{1,j}^{(l)})\\
&=& \pi_{\theta}(e_{i,1}^{(l)})  (\theta(e_{1,1}^{(l)}) - \epsilon)_+ d_l d_l^* \pi_{\theta}(e_{1,j}^{(l)})\\
& \stackrel{\eqref{d_l 2}}{=}& \pi_{\theta}(e_{i,1}^{(l)})  (\theta(e_{1,1}^{(l)}) - \epsilon)_+  \pi_{\theta}(e_{1,j}^{(l)})\\
&=& \pi_{\theta}(e_{i,1}^{(l)}) (\theta(e_{1,j}^{(l)}) - \epsilon)_+ \\
&=& (\theta(e_{i,j}^{(l)}) - \epsilon)_+ .
\end{eqnarray*}
 
The fact that $ss^* (\theta(e_{i,j}^{(l)}) - \epsilon)_+ = (\theta(e_{i,j}^{(l)}) - \epsilon)_+ $ follows from a nearly identical calculation.

For \ref{g(sigma)}, again it suffices to show the case $x = e_{i,j}^{(l)}$.

\begin{eqnarray*}
\lefteqn{s^* (\theta(e_{i,j}^{(l)}) - \epsilon)_+ s}\\
 &=& (\pi_{\sigma}(e_{i,1}^{(l)}) d_l^* \pi_{\theta}(e_{1,i}^{(l)}))  (\theta(e_{i,j}^{(l)}) - \epsilon)_+ (\pi_{\theta}(e_{j,1}^{(l)}) d_l  \pi_{\sigma}(e_{1,j}^{(l)}))\\
&=&  \pi_{\sigma}(e_{i,1}^{(l)}) d_l^*  (\theta(e_{1,1}^{(l)}) - \epsilon)_+ d_l  \pi_{\sigma}(e_{1,j}^{(l)})\\
&\stackrel{\eqref{g commutes}}{=}&  \pi_{\sigma}(e_{i,1}^{(l)}) g_{0, \beta_1}(\sigma)(e_{1,1}^{(l)})d_l^*  (\theta(e_{1,1}^{(l)}) - \epsilon)_+ d_l  \pi_{\sigma}(e_{1,j}^{(l)})\\
&=&   g_{0, \beta_1}(\sigma)(e_{i,j}^{(l)}) \pi_{\sigma}(e_{j,1}^{(l)})d_l^*  (\theta(e_{1,1}^{(l)}) - \epsilon)_+ d_l  \pi_{\sigma}(e_{1,j}^{(l)}) \\
&=&  g_{0, \beta_1}(\sigma)(e_{i,j}^{(l)}) \pi_{\sigma}(e_{j, 1}^{(l)})d_l^* \pi_{\theta}(e_{1,j}^{(l)}) (\theta(e_{j,j}^{(l)}) - \epsilon)_+ \pi_{\theta}(e_{j,1}^{(l)})  d_l  \pi_{\sigma}(e_{1,j}^{(l)}) \\
&=&  g_{0, \beta_1}(\sigma)(e_{i,j}^{(l)}) s^* (\theta(1_F) - \epsilon)_+ s.
\end{eqnarray*}

Similarly, $s^* (\theta(e_{i,j}^{(l)}) - \epsilon)_+ s = s^* (\theta(1_F) - \epsilon)_+ s g_{0, \beta_1}(\sigma)(e_{i,j}^{(l)})$.
\end{nproof}
\en

%%%%%%%%%%%%%%%%%%%%%%%%theorem
\bn
\begin{ntheorem}
\label{main result}
Let $A$ be a separable simple unital locally recursive subhomogeneous $\mathrm{C}^{*}$-algebra. Suppose that an approximating recursive subhomogeneous algebra $B$ can always be chosen to have an $(\mathcal{F},\eta)$-connected  recursive subhomogeneous decomposition 
\[
[B_{l},X_{l},\Omega_{l},r_{l},\phi_{l}]_{l=1}^{R}
\]
along which projections can be lifted and such that  $X_{l} \setminus \Omega_{l} \neq \emptyset$ for $l \ge 1$. Suppose further that $A$ has exactly $n$ extreme tracial states $\tau_0, \dots, \tau_{n-1} \in T(A)$ satisfying $(\tau_i)_* = (\tau_j)_*$ for every $i, j \in \{ 0, \dots, n-1\}$. Then, for any finite subset $\F \subset A_+^1$ and any $0 < \epsilon 1$, there are a partial isometry $s \in A \otimes \Q$, a finite dimensional $\mathrm{C}^{*}$-subalgebra $F \subset \Q$ and a $^{*}$-homomorphism $$\Phi : \mathcal{C}([0,1]) \otimes F  \to A \otimes \Q$$ such that
$$ss^* = \Phi(1_{\mathcal{C}([0,1])} \otimes 1_F)$$
and
\begin{enumerate}
\item $\| s^*s (a \otimes 1_{\Q}) - (a \otimes 1_{\Q}) s^*s \| < \epsilon$ for all $a \in \F$,
\item $\dist(s^*s (a \otimes 1_{\Q}) s^*s, s^*\Phi(\mathcal{C}([0,1]) \otimes F)s) < \epsilon$ for all $a \in \F$,
\item $\tau \otimes \tau_{\Q} (s^*s) \geq \frac{1}{2(n+2)}$ for all $\tau \in T(A)$.
\end{enumerate} 
\end{ntheorem}

\begin{nproof}
Let $\F$ and $\epsilon$ be given. Since $A$ is locally recursive subhomogeneous, we may assume, by taking a sufficiently good approximation, that $\F \subset B$ for some recursive subhomogeneous $\mathrm{C}^{*}$-algebra $B$.

We may furthermore assume that $1_A \subset \F$ so that $\tau^{(B)}_i  := \tau_i|_B \in T(B)$ are faithful and, as states on $K_0(B)$, satisfy $(\tau^{(B)}_i)_* = (\tau^{(B)}_j)_*$  for all $i, j \in \{ 0, \dots, n-1\}$.  

Let  $\eta > 0$ and $\beta_1 < \frac{1}{8}$ be so small that
\begin{equation} \label{eta_0}
\eta < \textstyle \frac{\beta_1}{6} \cdot \epsilon.
\end{equation}

We may apply Lemma~\ref{E} with respect to $\eta$ and $\F$ to get 
\[
0=K_{0} < K_{1 }< \ldots < K_{n-1} = K \in \mathbb{N}
\]
and pairwise orthogonal weighted $(\mathcal{F},\eta)$-excisors 
\[
(Q_{\frac{j}{K}},\rho_{\frac{j}{K}},\tilde{\sigma}_{\frac{j}{K}},\kappa_{\frac{j}{K}}), \; j \in \{0,\ldots,K\},
\]
implementing $(\mathcal{F},\eta)$-bridges 
\begin{eqnarray*}
(Q_{\frac{K_{0}}{K}},\rho_{\frac{K_{0}}{K}},\tilde{\sigma}_{\frac{K_{0}}{K}},\kappa_{\frac{K_{0}}{K}}) &  \sim_{(\mathcal{F},\eta)}  \ldots \sim_{(\mathcal{F},\eta)} & (Q_{\frac{K_{m}}{K}},\rho_{\frac{K_{m}}{K}},\tilde{\sigma}_{\frac{K_{m}}{K}},\kappa_{\frac{K_{m}}{K}}) \\ 
& \sim_{(\mathcal{F},\eta)} \ldots \sim_{(\mathcal{F},\eta)} &  (Q_{\frac{K_{n-1}}{K}},\rho_{\frac{K_{n-1}}{K}},\tilde{\sigma}_{\frac{K_{n-1}}{K}},\kappa_{\frac{K_{n-1}}{K}}) ,
\end{eqnarray*}
and such that, for each projection $q \in Q_{\frac{K_{i}}{K}}$, $i \in \{0,\ldots,n-1\}$, from \eqref{E9} we have
\begin{equation}
 \label{trace compare}
\tau_i \otimes \tau_{\mathcal{Q}} (\tilde{\sigma}_{\frac{K_{i}}{K}}(q)) \ge \textstyle \frac{1}{n+1} \cdot \tau_{\mathcal{Q}} \kappa_{\frac{K_{i}}{K}}(q).
\end{equation}

Let $0 < \alpha_1< \alpha_2 < \frac{1}{2(n-1)}$ and choose
\begin{equation}
\label{delta}
0 < \delta < \textstyle\frac{2}{3}
\end{equation}
to apply Lemma~\ref{interval map} with
\begin{equation} \label{delta_0}
\textstyle 0<  \delta_0 < \frac{(n-1)\delta \alpha_2}{2n}
\end{equation}
 to get a $^{*}$-homomorphism 
\[
\phi : \mathcal{C}([0,1]) \to A \otimes \Q
\]
satisfying 
\[
\tau_i \otimes \tau_{\Q}(\phi(\tilde{\gamma}_i)) \geq 1 - \delta_0,
\]
and
\[
0 < \tau_j \otimes \tau_{\Q}(\phi(\tilde{\gamma}_i)) < \delta_0
\]
for $i \neq j \in \{0, \dots, n-1\}$, where 
 \[
 \textstyle \tilde{\gamma}_i(t) = \left\{ \begin{array}{ccl} 0, & \, & t \in [0,1]\cap ((-\infty, \frac{i-1}{n-1}] \cup [\frac{i+1}{n-1}, \infty))\\  1, & \, & t = \frac{i}{n-1} \\ \text{ linear, } && \text{elsewhere.} \end{array} \right. 
 \]

For $i \in \{0, \dots, n-1\}$, define $\hat{\gamma}_i \in \mathcal{C}([0,1])$ by 
$$\hat{\gamma}_i= g_{\frac{i}{n-1} - \alpha_2, \frac{i}{n-1} - \alpha_1} - g_{\frac{i}{n-1} + \alpha_1, \frac{i}{n-1} + \alpha_2}$$
and for $i \in \{0, \dots, n-2\}$ define $\gamma_{i, i+1} \in \mathcal{C}([0,1])$ by
$$\gamma_{i,i+1} = g_{\frac{i}{n-1}+ \alpha_1, \frac{i}{n-1} + \alpha_2} - g_{\frac{i+1}{n-1} - \alpha_2, \frac{i+1}{n-1} - \alpha_1},$$
where we set $g_{ - \alpha_2, - \alpha_1} = 1$ and $g_{1+ \alpha_1, 1+\alpha_2} = 0$.
Note that 
\[ \textstyle \hat{\gamma}_{n-1} +  \sum_{i = 0}^{n-2} \hat{\gamma}_i + \gamma_{i, i +1} = 1_{\mathcal{C}([0,1])}
\]

We will now estimate the traces of the $\phi(\hat{\gamma}_i),$ and $\phi(\gamma_{i, i+1})$. We have 
\[0 \leq \hat{\gamma}_{i-1}(t), \gamma_{i-1,i}(t) \leq  \frac{1}{(n-1) \alpha_2} \cdot \tilde{\gamma}_{i-1}(t)\]
for all $t \in [0,1]$, for all $i \in \{1,\dots, n-1\}$, so 
\begin{eqnarray*}
\tau_i \otimes \tau_{\Q}(\phi(\gamma_{i-1,i})), \tau_i \otimes \tau_{\Q}(\phi(\hat{\gamma}_{i-1})) &\leq&  \textstyle \frac{1}{(n-1) \alpha_2} \cdot \tau_i \otimes \tau_{\Q}(\phi(\tilde{\gamma}_{i-1})) \\
&<& \textstyle \frac{\delta_0}{(n-1) \alpha_2} \\
&\stackrel{\eqref{delta_0}}{<}& \textstyle \frac{\delta}{2n}.
\end{eqnarray*}
One similarly shows that
\[
\tau_i \otimes \tau_{\Q}(\phi(\gamma_{i, i+1})), \tau_i \otimes \tau_{\Q}(\phi(\hat{\gamma}_{i+1})) 
\leq \textstyle \frac{1}{(n-1) \alpha_2} \cdot \tau_i\otimes \tau_{\Q}(\phi(\tilde{\gamma}_{i+1})) 
<\textstyle \frac{\delta}{2n}.
\]
It follows that 
\begin{equation} \label{trace large}
\textstyle \tau_i\otimes \tau_{\Q}(\phi(\hat{\gamma}_i)) = \tau_i \otimes \tau_{\Q}(1 - \sum_{j=0}^{n-2} \phi(\gamma_{j,j+1}) - \sum_{j = 0, j \neq i}^{n-1} \phi(\hat{\gamma}_j)) > 1 - \delta.
\end{equation}

Let $t_0 < t_1 < \dots < t_K$ be a partition of the interval $[0,1]$ satisfying
$$\textstyle t_{K_{i-1}} = \frac{i-1}{n-1}, \quad  t_{K_{i-1}+1} = \frac{i-1}{n-1} + \alpha_1, \quad t_{K_{i-1} + 2} = \frac{i-1}{n-1}+ \alpha_2$$
and
$$\textstyle t_{K_i - 2} = \frac{i}{n-1} - \alpha_2, \quad t_{K_i -1} = \frac{i}{n-1} - \alpha_1, \quad t_{K_i} = \frac{i}{n-1}$$
for $i \in \{1, \dots, n-1\}$. 
When $j = K_i$ for some $i \in \{0, \dots, n-1\}$, set 
\[\gamma_{\frac{j}{K}} := \hat{\gamma}_j.\]
When $j \in \{0, \dots, K\} \setminus \{K_0, \dots, K_{n-1}\}$, define
 $$ \gamma_{\frac{j}{K}}(t)= \left\{
\begin{array}{cl}
0, & 0 \leq t \leq t_j  \text{ and } t \geq t_{j+2}\\
1, &  t = t_{j+1} \\
\text{linear,} & t_j \leq t \leq t_{j+1} \text{ and }  t_{j+1} \leq t \leq t_{j+2}, \\
\end{array} \right. $$
so that the $\gamma_{\frac{j}{K}}$ are a partition of unity corresponding to $t_0 < t_1 < \cdots < t_K$.

Let $p \in \Q$ be a projection satisfying $\tau_{\Q}(p) = \frac{1}{n+2}$. Then by  \eqref{trace compare} and the choice of $\delta$ we have, for each $0 \leq j \leq K$, that
$$\tau \otimes \tau_{\Q}(\phi(\gamma_{\frac{j}{K}}) \otimes \kappa_{\frac{j}{K}}(q) \otimes p) < \tau \otimes \tau_{\Q}(\tilde{\sigma}_{\frac{j}{K}}(q))$$
for all $\tau \in T(A \otimes \Q)$ and for all projections $q \in Q_{\frac{j}{K}}$.

Define c.p.c.\ order zero maps  
\[
\theta_{\frac{j}{K}} : Q_{\frac{j}{K}} \to A \otimes \Q \otimes \Q \otimes \Q \cong A \otimes \Q
\]
by
\begin{equation}
\label{theta}
 \theta_{\frac{j}{K}}(a) = \phi(\gamma_{\frac{j}{K}}) \otimes \kappa_{\frac{j}{K}}(a) \otimes p.
 \end{equation}

Each finite-dimensional $\mathrm{C}^{*}$-algebra $Q_{\frac{j}{K}}$, $j\in \{0, \dots, K\}$ can be written as a sum of $L^{(j)} \in \mathbb{N}$ matrix algebras,  $Q_{\frac{j}{K}} = M_{r^{(j)}_1} \oplus \cdots \oplus M_{r^{(j)}_{L^{(j)}}}$ for some  $r^{(j)}_1, \dots, r^{(j)}_{L^{(j)}} \in \mathbb{N}$.
  
Note that for every $n \in \mathbb{N}$ and every $a \in (Q_{\frac{j}{K}})_+$ we have that
\[
 \tau \otimes \tau_{\Q}(\theta_{\frac{j}{K}}(a)^{1/n}) \leq \tau_{\Q} \otimes \tau_{\Q}((\kappa_{\frac{j}{K}}(a) \otimes p)^{1/n}) \text{ for all } \tau \in T(A).
 \]
Thus we see that for every projection $q \in Q_{\frac{j}{K}}$
\begin{eqnarray}
\nonumber d_{\tau}(\theta_{\frac{j}{K}}(q)) &<&  d_{\tau_{\Q} \otimes \tau_{\Q}}(\kappa_{\frac{j}{K}}(q) \otimes p) \\
\nonumber&=& \tau_{\Q} \otimes \tau_{\Q}(\kappa_{\frac{j}{K}}(q) \otimes p) \\
\nonumber&=& \textstyle \frac{1}{n+2}  \tau_{\Q}(\kappa_{\frac{j}{K}}(q)) \\
\nonumber&<& \textstyle \frac{1}{n+1} \tau_{\Q}(\kappa_{\frac{j}{K}}(q)) \\
&\stackrel{\eqref{trace compare}}{\leq} &
\tau(\tilde{\sigma}_{\frac{j}{K}}(q)) \label{d_tau compare}
\end{eqnarray}
 for all $\tau \in T(A \otimes \Q).$
 
Define order zero maps by
$$\sigma_{\frac{j}{K}} = g_{0, \beta_1}(\tilde{\sigma}_{\frac{j}{K}}).$$
Note that 
\begin{equation} 
\label{g estimate}
\textstyle g_{0, \beta_1}(t) \leq \frac{1}{\beta_1} t, \text{ for all } t \in [0,1],
\end{equation}
thus
\begin{eqnarray}
\nonumber &&  \| \sigma_{\frac{j}{K}}(1_{Q_{\frac{j}{K}}}) (b \otimes 1_{\Q}) - \sigma_{\frac{j}{K}}(\rho_{\frac{j}{K}}(b)) \| \\
\nonumber &=& \| g_{0, \beta_1}(\tilde{\sigma}_{\frac{j}{K}}(1_{Q_{\frac{j}{K}}})) \pi_{\tilde{\sigma}_{\frac{j}{k}}}(1_{Q_{\frac{j}{K}}})( b \otimes 1_{\Q}) -  g_{0, \beta_1}(\tilde{\sigma}_{\frac{j}{K}}(1_{Q_{\frac{j}{K}}}))\pi_{\tilde{\sigma}_{\frac{j}{k}}}(\rho_{\frac{j}{K}}(b)) \| \\
\nonumber &\stackrel{\eqref{g estimate}}{\leq}& \textstyle \frac{1}{\beta_1} \| \tilde{\sigma}_{\frac{j}{K}}(1_{Q_{\frac{j}{K}}}) (b \otimes 1_{\Q}) - \tilde{\sigma}_{\frac{j}{K}}(\rho_{\frac{j}{K}}(b)) \| \\
&\stackrel{\eqref{eta_0}}{\leq}& \textstyle \frac{\epsilon}{6}. \label{sig - sig rho est}
\end{eqnarray}

Let $\eta_1 > 0$ be so small that if $a \in A$ and $p$ is a projection such that $\| a^*a - p\| < \eta_1$ then there is $v \in A$ such that $v^*v = p$ and $\|v - a \| < \frac{\epsilon}{12}.$ 

Since \eqref{d_tau compare} holds for $\theta_{\frac{j}{K}}$ and $\tilde{\sigma}_{\frac{j}{K}}$ for every $j \in \{ 0 , \dots, K \}$, we may apply Lemma~\ref{get s} with 
\begin{equation}
\label{eta_1}
\textstyle \eta_0 = \min \{ \frac{\epsilon}{96}, \frac{1}{4}, \frac{\eta_1}{3}\}
\end{equation}
in place of $\epsilon$ to each $j \in \{ 0, \dots, K\}$ to get elements
\[
s_{\frac{j}{K}} \in A \otimes \Q
\]
satisfying
\begin{equation}
s_{\frac{j}{K}}^* s_{\frac{j}{K}} \in \her(f_{\beta_1', \beta_2} (\tilde{\sigma}_{\frac{j}{K}})(1_{Q_{\frac{j}{K}}})) \label{orthogonal} 
\end{equation}
with $0< \beta_1' < \min\{(\tilde{\sigma}_{\frac{j}{K}}(q)) -d_{\tau}(\theta_{\frac{j}{K}}(q)), \beta_1\}$, and 
\begin{eqnarray}
\lefteqn{s_{\frac{j}{K}} s_{\frac{j}{K}}^* (\theta_{\frac{j}{K}}(a) - \eta_0)_+}  \label{s_js_j* a unit}  \\
\nonumber &=& (\theta_{\frac{j}{K}}(a) - \eta_0)_+s_{\frac{j}{K}} s_{\frac{j}{K}}^*\\
 \nonumber &=& (\theta_{\frac{j}{K}}(a)- \eta_0)_+ \text{ for all }a \in Q_{\frac{j}{K}},
 \end{eqnarray}
 and
 \begin{eqnarray}
 \lefteqn{\tilde{\sigma}_{\frac{j}{K}}(a)s_{\frac{j}{K}}^* (\theta_{\frac{j}{K}}(1_{Q_{\frac{j}{K}}})- \eta_0)_+ s_{\frac{j}{K}}} \label{swap}\\
 \nonumber  &=& s_{\frac{j}{K}}^* (\theta_{\frac{j}{K}}(1_{Q_{\frac{j}{K}}})- \eta_0)_+ s_{\frac{j}{K}}\tilde{\sigma}_{\frac{j}{K}}(a)\\
\nonumber   &=& s_{\frac{j}{K}}^* (\theta_{\frac{j}{K}}(a)- \eta_0)_+ s_{\frac{j}{K}} \text{ for all } a \in Q_{\frac{j}{K}}.
 \end{eqnarray}

If $f$ is a continuous function then upon approximating by polynomials we have
\begin{equation}
s_{\frac{j}{K}} s_{\frac{j}{K}}^*f((\theta(1_{Q_{\frac{j}{K}}}) - \eta_0)_+) = f((\theta(1_{Q_{\frac{j}{K}}}) - \eta_0)_+) s_{\frac{j}{K}} s_{\frac{j}{K}}^* = f((\theta(1_{Q_{\frac{j}{K}}}) - \eta_0)_+) \label{commutes with functions}.
\end{equation}
Moreover 
$$f(s_{\frac{j}{K}}^* \theta(1_{Q_{\frac{j}{K}}}) - \eta_0)_+s_{\frac{j}{K}}) = s_{\frac{j}{K}}^* f((\theta(1_{Q_{\frac{j}{K}}}) - \eta_0)_+) s_{\frac{j}{K}},$$
hence 
\begin{equation}
\label{sigma(a) commutes}
\sigma_{\frac{j}{K}}(a) s_{\frac{j}{K}}^* f((\theta(1_{Q_{\frac{j}{K}}}) - \eta_0)_+)s_{\frac{j}{K}} = s_{\frac{j}{K}}^* f((\theta(1_{Q_{\frac{j}{K}}}) - \eta_0)_+)s_{\frac{j}{K}}\sigma_{\frac{j}{K}}(a) 
\end{equation}
for all $a \in Q_{\frac{j}{K}}$.
Now put 
\begin{equation} \label{s def}
\textstyle \tilde{s} = \sum_{j=0}^K ((\theta_{\frac{j}{K}}(1_{Q_{\frac{j}{K}}}) - \eta_0)_+)^{1/2} s_{\frac{j}{K}}.
\end{equation}

Since $\Q$ is a UHF algebra, there is a finite-dimensional $\mathrm{C}^{*}$-algebra $F \subset \Q$ such that 
\begin{equation}
\label{f-d alg}
\| 1_F - 1_{\Q}\| < \eta_0
\end{equation}
and 
\[ \{ \kappa_{\frac{j}{K}} \circ \rho_{\frac{j}{K}}(a) \mid a \in \F, 0 \leq j \leq K\} \subset_{\eta_0} F.\]
Let $\iota : F \into \Q$ be the inclusion map and set
\begin{equation}
\label{Phi}
\Phi := \phi(\cdot) \otimes \iota(\cdot) \otimes p : \mathcal{C}([0,1]) \otimes F \to A \otimes \Q \otimes \Q \otimes \Q\, (\cong A \otimes \Q).\end{equation}

Since the $\tilde{\sigma}_{\frac{j}{K}}$ have orthogonal images, it follows from \eqref{orthogonal} above that
\begin{equation}
\label{ss*=0}
s_{\frac{j}{K}} s_{\frac{j'}{K}}^* = 0,
\end{equation}
and
\begin{equation}
\label{s x sigma}
s_{\frac{j}{K}} \sigma_{\frac{j'}{K}}(a) = 0 \text{ for all } a \in Q_{\frac{j'}{K}}
\end{equation}
whenever $j \neq j'$. 

It follows from \eqref{ss*=0} that
\begin{eqnarray*}
\tilde{s}\tilde{s}^* &=& \textstyle \sum_{j=0}^K((\theta_{\frac{j}{K}}(1_{Q_{\frac{j}{K}}}) - \eta_0)_+)^{1/2} s_{\frac{j}{K}} s_{\frac{j}{K}}^*((\theta_{\frac{j}{K}}(1_{Q_{\frac{j}{K}}}) - \eta_0)_+)^{1/2}\\
&\stackrel{\eqref{s_js_j* a unit}}{=}& \textstyle \sum_{j=0}^K (\theta_{\frac{j}{K}}(1_{Q_{\frac{j}{K}}}) - \eta_0)_+,
\end{eqnarray*}
and we estimate
\begin{eqnarray*}
\lefteqn{\| \tilde{s}\tilde{s}^* - \Phi(1_{\mathcal{C}([0,1])} \otimes 1_F)\|}\\
 &=&\textstyle \| \sum_{j=0}^K  (\theta_{\frac{j}{K}}(1_{Q_{\frac{j}{K}}}) - \eta_0)_+ - \Phi(1_{\mathcal{C}([0,1])} \otimes 1_F) \| \\
&\stackrel{\eqref{theta}}{\leq}&\textstyle \| \sum_{j=0}^K \phi(\gamma_{\frac{j}{K}}) \otimes \kappa_{\frac{j}{K}}(1_{Q_{\frac{j}{K}}}) \otimes p - \Phi(1_{\mathcal{C}([0,1])} \otimes 1_F) \| + 2 \cdot \eta_0 \\
&\stackrel{\eqref{Phi}}{=}&\textstyle \| \phi(1_{\mathcal{C}([0,1]}) \otimes 1_{\Q} \otimes p - \phi(1_{\mathcal{C}([0,1])}) \otimes \iota(1_F) \otimes p\| + 2 \cdot \eta_0 \\
&\stackrel{\eqref{f-d alg}}{=}& 3 \cdot \eta_0\\
&\stackrel{\eqref{eta_1}}{<}& \eta_1.
\end{eqnarray*}
By our choice of $\eta_1$ there is an honest partial isometry $s \in A \otimes \Q$ satisfying 
$$ss^* = \Phi(1_{\mathcal{C}([0,1])} \otimes 1_F)$$
and 
\begin{equation}
\label{s - tilde-s}
\| \tilde{s} - s \| < \textstyle \frac{\epsilon}{12}.
\end{equation}

Let $a \in \F$ and consider the element $\sum_{j=0}^K \sigma_{\frac{j}{K}}(\rho_{\frac{j}{K}}(a)) \in A \otimes Q$. We will use this to estimate $\| a \tilde{s}^*\tilde{s} - \tilde{s}^*\tilde{s} a \|$. 
Note that since the functions $\gamma_{\frac{j}{K}}$ and $\gamma_{\frac{j'}{K}}$ are pairwise orthogonal whenever $| j - j' | \geq2$ we have that $\theta_{\frac{j}{K}}(1_{Q_{\frac{j}{K}}})\theta_{\frac{j'}{K}}(1_{Q_{\frac{j'}{K}}}) = 0$ whenever $|j - j'| \geq 2$ whence
\begin{equation}
\label{ortho theta}
(\theta_{\frac{j}{K}}(1_{Q_{\frac{j}{K}}}) - \eta_0)_+^{1/2} s_{\frac{j}{K}})^* ((\theta_{\frac{j'}{K}}(1_{Q_{\frac{j'}{K}}}) - \eta_0)_+^{1/2}s_{\frac{j'}{K}}) = 0
\end{equation}
whenever $|j - j'| \geq 2.$
We calculate
\begin{eqnarray}
\nonumber  \lefteqn{\textstyle (\sum_{j=0}^K \sigma_{\frac{j}{K}}(\rho_{\frac{j}{K}}(a)))\tilde{s}^*\tilde{s}} \\
\nonumber &\stackrel{\eqref{s def}}{=}&\textstyle (\sum_{j=0}^K \sigma_{\frac{j}{K}}(\rho_{\frac{j}{K}}(a)))(\sum_{j=0}^K  s_{\frac{j}{K}}^*(\theta_{\frac{j}{K}}(1_{Q_{\frac{j}{K}}}) - \eta_0)_+^{1/2})  (\sum_{j=0}^K (\theta_{\frac{j}{K}}(1_{Q_{\frac{j}{K}}}) - \eta_0)_+^{1/2}s_{\frac{j}{K}})\\
\nonumber &\stackrel{\eqref{s x sigma}}{=}& \textstyle (\sum_{j=0}^K \sigma_{\frac{j}{K}}(\rho_{\frac{j}{K}}(a)) s_{\frac{j}{K}}^*(\theta_{\frac{j}{K}}(1_{Q_{\frac{j}{K}}}) - \eta_0)_+^{1/2})  \cdot (\sum_{j=0}^K (\theta_{\frac{j}{K}}(1_{Q_{\frac{j}{K}}}) - \eta_0)_+^{1/2}s_{\frac{j}{K}})\\
 &\stackrel{\eqref{ortho theta}}{=}&\textstyle \sum_{j=0}^{K} \sigma_{\frac{j}{K}}(\rho_{\frac{j}{K}}(a)) s_{\frac{j}{K}}^*(\theta_{\frac{j}{K}}(1_{Q_{\frac{j}{K}}}) - \eta_0)_+^{1/2} \label{sigma ss*}\\
 \nonumber &&\textstyle \cdot  (\sum_{\{j' \,\mid \, |j-j'| <2\}} \theta_{\frac{j'}{K}}(1_{Q_{\frac{j'}{K}}}) - \eta_0)_+^{1/2}s_{\frac{j'}{K}}).
\end{eqnarray}

A similar calculation yields
\begin{eqnarray}
 \lefteqn{\textstyle \tilde{s}^*\tilde{s}(\sum_{j=0}^K \sigma_{\frac{j}{K}}(\rho_{\frac{j}{K}}(a)))} \label{ss* sigma}\\
\nonumber &=&\textstyle \sum_{j=0}^{K} s_{\frac{j}{K}}^*(\theta_{\frac{j}{K}}(1_{Q_{\frac{j}{K}}}) - \eta_0)_+^{1/2}  \textstyle \cdot (\sum_{\{j' \,\mid \, |j-j'| <2\}}( \theta_{\frac{j'}{K}}(1_{Q_{\frac{j'}{K}}}) - \eta_0)_+^{1/2} s_{\frac{j'}{K}} \sigma_{\frac{j'}{K}} (\rho_{\frac{j'}{K}}(a))).
\end{eqnarray}
Thus
%%%%%%%%%%%%%%%%%%%%%%%%%%%
\begin{eqnarray*}
\lefteqn{\textstyle \| \tilde{s}^*\tilde{s} (\sum_{j=0}^K \sigma_{\frac{j}{K}}(\rho_{\frac{j}{K}}(a))) -  (\sum_{j=0}^K \sigma_{\frac{j}{K}}(\rho_{\frac{j}{K}}(a))) \tilde{s}^*\tilde{s}\|}\\
&\stackrel{\eqref{sigma ss*},\eqref{ss* sigma}}{\leq}&\textstyle \| \sum_{j=0}^K s_{\frac{j}{K}}^*(\theta_{\frac{j}{K}}(1_{Q_{\frac{j}{K}}}) - \eta_0)_+s_{\frac{j}{K}} \sigma_{\frac{j}{K}}(\rho_{\frac{j}{K}}(a)) \\
&&\textstyle -  \sigma_{\frac{j}{K}}(\rho_{\frac{j}{K}}(a))s_{\frac{j}{K}}^*(\theta_{\frac{j}{K}}(1_{Q_{\frac{j}{K}}}) - \eta_0)_+s_{\frac{j}{K}} \| \\
&&\textstyle + \| \sum_{j=0}^K \sum_{\{j' \, \mid \, |j-j'|=1\}} s_{\frac{j}{K}}^*(\theta_{\frac{j}{K}}(1_{Q_{\frac{j}{K}}}) - \eta_0)_+^{1/2}\\
&& (\theta_{\frac{j'}{K}}(1_{Q_{\frac{j'}{K}}}) - \eta_0)_+^{1/2}s_{\frac{j'}{K}}\sigma_{\frac{j'}{K}}(\rho_{\frac{j'}{K}}(a)) - \sigma_{\frac{j}{K}}(\rho_{\frac{j}{K}}(a) )s_{\frac{j}{K}}^*(\theta_{\frac{j}{K}}(1_{\Q_{\frac{j}{K}}}) - \eta_0)_+^{1/2}  \\
&& (\theta_{\frac{j'}{K}}(1_{\Q_{\frac{j'}{K}}}) - \eta_0)_+^{1/2} s_{\frac{j'}{K}}\|\\
&\stackrel{\eqref{swap}}{=}&\textstyle \| \sum_{j=0}^{K} \sum_{\{j' \, \mid \, |j-j'|=1\}} s_{\frac{j}{K}}^*(\theta_{\frac{j}{K}}(1_{Q_{\frac{j}{K}}}) - \eta_0)_+^{1/2}(\theta_{\frac{j'}{K}}(\rho_{\frac{j'}{K}}(a)) - \eta_0)_+^{1/2}s_{\frac{j'}{K}}\\
&& - s_{\frac{j}{K}}^*(\theta_{\frac{j}{K}}(\rho_{\frac{j}{K}}(a)) - \eta_0)_+^{1/2}(\theta_{\frac{j'}{K}}(1_{Q_{\frac{j'}{K}}}) - \eta_0)_+^{1/2}s_{\frac{j'}{K}} \| \\
&\leq&\textstyle \|  \sum_{i=0}^{D_1} s_{\frac{2i}{K}}^*(\theta_{\frac{2i}{K}}(1_{Q_{\frac{j}{K}}}) - \eta_0 )_+^{1/2}(\theta_{\frac{2i+1}{K}}(\rho_{\frac{2i+1}{K}}(a)) - \eta_0)_+^{1/2} s_{\frac{2i+1}{K}}\\
&& - s_{\frac{2i}{K}}^*(\theta_{\frac{2i}{K}}(\rho_{\frac{2i}{K}}(a)) - \eta_0 )_+^{1/2}(\theta_{\frac{2i+1}{K}}(1_{Q_{\frac{2i+1}{K}}}) - \eta_0)_+^{1/2} s_{\frac{2i+1}{K}}\\
&& + s_{\frac{2i+1}{K}}^*(\theta_{\frac{2i+1}{K}}(1_{Q_{\frac{j}{K}}}) - \eta_0 )_+^{1/2}(\theta_{\frac{2i}{K}}(\rho_{\frac{2i}{K}}(a)) - \eta_0)_+^{1/2} s_{\frac{2i}{K}}\\
&& - (\theta_{\frac{2i+1}{K}}(\rho_{\frac{2i+1}{K}}(a)) - \eta_0 )_+^{1/2}(\theta_{\frac{2i}{K}}(1_{Q_{\frac{2i}{K}}}) - \eta_0)_+^{1/2} s_{\frac{2i}{K}} \|\\
&&\textstyle + \| \sum_{i=0}^{D_2} s_{\frac{2i+1}{K}}^*(\theta_{\frac{2i+1}{K}}(1_{Q_{\frac{2i+1}{K}}}) - \eta_0 )_+^{1/2}(\theta_{\frac{2i+2}{K}}(\rho_{\frac{2i+2}{K}}(a)) - \eta_0)_+^{1/2} s_{\frac{2i+2}{K}} \\
&& - s_{\frac{2i+1}{K}}^* (\theta_{\frac{2i+1}{K}}(\rho_{\frac{2i}{K}}(a)) - \eta_0 )_+^{1/2}(\theta_{\frac{2i+2}{K}}(1_{Q_{\frac{2i+2}{K}}}) - \eta_0)_+^{1/2} s_{\frac{2i+2}{K}}\\
&& + s_{\frac{2i+2}{K}}^*(\theta_{\frac{2i+2}{K}}(1_{Q_{\frac{j}{K}}}) - \eta_0 )_+^{1/2}(\theta_{\frac{2i+1}{K}}(\rho_{\frac{2i+1}{K}}(a)) - \eta_0)_+^{1/2} s_{\frac{2i+1}{K}}\\
&& - s_{\frac{2i+2}{K}}^*(\theta_{\frac{2i+2}{K}}(\rho_{\frac{2i+2}{K}}(a)) - \eta_0 )_+^{1/2}(\theta_{\frac{2i+1}{K}}(1_{Q_{\frac{2i}{K}}}) - \eta_0)_+^{1/2} s_{\frac{2i+1}{K}} \|, \\
\end{eqnarray*}
where 
$$ D_1 = \left\{ \begin{array}{cl}
\textstyle \frac{D}{2} - 1, & \text{ if $D$ is even }\\
\textstyle \frac{D-1}{2}, & \text{ if $D$ is odd, }
\end{array} \right. \,
D_2 = \left\{ \begin{array}{cl}
\textstyle \frac{D}{2} - 1, & \text{ if $D$ is even }\\
\textstyle \frac{D-3}{2}, & \text{ if $D$ is odd.}
\end{array} \right. $$

Note that if $i$ and $i'$ are either both even or both odd, $i \neq i'$ we have
\begin{eqnarray*}
 \left(s_{\frac{i}{K}}^* x_0 s_{\frac{i+1}{K}} + s_{\frac{i+2}{K}}^* x_1 s_{\frac{i+1}{K}}\right) \cdot (s_{\frac{i'}{K}}^* x_2 s_{\frac{i'+1}{K}} + s_{\frac{i'+2}{K}}^* x_3 s_{\frac{i'+1}{K}})= 0,
\end{eqnarray*}
for any $x_0, \dots, x_3 \in A \otimes \Q$ since $| i +1 - i' | > 2$ implies
$$s_{\frac{i+1}{K}}(s_{\frac{i+1}{K}})^*s_{\frac{i'}{K}}(s_{\frac{i'}{K}})^* = 0.$$
Thus each sum in the norm estimates above consists of mutually orthogonal summands, allowing us to estimate

\begin{eqnarray}
\nonumber \lefteqn{\textstyle \| \tilde{s}^*\tilde{s} (\sum_{j=0}^K \sigma_{\frac{j}{K}}(\rho_{\frac{j}{K}}(a))) -  (\sum_{j=0}^K \sigma_{\frac{j}{K}}(\rho_{\frac{j}{K}}(a))) \tilde{s}^*\tilde{s}\|}\\
\nonumber &\leq& 2 \cdot \max_{j = 0 \dots, K} (\|(\theta_{\frac{j}{K}}(1_{Q_{\frac{j}{K}}}) - \eta_0 )_+^{1/2}(\theta_{\frac{j+1}{K}}((\rho_{\frac{j+1}{K}}(a)) - \eta_0)_+^{1/2}\\
\nonumber && - (\theta_{\frac{j}{K}}((\rho_{\frac{j}{K}}(a)) - \eta_0 )_+^{1/2})(\theta_{\frac{j+1}{K}}(1_{Q_{\frac{j+1}{K}}}) - \eta_0)_+^{1/2} \| \\
\nonumber && + \| (\theta_{\frac{j+1}{K}}(1_{Q_{\frac{j+1}{K}}}) - \eta_0 )_+^{1/2})(\theta_{\frac{j}{K}}(\rho_{\frac{j}{K}}(a)) - \eta_0)_+^{1/2}\\
\nonumber && - (\theta_{\frac{j+1}{K}}(\rho_{\frac{j+1}{K}}(a)) - \eta_0 )_+^{1/2})(\theta_{\frac{j}{K}}(1_{Q_{\frac{j}{K}}}) - \eta_0)_+^{1/2} \|) \\
\nonumber &\leq& 4 \cdot ( 4\eta_0^{1/2} + \max_{j=0, \dots, K}  \|  \theta_{\frac{j+1}{K}}(1_{Q_{\frac{j+1}{K}}})^{1/2}\theta_{\frac{j}{K}}(\rho_{\frac{j}{K}}(a))^{1/2} - \theta_{\frac{j+1}{K}}(\rho_{\frac{j+1}{K}}(a))^{1/2}\theta_{\frac{j}{K}}(1_{Q_{\frac{j}{K}}})^{1/2} \|) \\
\nonumber &\stackrel{\eqref{theta}}{<}& 4 \cdot (4 \eta_0^{1/2} + \max_{j=0, \dots, K}  \|  \kappa_{\frac{j+1}{K}}(1_{Q_{\frac{j+1}{K}}})^{1/2}\kappa_{\frac{j}{K}}(\rho_{\frac{j}{K}}(a))^{1/2} - \kappa_{\frac{j+1}{K}}(\rho_{\frac{j+1}{K}}(a))^{1/2}\kappa_{\frac{j}{K}}(1_{Q_{\frac{j}{K}}})^{1/2} \|) \\
\nonumber &\stackrel{\eqref{L1}}{<}& 16  \eta_0^{1/2} + 4 \eta^{1/2} \\
\nonumber &\stackrel{\eqref{eta_1}}{<}& \epsilon/6 + \epsilon/12\\
&\stackrel{\eqref{eta_0}}{<}& \epsilon/4. \label{s*s approx commutes}
\end{eqnarray}

It is straightforward to check that  $s_{\frac{j}{K}} =  s_{\frac{j}{K}} \sigma_{\frac{j}{K}}(1_{Q_{\frac{j}{K}}})$. Then note that
\begin{eqnarray}
\nonumber \| s_{\frac{j}{K}} \sigma_{\frac{j}{K}}(\rho_{\frac{j}{K}}(a))) - s_{\frac{j}{K}}( a \otimes 1_{\Q}) \| &\leq& \| \sigma_{\frac{j}{K}}(\rho_{\frac{j}{K}}(a)) - \sigma_{\frac{j}{K}}(1_{Q_{\frac{j}{K}}})(a \otimes 1_{\Q}) \| \\
&\stackrel{\eqref{sig - sig rho est}}{<}& \epsilon/6. \label{s a approx}
\end{eqnarray}
Thus
\begin{eqnarray*}
\lefteqn{\|s^*s (a \otimes 1_{\Q}) - (a \otimes 1_{\Q}) s^*s\|} \\
&\leq& 4 \cdot \| \tilde{s}- s\| + \|\tilde{s}^*\tilde{s} (a \otimes 1_{\Q}) - (a \otimes 1_{\Q}) \tilde{s}^*\tilde{s}\|\\
&\stackrel{\eqref{s - tilde-s}}{\leq}& \epsilon/3 + 2 \cdot  \max_j   \| s_{\frac{j}{K}} \sigma_{\frac{j}{K}}(\rho_{\frac{j}{K}}(a))) - s_{\frac{j}{K}}( a \otimes 1_{\Q}) \| \\
&& +  \|\textstyle \tilde{s}^*\tilde{s} (\sum_{j=0}^K \sigma_{\frac{j}{K}}(\rho_{\frac{j}{K}}(a))) -  (\sum_{j=0}^K \sigma_{\frac{j}{K}}(\rho_{\frac{j}{K}}(a))) \tilde{s}^*\tilde{s}\|\\
&\stackrel{\eqref{s a approx}, \eqref{s*s approx commutes}}{<}& \epsilon/3 + \epsilon/3 + \epsilon/4 \\
&<& \epsilon.
\end{eqnarray*}

For (ii), we calculate, for $a \in \F$,

\begin{eqnarray*}
\lefteqn{\textstyle \tilde{s} (\sum_{j=0}^K \sigma_{\frac{j}{K}}(\rho_{\frac{j}{K}}(a)) ) \tilde{s}^*} \\
%&=&\textstyle (\sum_{j=0}^K (\theta(1_{Q_{\frac{j}{K}}}) - \eta_0)_+^{1/2} s_{\frac{j}{K}} \sigma_{\frac{j}{K}}(\rho_{\frac{j}{K}}(a)) ) s^* \\
&\stackrel{\eqref{s x sigma}}{=}&\textstyle \sum_{j=0}^K (\theta(1_{Q_{\frac{j}{K}}}) - \eta_0)_+^{1/2} s_{\frac{j}{K}} \sigma_{\frac{j}{K}}(\rho_{\frac{j}{K}}(a)) s_{\frac{j}{K}}^* (\theta(1_{Q_{\frac{j}{K}}}) - \eta_0)_+^{1/2} \\
&\stackrel{\eqref{commutes with functions}}{=}&\textstyle \sum_{j=0}^K  s_{\frac{j}{K}} s_{\frac{j}{K}}^* (\theta(1_{Q_{\frac{j}{K}}}) - \eta_0)_+^{1/2} s_{\frac{j}{K}} \sigma_{\frac{j}{K}}(\rho_{\frac{j}{K}}(a)) s_{\frac{j}{K}}^* (\theta(1_{Q_{\frac{j}{K}}}) - \eta_0)_+^{1/2} \\
&\stackrel{\eqref{sigma(a) commutes}}{=}&\textstyle \sum_{j=0}^K  s_{\frac{j}{K}} \sigma_{\frac{j}{K}}(\rho_{\frac{j}{K}}(a))  s_{\frac{j}{K}}^* (\theta(1_{Q_{\frac{j}{K}}}) - \eta_0)_+^{1/2} s_{\frac{j}{K}} s_{\frac{j}{K}}^* (\theta(1_{Q_{\frac{j}{K}}}) - \eta_0)_+^{1/2} \\
&\stackrel{\eqref{commutes with functions}}{=}&\textstyle \sum_{j=0}^K  s_{\frac{j}{K}} \sigma_{\frac{j}{K}}(\rho_{\frac{j}{K}}(a))  s_{\frac{j}{K}}^* (\theta(1_{Q_{\frac{j}{K}}}) - \eta_0)_+s_{\frac{j}{K}} s_{\frac{j}{K}}^*\\
&\stackrel{\eqref{swap}}{=}&\textstyle \sum_{j=0}^K  s_{\frac{j}{K}} s_{\frac{j}{K}}^* (\theta(\rho_{\frac{j}{K}}(a)) - \eta_0)_+s_{\frac{j}{K}} s_{\frac{j}{K}}^*\\ 
&\stackrel{\eqref{s_js_j* a unit}}{=}&\textstyle \sum_{j=0}^K (\theta(\rho_{\frac{j}{K}}(a)) - \eta_0)_+ \\  
&\stackrel{\eqref{theta}}{=}&\textstyle \sum_{j=0}^K (\phi(\gamma_{\frac{j}{K}}) \otimes \kappa_{\frac{j}{K}}(\rho_{\frac{j}{K}}(a)) \otimes p - \eta_0)_+. \label{shs estimate}
\end{eqnarray*}

%$\mathfrak{insane}$

Define $h \in \mathcal{C}([0,1]) \otimes \Q \otimes \Q$ by
$$\textstyle h := \sum_{j=0}^K \phi(\gamma_{\frac{j}{K}}) \otimes \kappa_{\frac{j}{K}}(\rho_{\frac{j}{K}}(a)) \otimes p.$$

Let $a_{\frac{j}{K}} \in F$ be elements satisfying 
\begin{equation}
\label{approx F}
\| a_{\frac{j}{K}} - \kappa_{\frac{j}{K}}(\rho_{\frac{j}{K}}(a)) \| < \eta_0,
\end{equation}
and put 
\[
h' := \textstyle \sum_{j=0}^K \Phi(\gamma_{\frac{j}{K}} \otimes a_{\frac{j}{K}}) \in \Phi(\mathcal{C}([0,1]) \otimes F).
\]
Then

\begin{eqnarray*}
\| h - h' \| &\leq& \| \textstyle \sum_{j \text{ even }} \phi(\gamma_{\frac{j}{K}}) \otimes (\kappa_{\frac{j}{K}}(\rho_{\frac{j}{K}}(a)) - a_{\frac{j}{K}}) \|\\
&& \textstyle + \| \sum_{j \text{ odd }} \phi(\gamma_{\frac{j}{K}}) \otimes (\kappa_{\frac{j}{K}}(\rho_{\frac{j}{K}}(a)) - a_{\frac{j}{K}}) \| \\
&\stackrel{\eqref{approx F}}{<}& 2\cdot \eta_0\\
&\stackrel{\eqref{eta_1}}{<}& \epsilon/48 \\
&<& \epsilon/4.
\end{eqnarray*}

We calculate
\begin{eqnarray*}
\lefteqn{\textstyle \| \tilde{s}^*h \tilde{s} - \tilde{s}^*\tilde{s} (\sum_{j=0}^K \sigma_{\frac{j}{K}}(\rho_{\frac{j}{K}}(a)) )\tilde{s}^*\tilde{s}\|}\\
&\leq&\textstyle \| h - \tilde{s} (\sum_{j=0}^K \sigma_{\frac{j}{K}}(\rho_{\frac{j}{K}}(a)) )\tilde{s}^*\|\\
&\stackrel{\eqref{shs estimate}}{<}& 2 \cdot \eta_0\\
&\stackrel{\eqref{eta_1}}{<} & \epsilon/4,
\end{eqnarray*}
so
\begin{eqnarray*}
\|s^* h' s - s^*s(a \otimes 1_{\Q}) s^* s \| &\leq& \textstyle 6 \cdot \|s - \tilde{s}\| + \|h - h'\| \\
&& +  \| \tilde{s}^*h \tilde{s} - \tilde{s}^*\tilde{s} (\textstyle{\sum_{j=0}^K \sigma_{\frac{j}{K}}(\rho_{\frac{j}{K}}(a)) )\tilde{s}^*\tilde{s}}\|\\
&<& \epsilon/2 + \epsilon/4 +  \epsilon/4\\
&=& \epsilon.
\end{eqnarray*}
This shows (ii).

Finally,
\begin{eqnarray*}
\tau(s^* s) &=& \tau(ss^*)\\
&=& \tau(\phi(1_{\mathcal{C}([0,1])}) \otimes 1_{F} \otimes p)\\
&\stackrel{\eqref{f-d alg}}{>}& \textstyle \frac{1 - \eta_0}{n+2} \cdot \tau(\phi(1_{\mathcal{C}([0,1])}) \\
&\geq& \textstyle \frac{1 - \eta_0}{n+2}  \cdot (\sum_{i=0}^{n-1} \tau(\phi(\gamma_i))) \\
&\stackrel{\eqref{trace large}}{\geq}& \textstyle \frac{1 - \eta_0}{n+2}  \cdot  (1 - \delta)\\
&\stackrel{\eqref{delta}, \eqref{eta_1}}{\geq}& \textstyle \frac{1}{n+2} \cdot \frac{3}{4} \cdot \frac{2}{3}\\
&>& \textstyle \frac{1}{2(n+2)},
\end{eqnarray*}
for all $\tau \in T(A \otimes \Q)$, showing that (iii) holds.
\end{nproof}
\en

%%%%%%%%%%%%%%%%%%%%%%%%%%%%%%%%%%%%%%%%%%%%%%%%
\section{Main result, applications and outlook} \label{outlook} 

\bn
\begin{ntheorem} 
\label{TAI}
Let $A$ be a separable simple unital locally recursive subhomogeneous $\mathrm{C}^{*}$-algebra with exactly $n >0$ extreme tracial states $\tau_0, \dots, \tau_{n-1} \in T(A)$ satisfying $(\tau_i)_* = (\tau_j)_*$ for all $i, j \in \{ 0, \dots, n-1\}$. Suppose that an approximating recursive subhomogeneous algebra $B$ can always be chosen to have an $(\mathcal{F},\eta)$-connected  recursive subhomogeneous decomposition 
\[
[B_{l},X_{l},\Omega_{l},r_{l},\phi_{l}]_{l=1}^{R}
\]
along which projections can be lifted and such that  $X_{l} \setminus \Omega_{l} \neq \emptyset$ for $l \ge 1$. 
Then $A \otimes \Q$ is TAI.
\end{ntheorem}

\begin{nproof} 
The class I contains the finite dimensional $\mathrm{C}^{*}$-algebras, is closed under direct sums and tensor products with finite dimensional $\mathrm{C}^{*}$-algebras, and every $\mathrm{C}^{*}$-algebra in I can be written as a universal $\mathrm{C}^*$-algebra with weakly stable relations. Thus we may apply Lemma~\ref{away from 0}, and it is enough to show that there is an $m \in \mathbb{N}$ such that, for any $\epsilon >0$ and any finite subset $\mathcal{F} \subset A\otimes \Q$, there exist a projection $p \in A \otimes \Q$ and a unital $\mathrm{C}^{*}$-subalgebra $B \subset p(A \otimes \uhf)p$ and $B \in \mathcal{S}$ such that:
\begin{enumerate}
\item $\|pb - bp \| < \epsilon$ for all $b \in \mathcal{F}$,
\item $\dist(p b p, B) < \epsilon$ for all $b \in \mathcal{F}$,
\item $\tau(p) > 1/m$ for all $\tau \in T(A \otimes \Q)$.
\end{enumerate}
By Lemma~\ref{special F} we need only consider finite subsets of the form $\F = \G \otimes \{1_{\Q}\}$ for $\G \subset A$. Now the result follows from Theorem~\ref{main result} with $m = 2(n+2)$. \end{nproof}
\en

\bn
\begin{ncor}
\label{dim1 TAI} Let $A$ be a separable simple unital locally recursive subhomogeneous $\mathrm{C}^{*}$-algebra with exactly $n >0$ extreme tracial states $\tau_0, \dots, \tau_{n-1} \in T(A)$ satisfying $(\tau_i)_* = (\tau_j)_*$ for all $i, j \in \{ 0, \dots, n-1\}$. Suppose that an approximating recursive subhomogeneous algebra $B$ can always be chosen to have a recursive subhomogeneous decomposition 
\[
[B_{l},X_{l},\Omega_{l},r_{l},\phi_{l}]_{l=1}^{R}
\]
such that $\mathrm{dim} X_{l} \le 1$ for $l \ge 2$. Then $A \otimes \Q$ is TAI.
\end{ncor}

\begin{nproof}
Follows from Theorem~\ref{TAI} with Corollary~\ref{dim1}.
\end{nproof}
\en

Let $A$ be a unital separable simple $\mathrm{C}^{*}$-algebra. The Elliott invariant of $A$ is given by 
$$\Ell(A) = (K_0(A), K_0(A)_+, [1_A], K_1(A), T(A), r_A),$$
where $(K_0(A), K_0(A)_+, [1_A])$ is the partially ordered $K_0$-group with positive cone $K_0(A)_+$ and order unit $[1_A]$, $K_1(A)$ the $K_1$-group of $A$, $T(A)$ the simplex of tracial states and $r_A : T(A) \to S(K_0(A))$ is the map given by $r_A(\tau)([p] - [q]) = \tau(p) - \tau(q)$. For two $\mathrm{C}^{*}$-algebras $A$ and $B$, we write $\Ell(A) \cong \Ell(B)$ if there are an order and unit preserving group homomorphism $\phi : (K_0(A), K_0(A)_+, [1_A]) \to (K_0(B), K_0(B)_+, [1_B])$, a group homomorphism $\psi: K_1(A) \to K_1(B)$ and a homeomorphism $\gamma: T(B) \to T(A)$ such that the following diagram commutes:

 \begin{displaymath}
\xymatrix{ T(B)  \ar[d]^{r_{B}} \ \ar[r]^{\gamma}  & T(A)\ar[d]_{r_A}  \\
S(K_0(B))  \ar[r]^{\cdot \circ \phi} & S(K_0(A)).   \\}
\end{displaymath}

\bn
\begin{nnotation} We let $\mathcal{A}$ denote the class of $\mathrm{C}^{*}$-algebras such that if $A \in \mathcal{A}$ then $A$ is a unital separable simple locally recursive subhomogeneous $\mathrm{C}^{*}$-algebra such that the approximating recursive subhomogeneous algebra $B$ can always be chosen to have an $(\mathcal{F},\eta)$-connected  recursive subhomogeneous decomposition 
\[
[B_{l},X_{l},\Omega_{l},r_{l},\phi_{l}]_{l=1}^{R}
\]
along which projections can be lifted and such that  $X_{l} \setminus \Omega_{l} \neq \emptyset$ for $l \ge 1$. 
\end{nnotation}
\en

\bn
\begin{ncor} 
\label{TAI any UHF} Let $A \in \mathcal{A}$ with exactly $n >0$ extreme tracial states $\tau_0, \dots, \tau_{n-1} \in T(A)$ satisfying $(\tau_i)_* = (\tau_j)_*$ for all $i, j \in \{ 0, \dots, n-1\}$. Let $\mathfrak{p}$ be a supernatural number and $M_{\mathfrak{p}}$ the associated UHF algebra. Then $A \otimes M_{\mathfrak{p}}$ is TAI.
 \end{ncor}

 \begin{nproof} This follows immediately from Theorem~\ref{TAI} and \cite[Theorem 3.6]{LinNiu:RationallyTAI} with \cite[Theorem 7.1 (b)]{Lin:ttr}, which shows that a simple unital $\mathrm{C}^{*}$-algebra is TAI if and only if it has tracial rank less than or equal to one. We note that \cite[Theorem 7.1 (b)]{Lin:ttr} uses Gong's decomposition theorem \cite{Gong:simple_AH}. To avoid this technical theorem, we observe that tracial rank less than or equal to one can be replaced by TAI in the statements of Lemma 3.4 and Theorem 3.6 of \cite{LinNiu:RationallyTAI} and that the proofs work in exactly the same way by simply replacing all $\mathrm{C}^{*}$-algebras of tracial rank less than or equal to one with $\mathrm{C}^{*}$-algebras that are TAI, and invoking the more general Lemma 2.3 of \cite{EllNiu:tracial_approx} instead of \cite[Proposition 3.6]{Lin:simpleNuclearTR1}.
 \end{nproof}
\en

\bn
\begin{nprop} Let $A$ be a separable simple unital locally recursive subhomogeneous $\mathrm{C}^{*}$-algebra. Then $A$ satisfies the UCT.
\end{nprop}

\begin{nproof} For any $\epsilon >0$ and any finite subset $\F \subset A$ we may approximate $A$ by a subhomogeneous $\mathrm{C}^{*}$-algebra $B$. Since $B$ is Type I, $B$ satisfies the UCT. Therefore the result follows immediately by appealing to Theorem 1.1 of \cite{Dad:UCT}.
\end{nproof}
\en

 \bn
\begin{ncor} \label{classification}
Let $A, B \in \mathcal{A}$ be $\mathrm{C}^{*}$-algebras, and let $n \in \mathbb{N} \setminus\{0\}$. Suppose there are exactly $n$ extreme tracial states $\tau_0, \dots, \tau_{n-1} \in T(A)$ satisfying $(\tau_i)_* = (\tau_j)_*$ for all $i, j \in \{ 0, \dots, n-1\}$. Then 
$$A \otimes \mathcal{Z} \cong B \otimes \mathcal{Z} \text{ if and only if } \Ell(A \otimes \mathcal{Z}) \cong \Ell(B \otimes \mathcal{Z}).$$ 
If, in addition, $A$ and $B$ have finite decomposition rank, then 
$$A \cong B \text{ if and only if } \Ell(A) \cong \Ell(B).$$ 
\end{ncor}

\begin{nproof}
$A \otimes \Q$ and $B \otimes \Q$ are TAI by Theorem~\ref{TAI}. Since $A$ and $B$ satisfy the UCT, the result follows by applying  \cite[Corollary 11.9]{Lin:asu-class}. Since $A$ and $B$ are separable, simple, nonelementary and unital, the second statement then follows from the fact that finite decomposition rank implies $\mathcal{Z}$-stability \cite[Theorem 5.1]{Winter:dr-Z-stable}.
\end{nproof}
\en

In \cite[Section 5]{Ell:invariant}, Elliott constructs examples of approximately subhomogeneous $\mathrm{C}^{*}$-algebras by attaching one-dimensional spaces to the circle. These examples exhaust the Elliott invariant in the weakly unperforated case. In that paper, the Elliott invariant of these algebras is computed but classification results are not given. In the case of finitely many traces inducing the same state on $K_0$, we are able to obtain classification by the results above. In particular this shows that Elliott's examples, assuming the restriction to finitely many tracial states inducing the same state on $K_0$, agree with the examples of \cite{LinNiu:RationallyTAI}; this was previously unknown.

\bn
\begin{ncor} \label{Ell blocks}
Let $A$ and $B$ be inductive limits of building block algebras defined in \cite[Section 5.1.2]{Ell:invariant} and let $n \in \mathbb{N} \setminus\{0\}$. Suppose there are exactly $n$ extreme tracial states $\tau_0, \dots, \tau_{n-1} \in T(A)$ satisfying $(\tau_i)_* = (\tau_j)_*$ for all $i, j \in \{ 0, \dots, n-1\}$ and  exactly $n$ extreme tracial states $\tau'_0, \dots, \tau'_{n-1} \in T(B)$ satisfying $(\tau'_i)_* = (\tau'_j)_*$ for all $i, j \in \{ 0, \dots, m-1\}$. Then $A \otimes \Q$ and $B \otimes \Q$ are TAI and we have
$$A  \cong B \text{ if and only if } \Ell(A) \cong \Ell(B).$$ 
\end{ncor}

\begin{nproof}
By definition,  $A$ and $B$ can be written as inductive limits $A= \Lim A_n$ and $B = \Lim B_n$ where $A_n$ and $B_n$ are recursive subhomogeneous $\mathrm{C}^{*}$-algebras of topological dimension less than or equal to one. It follows from Corollary~\ref{dim1} that $A, B \in \mathcal{A}$ and thus by the assumptions on the tracial state spaces, we have $A \otimes \Q$ and $B \otimes \Q$ are TAI by Corollary~\ref{dim1 TAI}. Since the approximating algebras $A_n$ and $B_n$ all have dimension less than or equal to one, both $A$ and $B$ have finite decomposition rank. Thus classification follows from Corollary~\ref{classification}.
\end{nproof}
\en

At least to some extent, we are also able to apply our result in the context of $\mathrm{C}^{*}$-algebras of minimal dynamical systems.  In \cite{LinMat:CantorT1}, Lin and Matui study minimal dynamical systems on the product of the Cantor set and $\mathbb{T}$.  Let $X$ be the Cantor set and let $\xi: X \to \mathbb{T}$ be a continuous map. Then we can define $R_{\xi} : X \to \homeo(\mathbb{T})$ by $R_{\xi}(x)(t) = t + \xi(x)$ for $x \in X$ and $t \in \mathbb{T}$. If $\alpha:  X \to X$ is a homeomorphism of the Cantor set $X$, then 
$$ \alpha \times R_{\xi} : X \times \mathbb{T} \to X \times \mathbb{T} : (x, t) \mapsto (\alpha(x), R_{\xi}(x)(t))$$
is a homeomorphism of $X \times \mathbb{T}$.

In the case that the homeomorphisms $\alpha \times R_{\xi}$ are minimal, Lin and Matui show that the crossed products $\mathcal{C}(X \times \mathbb{T}) \rtimes_{\alpha \times R_{\xi}} \mathbb{Z}$ are tracially approximately finite or have tracial rank one and hence classifiable as they satisfy the UCT \cite[Theorem 4.3]{LinMat:CantorT1}. Under the additional assumption of finitely many extreme tracial states, all of which induce the same state at the level of $K_0$, our Theorem~\ref{TAI} offers an alternative route to the same result.

\bn
\begin{ncor}
Let $(X, \alpha)$ and $(Y, \beta)$ be Cantor dynamical systems, $\xi : X \to \mathbb{T}$ and $\zeta: Y \to \mathbb{T}$ continuous maps and suppose that $\alpha \times R_{\xi}$ and $\beta \times R_{\zeta}$ are minimal. Put $A := \mathcal{C}(X \times \mathbb{T}) \rtimes_{\alpha \times R_{\xi}}  \mathbb{Z}$ and $B :=  \mathcal{C}(Y \times \mathbb{T}) \rtimes_{\alpha \times R_{\zeta}} \mathbb{Z}$. Suppose $T(A)$ and $T(B)$ each have finitely many extreme points such that $[\tau_A]_* = [\tau_A']_*$ in $K_0(A)$ for every extreme point $\tau_A, \tau_A'$ and  $[\tau_B]_* = [\tau_B']_*$ in $K_0(B)$ for every extreme point $\tau_B, \tau_B'$. Then 
$$A \cong B \text{ if and only if } \Ell(A) \cong \Ell(B).$$ 
\end{ncor}

\begin{nproof}
Let $A$ and $B$ be as above. Let $u$ and $v$ be the canonical unitaries inducing the actions of $\mathbb{Z}$ in $A$ and $B$, respectively. For $x \in X \times \mathbb{T}$ and $y \in Y \times \mathbb{T}$, define
$$A_x := \mathrm{C}^*(\mathcal{C}(X \times \mathbb{T}), u \mathcal{C}_0((X\times \mathbb{T}) \setminus \{x\})),$$
and similarly, 
$$B_y := \mathrm{C}^*(\mathcal{C}(Y \times \mathbb{T}), v \mathcal{C}_0((Y \times \mathbb{T}) \setminus \{y\})).$$
These are generalizations of the subalgebras introduced by Putnam in \cite{Putnam:MinHomCantor}. By \cite[Section 3]{LinQPhil:KthoeryMinHoms} $A_x$ and $B_y$ can be written as inductive limits $A_x = \lim A_x^{(n)}$ and $B_y = \lim B_y^{(n)}$ where $A_x^{(n)}$ and $B_y^{(n)}$ are recursive subhomogeneous $\mathrm{C}^{*}$-algebras of topological dimensions $\dim(X \times \mathbb{T})$ and $\dim(Y \times \mathbb{T})$, respectively. Hence by Proposition~\ref{B} the recursive subhomogeneous algebras can be chosen to have $(\F, \eta)$-connected decompositions  with base spaces of dimension less than or equal to one. It follows from Corollary~\ref{dim1} that projections can be lifted along the recursive subhomogeneous decompositions. 

We have affine homeomorphisms $T(A_x) \cong T(A)$, $T(B_y) \cong T(B)$ and order isomorphisms $K_0(A_y) \cong K_0(A)$, $K_0(B_y) \cong K_0(B_y)$ \cite[Theorem 1.2 (2), (4)]{LinQPhil:KthoeryMinHoms} so the requirements for Corollary~\ref{dim1 TAI} are satisfied, hence with Corollary~\ref{TAI any UHF} we see that $A_x \otimes M_{\mathfrak{p}}$ and $B_y \otimes M_{\mathfrak{p}}$ are TAI for any supernatural number $\mathfrak{p}$. From \cite[Theorem 4.5]{StrWin:Z-stab_min_dyn} this implies that $A \otimes M_{\mathfrak{p}}$ and $B \otimes M_{\mathfrak{p}}$ are both TAI.

Since $A$ and $B$ satisfy the UCT and are $\mathcal{Z}$-stable \cite[Theorem B]{TomsWinter:PNAS} (also see \cite[Theorem 0.2]{TomsWinter:minhom}), as in the proof of Corollary~\ref{Ell blocks}, the result now follows from \cite[Corollary 11.9]{Lin:asu-class} 
\end{nproof}
\en

Since our main theorem does not require that projections separate tracial states (indeed, we assume this is not the case if there is more than one tracial state), this classification result is not covered by \cite[Theorem A]{TomsWinter:PNAS} (see also Theorem~0.1 of \cite{TomsWinter:minhom}). Despite the fact that all tracial states induce the same state on $K_0$, our result does not cover Connes' odd sphere examples \cite[Section 3]{Con:Thom}. In this case, the dimension of base spaces of the standard RSH decomposition of $A_y$ will not have dimension less than or equal to one, and so it is not clear that any decomposition can be found which satisfies the projection lifting requirement.  At present this is needed to guarantee that the solving of the linear system in Proposition~\ref{J} produces a projection, but it is possible that the lifting result is more than is necessary.  It seems promising that our current techniques can be adapted to cover more general $\mathrm{C}^{*}$-algebras, specifically more complicated $\mathrm{C}^{*}$-algebras of minimal dynamical systems including the Connes spheres.

In subsequent work we will pass from finitely many extremal traces to arbitrary trace spaces. This requires an additional set of UHF slicing tools as well as certain lifting results for maps between Cuntz semigroups (see for example \cite{CiuEll:InvSR1, RobSan:ClassofHoms, CiuEllSan:Cu1dim, Robert:NCCW}).

\bibliographystyle{amsplain}
\providecommand{\bysame}{\leavevmode\hbox to3em{\hrulefill}\thinspace}
\providecommand{\MR}{\relax\ifhmode\unskip\space\fi MR }
\providecommand{\MRhref}[2]{%
  \href{http://www.ams.org/mathscinet-getitem?mr=#1}{#2}
}
\providecommand{\href}[2]{#2}

\end{document}